\documentclass[reqno]{amsart}
\usepackage{setspace}
\usepackage{graphicx}
\usepackage{amsmath}
\usepackage{amsfonts}
\usepackage{amssymb}
\usepackage{amsthm}
\usepackage{mathtools}
\usepackage{subfigure}
\usepackage{enumerate}
\usepackage{enumitem}

\begin{document}
{\large
\centerline{\Huge\bf Characterization of the Three-Dimensional}

\medskip
\centerline{\Huge\bf Fivefold Translative Tiles}

\bigskip\medskip
\centerline{\large Mei Han, Kirati Sriamorn, Qi Yang and Chuanming Zong \footnote{All Mei Han, Kirati Sriamorn and Qi Yang are first authors and Chuanming Zong is the corresponding author.}
}

\bigskip\medskip
\centerline{\begin{minipage}{13cm}
{\bf Abstract.} This paper proves the following statement: {\it If a convex body can form a fivefold translative tiling in $\mathbb{E}^3$, it must be a parallelotope, a hexagonal prism, a rhombic dodecahedron, an elongated dodecahedron, a truncated octahedron, a cylinder over a particular octagon, or a cylinder over a particular decagon, where the octagon and the decagon are fivefold translative tiles in $\mathbb{E}^2$. Furthermore, it presents  an example of multiple tiles in $\mathbb{E}^3$ with multiplicity at most 10 which is neither a parallelohedron nor a cylinder.}
\end{minipage}}



\vspace{2.5cm}
\bigskip
\centerline{\LARGE\bf Contents}

\bigskip

\begin{flushright} 1 \end{flushright}
\vspace{-0.5cm}
1. Introduction

\begin{flushright} 3 \end{flushright}
\vspace{-0.5cm}
2. Basic known results

\begin{flushright} 5 \end{flushright}
\vspace{-0.5cm}
3. Definitions and lemmas

\begin{flushright} 7 \end{flushright}
\vspace{-0.5cm}
4. Some properties of zonotopes

\begin{flushright} 15 \end{flushright}
\vspace{-0.5cm}
5. Projections of fivefold translative tiles

\begin{flushright} 33 \end{flushright}
\vspace{-0.5cm}
6. Non-proper fivefold translative tilings

\begin{flushright} 44 \end{flushright}
\vspace{-0.5cm}
7. Proper fivefold translative tilings

\begin{flushright} 51 \end{flushright}
\vspace{-0.5cm}
8. Proof of Theorem 1.1

\begin{flushright} 51 \end{flushright}
\vspace{-0.5cm}
9. A non-trivial multiple lattice tile

\begin{flushright} 52 \end{flushright}
\vspace{-0.5cm}
Acknowledgements

\begin{flushright} 52 \end{flushright}
\vspace{-0.5cm}
References

\vspace{1cm}
\noindent
{\LARGE\bf 1. Introduction}

\bigskip\noindent
Let $K$ be a convex body with interior $int(K)$ and boundary $\partial (K)$, and let $X$ be a discrete set, both in $\mathbb{E}^n$. We call $K+X$ a {\it translative tiling} of $\mathbb{E}^n$ and call $K$ a {\it translative tile} if $K+X=\mathbb{E}^n$ and the translates $int(K)+{\bf x}_i$ are pairwise disjoint. In other words, if $K+X$ is both a packing in $\mathbb{E}^n$ and a covering of $\mathbb{E}^n$. In particular, we call $K+X$ a {\it lattice tiling} of $\mathbb{E}^n$ and call $K$ a {\it lattice tile} if $X$ is an $n$-dimensional lattice. Apparently, a translative tile must be a convex polytope. Usually, a lattice tile is called a {\it parallelohedron}.

In 1885, Fedorov \cite{fedo} characterized the two- and three-dimensional lattice tiles: {\it A two-dimensional lattice tile is either a parallelogram or a centrally symmetric hexagon; A three-dimensional lattice tile must be a parallelotope, a hexagonal prism, a rhombic dodecahedron, an elongated dodecahedron, or a truncated octahedron.} The situations in higher dimensions turn out to be very complicated. Through the works of Delone \cite{delo}, $\check{S}$togrin \cite{stog} and Engel \cite{enge}, we know that there are exact $52$ combinatorially different types of parallelohedra in $\mathbb{E}^4$. A computer classification for the five-dimensional Dirichlet-Voronoi cells of lattices was published by Dutour Sikiri$\acute{\rm c}$, Garber, Sch$\ddot{\rm u}$rmann and Waldmann \cite{dgsw} only in 2016. They are $110 244$ combinatorially different types.

Fedorov's discovery inspired D. Hilbert \cite{hilbert} asking the following question as a part of his 18th problem: {\it Whether polyhedra also exist which do not appear as fundamental regions of groups of motions, by means of which nevertheless by a suitable juxtaposition of congruent copies a complete filling up of all space is possible.} The answer to Hilbert's question is \lq\lq yes" even in $\mathbb{E}^2$ (see \cite{rao,zong20}). However, when $n\ge 3$, the structures of tiles and fundamental regions of groups of motions in $\mathbb{E}^n$ are complicated and still far from well understood (see \cite{Santos, Schulte}).

Let $\Lambda $ be an $n$-dimensional lattice. The {\it Dirichlet-Voronoi cell} of $\Lambda $ centered at the origin ${\bf o}$ is defined by
$$D=\left\{ {\bf x}: {\bf x}\in \mathbb{E}^n,\ \| {\bf x}, {\bf o}\| = \| {\bf x}, \Lambda \|\right\},$$
where $\| X, Y\|$ denotes the Euclidean distance between $X$ and $Y$. Clearly, $D+\Lambda $ is a lattice tiling and the Dirichlet-Voronoi cell $D$ is a parallelohedron. In 1908, Voronoi \cite{voro} made a conjecture that {\it every parallelohedron is a linear image of the Dirichlet-Voronoi cell of a suitable lattice.} In $\mathbb{E}^2$, $\mathbb{E}^3$ and $\mathbb{E}^4$, this conjecture was confirmed by Delone \cite{delo} in 1929. Recently, a proof for the five-dimensional Voronoi conjecture was announced by Garber and Magazinov \cite{Garber}. Combining the results of \cite{dgsw} and \cite{Garber}, all five-dimensional parallelohedra can be completely enumerated. In higher dimensions, both Voronoi's conjecture and the problem of parallelohedra classification are still open.

To characterize the translative tiles is another fascinating problem. First it was shown by Minkowski \cite{mink} in 1897 that {\it every translative tile must be centrally symmetric}. In 1954, Venkov \cite{venk} and Aleksandrov \cite{alek} proved that {\it all translative tiles are parallelohedra.} Later, a new proof for this beautiful result was independently discovered by McMullen \cite{mcmu2} (see also \cite{zong96}).

Let $X$ be a discrete multiset in $\mathbb{E}^n$ and let $k$ be a positive integer. We call $K+X$ a {\it $k$-fold translative tiling} of $\mathbb{E}^n$ and call $K$ a {\it $k$-fold translative tile} if every point ${\bf x}\in \mathbb{E}^n$ belongs to at least $k$ translates of $K$ in $K+X$ and every point ${\bf x}\in \mathbb{E}^n$ belongs to at most $k$ translates of $int(K)$ in $int(K)+X$. In other words, if $K+X$ is both a $k$-fold packing in $\mathbb{E}^n$ and a $k$-fold covering of $\mathbb{E}^n$. In particular, we call $K+\Lambda$ a {\it $k$-fold lattice tiling} of $\mathbb{E}^n$ and call $K$ a {\it $k$-fold lattice tile} if $\Lambda $ is an $n$-dimensional lattice. Apparently, a $k$-fold translative tile must be a convex polytope (see \cite{kolo-m}).

Multiple tilings were first investigated by Furtw\"angler \cite{furt} in 1936 as a generalization of Minkowski's conjecture on cube tilings. Let $C$ denote the $n$-dimensional unit cube. Furtw\"angler made a conjecture that {\it every $k$-fold lattice tiling $C+\Lambda$ has twin cubes. In other words, every multiple lattice tiling $C+\Lambda$ has two cubes which share a whole facet.} In the same paper, he proved the two- and three-dimensional cases. Unfortunately, when $n\ge 4$, this beautiful conjecture was disproved by Haj\'os \cite{hajo} in 1941. In 1979, Robinson \cite{robi} determined all the integer pairs $\{ n,k\}$ for which Furtw\"angler's conjecture is false. We refer to Zong \cite{zong05,zong06} for detailed accounts on this fascinating problem.

\medskip
Recent years, multiple tilings has been studied by Bolle \cite{boll}, Gravin, Robins and Shiryaev \cite{grs}, Gravin, Kolountzakis, Robins and Shiryaev \cite{gkrs}, Grepstad and Lev \cite{grepstad}, Kolountzakis \cite{kolo}, Lev and Liu \cite{lev-liu}, Liu \cite{liu}, Yang and Zong \cite{yz1,yz2}, Zong \cite{zong20, zong-x} and others. Clearly, one of the most important and natural problems in multiple tilings is to classify or characterize the $n$-dimensional $k$-fold translative tiles and the $n$-dimensional $k$-fold lattice tiles (see open questions 1-4 at the end of Gravin, Robins and Shiryaev \cite{grs}). In the plane, it was proved by Yang and Zong \cite{yz1,yz2,zong20,zong-x} that, {\it besides parallelograms and centrally symmetric hexagons, there is no other two-, three- or fourfold translative tile. However, there are three classes of other fivefold translative tiles and three classes of other sixfold lattice tiles.} In the three-dimensional Euclidean space $\mathbb{E}^3$, it was shown by Han, Sriamorn, Yang and Zong \cite{hsyz} that, {\it besides parallelotopes, hexagonal prisms, rhombic dodecahedra, elongated dodecahedra and truncated octahedra, there is no other two-, three- or fourfold translative tile.} This paper proves the following result.

\medskip\noindent
{\bf Theorem 1.1.} {\it In $\mathbb{E}^3$, a convex body can form a fivefold translative tiling if and only if it is a parallelotope, a hexagonal prism, a rhombic dodecahedron, an elongated dodecahedron, a truncated octahedron, a cylinder over a particular octagon, or a cylinder over a particular decagon, where the octagon and the decagon are fivefold translative tiles in $\mathbb{E}^2$.}

\medskip
Furthermore, based on the results of Gravin, Robins and Shiryaev \cite{grs} and Lev and Liu \cite{lev-liu}, we discover an example of multiple tiles in $\mathbb{E}^3$ with multiplicity at most 10 which is neither a parallelohedron nor a cylinder.

\medskip\noindent
{\bf Remark 1.1.} At the end of \cite{grs}, Gravin, Robins and Shiryaev listed \lq\lq Classify the combinatorial types of all polytopes which $k$-tile $\mathbb{E}^d$ by translations." as a fascinating open question. Clearly, Theorem 1.1 resolves the $d=3$ and $k=5$ case.

\medskip\noindent
{\bf Remark 1.2.} By the results in \cite{yz2,zong-x}, we have known that all fivefold translative tiles in $\mathbb{E}^2$ are fivefold lattice tiles. Together with Theorem 1.1, we can easily get that all fivefold translative tiles in $\mathbb{E}^3$ are also fivefold lattice tiles.

\medskip\noindent
{\bf Remark 1.3.} This paper is a continuation of Han, Sriamorn, Yang and Zong \cite{hsyz}. For the completeness of this paper, we repeat its introduction and some basic lemmas.

\vspace{1cm}
\noindent
{\LARGE\bf 2. Basic known results}

\bigskip\noindent
In this section, we recall some known concepts and results which will be useful for this paper.

\medskip
In 1885, E. S. Fedorov studied the two- and three-dimensional lattice tiles. He proved the following result.

\smallskip\noindent
{\bf Lemma 2.1 (Fedorov \cite{fedo}).} {\it A two-dimensional lattice tile is either a parallelogram or a centrally symmetric hexagon; A three-dimensional lattice tile must be a parallelotope, a hexagonal prism, a rhombic dodecahedron, an elongated dodecahedron, or a truncated octahedron.}

\smallskip
Tilings in higher dimensions have been studied by Minkowski \cite{mink}, Voronoi \cite{voro}, Delone \cite{delo}, Venkov \cite{venk}, Alexsandrov \cite{alek}, McMullen \cite{mcmu2} and many others. This paper needs the following concepts and results thereof.

\smallskip
\noindent
{\bf Definition 2.1.} Let $P$ denote an $n$-dimensional centrally symmetric convex polytope with centrally symmetric facets and let $E$ denote an $(n-2)$-dimensional face of $P$. We call the collection of all those facets of $P$ which contain a translate of $E$ as a subface a {\it belt} of $P$.

\medskip\noindent
{\bf Lemma 2.2 (Venkov \cite{venk} and McMullen \cite{mcmu2}).} {\it  Let $K$ be an $n$-dimensional convex body. The following three statements are equivalent to each other:}
\begin{enumerate}
\item {\it $K$ is a translative tile;}
\item {\it $K$ is a centrally symmetric polytope with centrally symmetric facets such that each belt contains four or six facets;}
\item {\it $K$ is a parallelohedron.}
\end{enumerate}

\medskip\noindent
{\bf Definition 2.2.} Let $P$ be an $n$-dimensional convex polytope. We call it a {\it zonotope} if it is a Minkowski sum of finite number of segments.
In other words,
$$P=S_0+S_1+\ldots +S_w,$$
where $w$ is an integer and $S_i$ are segments in $\mathbb{E}^n$. Without loss of generality, assume that $\mathbb{S}_P=\{S_0, S_1, S_2,.., S_w\}$ is a set of pairwise linear independent segments in $\mathbb{E}^n$ and we call $\mathbb{S}_P$ a {\it generator set} of $P$.

\medskip
Zonotopes are a class of important subjects in convex geometry and have been studied by many authors. For a survey on this topic, we refer to \cite{mcmu1} and \cite{Schneider}. In 2012, N. Gravin, S. Robins and D. Shiryaev studied multiple tilings in general dimensions and discovered the following result.

\medskip\noindent
{\bf Lemma 2.3 (Gravin, Robins and Shiryaev \cite{grs}).} {\it An $n$-dimensional $k$-fold translative tile is a centrally symmetric polytope with centrally symmetric facets. In particular, in $\mathbb{E}^3$, every $k$-fold translative tile is a zonotope.}

\medskip
In 2022, Han, Sriamorn, Yang and Zong \cite{hsyz} introduced dihedral adjacent wheel into the multiple tiling study in $\mathbb{E}^3$, which could be regarded as a generalization of adjacent wheel firstly introduced into the multiple tiling in the plane by Yand and Zong \cite{yz2} in 2021.

\medskip
Let $P$ be a three-dimensional zonotope centered at {\bf o}. Assume that $P+X$ is a $k$-fold translative tiling of $\mathbb{E}^3$, where X=$\{{\bf x}_1, {\bf x}_2, ...\}$ is a discrete multiset with ${\bf x}_1={\bf o}$. Let $E$ be an edge of $P$, and let $B(E)$ be the belt determined by $E$, which consists of $2m$ facets $F_1$, $F_2$, .., $F_{2m}$ enumerated in a circular order, then $B(E)=\{F_1, F_2,.., F_{2m}\}$. Apparently, $F_{i+m}=-F_i$ for $1\le i\le m$. Let $E_1$, $E_2$, ..., $E_{2m}$ be the translates of $E$ such that $E_i, E_{i+1}\subset F_i$ and $E_1=E$, $E_{2m+1}=E_1$. Moreover, let $A(E)=\{E_1, E_2,..,E_{2m}\}$ and let $\#B(E)$ denote the number of facets in $B(E)$.

Assume that $W$ is a bounded set in $\mathbb{E}^3$. As usual, let $relint(W)$ and $\overline{W}$ denote its relative interior and its closure, respectively. And we define $C(E)=\overline{\partial(P)\setminus B(E)}$ {\it the double caps} of $P$ in direction $E$.

Assume that ${\bf p}\in E\setminus\{C(E)+X\}$. Let $X^{\bf p}$ denote the subset of $X$ consisting of all points ${\bf x}_i$ such that $${\bf p}\in\partial(P)+{\bf x}_i.$$
Since $P+X$ is a multiple tiling, the set $X^{\bf p}$ can be divided into disjoint subsets $X_1^{\bf p}$, $X_2^{\bf p}$,..., $X_t^{\bf p}$ such that the translates in $P+X_j^{\bf p}$ can be enumerated as $P+{\bf x}_1^j$, $P+{\bf x}_2^j$,..., $P+{\bf x}_{s_j}^j$ satisfying the following conditions:
\begin{enumerate}
\item ${\bf p}\in\partial(P)+{\bf x}_i^j$ holds for all $i=1,2,..., s_j$;
\item Let $\alpha_i^j$ denote the dihedral angle of $P+{\bf x}_i^j$ at the edge which contains {\bf p}. If ${\bf p}\in relint(F_l)+{\bf x}_i^j$ holds for some $l$, then $\alpha_i^j=\pi$. Assume that the dihedral angle is bounded by two half-planes $N_{i,1}^j$ and $N_{i,2}^j$ such that $$N_{i,2}^j=N_{i+1,1}^j$$ holds for all $i=1,2,..., s_j$, where $N_{s_j+1,1}^j=N_{1,1}^j$.
\end{enumerate}
For convenience, we call such a sequence $P+{\bf x}_1^j$, $P+{\bf x}_2^j$,..., $P+{\bf x}_{s_j}^j$ a {\it dihedral adjacent wheel} at {\bf p} and we call such a sequence $\alpha_1^j$, $\alpha_2^j$,..., $\alpha_{s_j}^j$ {\it the structure of the dihedral adjacent wheel} at {\bf p}. In other words, if {\bf p} belongs to the relative interior of an edge of a tile then we follow this tile around, moving from tile to tile, until it closes up again. It is easy to see that $$\sum_{i=1}^{s_j}\alpha_i^j=2\omega_j\cdot\pi$$
hold for positive integers $\omega_j$. Then we define $$\varpi({\bf p})=\sum_{j=1}^t w_j=\frac{1}{2\pi}\sum_{j=1}^t\sum_{i=1}^{s_j}\alpha_i^j$$
and $$\varphi(\bf p)=\#\{{\bf x}_i:{\bf x}_i\in{\it X}, {\bf p}\in{\it int(P)}+{\bf x}_i\}.$$ In other words, $\varpi({\bf p})$ is the tiling multiplicity produced by the boundary and $\varphi({\bf p})$ is the tiling multiplicity produced by the interior. For every ${\bf p}\in E\setminus\{C(E)+X\}$, we define $$k({\bf p})=\varpi({\bf p})+\varphi({\bf p}).$$

Assume that the tiling multiplicity $k=5$, we naturally have $k({\bf p})=\varpi({\bf p})+\varphi({\bf p})=5$. For the fivefold translative tiling $P+X$ of $\mathbb{E}^3$, we have the following lemmas useful in this paper.

\medskip\noindent
{\bf Lemma 2.4 (Han, Sriamorn, Yang and Zong \cite{hsyz}).} {\it Every belt of a fivefold translative tile in $\mathbb{E}^3$ has at most ten facets. In other words, every belt of a fivefold translative tile in $\mathbb{E}^3$ has four, six, eight or ten facets.}

\medskip\noindent
{\bf Lemma 2.5 (Yang and Zong \cite{yz2}).} {\it A convex domain can form a fivefold translative tiling of the Euclidean plane if and only if it is a parallelogram, a centrally symmetric hexagon, under a suitable affine linear transformation, a centrally symmetric octagon with vertices ${\bf v}_1=(\frac{3}{2}-\frac{5\alpha}{4},-2)$, ${\bf v}_2=(-\frac{1}{2}-\frac{5\alpha}{4},-2)$, ${\bf v}_3=(\frac{\alpha}{4}-\frac{3}{2},0)$, ${\bf v}_4=(\frac{\alpha}{4}-\frac{3}{2},1)$, ${\bf v}_5=-{\bf v}_1$, ${\bf v}_6=-{\bf v}_2$, ${\bf v}_7=-{\bf v}_3$, and ${\bf v}_8=-{\bf v}_4$, where $0< \alpha <\frac{2}{3}$, or with vertices ${\bf v}_1=(2-\beta, -3)$, ${\bf v}_2=(-\beta, -3)$, ${\bf v}_3=(-2, -1)$, ${\bf v}_4=(-2, 1)$, ${\bf v}_5=-{\bf v}_1$, ${\bf v}_6=-{\bf v}_2$, ${\bf v}_7=-{\bf v}_3$, and ${\bf v}_8=-{\bf v}_4$, where $0< \beta \le 1$, or a centrally symmetric decagon with ${\bf u}_1=(0,1)$, ${\bf u}_2=(1,1)$, ${\bf u}_3=(\frac{3}{2}, \frac{1}{2})$, ${\bf u}_4=(\frac{3}{2}, 0)$, ${\bf u}_5=(1, -\frac{1}{2})$, ${\bf u}_6=-{\bf u}_1$, ${\bf u}_7=-{\bf u}_2$, ${\bf u}_8=-{\bf u}_3$, ${\bf u}_9=-{\bf u}_4$, and ${\bf u}_{10}=-{\bf u}_5$ as the middle points of its edges.}

\vspace{1cm}
\noindent
{\LARGE\bf 3. Definitions and lemmas}

\bigskip\noindent
In Section 2, we have introduced several definitions and notations which have been used in \cite{hsyz} and still useful in this paper. In order to simplify the following discussion, we will introduce some new definitions and notations in this section. Firstly, we recall a definition in \cite{hsyz}.

\medskip\noindent
{\bf Definition 3.1 (\bf Han, Sriamorn, Yang and Zong \cite{hsyz}).} Let $F$ be a centrally symmetric polygon with an edge $E$ and let {\bf g} be a vector such that $E+{\bf g}$ is also an edge of $F$. Then we call {\bf g} {\it the translation vector of E in F}. For ${\bf v}\in relint(E)$, we call ${\bf v}+{\bf g}$ {\it the corresponding point of {\bf v} in $F$}.

\medskip
In Section 2, we have assumed that the three-dimensional fivefold translative tile $P$ is a zonotpe centered at {\bf o}, and $P$ has an edge $E$ and a belt $B(E)=\{F_1,F_2,...,F_{2m}\}$ determined by $E$. Without loss of generality, let $\mathbb{S}_P=\{S_0, S_1,..., S_w\}$ be a generator set of $P$, where $\mathbb{S}_P$ is a set of pairwise linear independent segments centered at {\bf o} in $\mathbb{E}^3$. For convenience, let $E$ be a translate of $S_0$. Then we introduce the following definitions.

\medskip\noindent
{\bf Definition 3.2.}  Let $$F_i^0=S_0+S_1^i+\cdots+S_{w_i}^i,$$ where $\{S_0, S_1^i, \cdots,S_{w_i}^i\}\subset\mathbb{S}_P$ for $1\le i\le m$. We call $F_i^0$ an {\it original facet} of $P$ corresponding to $F_i$ if there is vector ${\bf g}_{F_i^0}$ such that $F_i=F_i^0-{\bf g}_{F_i^0}$ where $F_i\in B(E)$. Apparently, $F_{i+m}=F_i+2{\bf g}_{F_i^0}$. And similar to Definition 3.1, we call ${\bf g}_i=2{\bf g}_{F_i^0}$ the {\it translation vector} of $F_i$ in $P$. For convenience, write $\{S_0, S_1^i, \cdots,S_{w_i}^i\}$ as $\mathbb{S}_{F_i^0}$, the generator set of $F_i^0$. Moreover, we assume $R_i=S_1^i+S_2^i+\cdots+S_{w_i}^i$ and $\mathbb{S}_{R_i}=\{S_1^i,S_2^i,...,S_{w_i}^i\}$.

\begin{figure}[h!]
\centering
\includegraphics[scale=0.65]{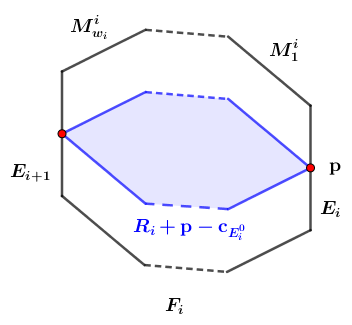}
\caption{The facet $F_i$ of $P$}\label{Fi}
\end{figure}

Let $E_i^0$ and $E_{i+1}^0$ be the two edges of $F_i^0$ which are translates of $E$, then it is easy to see that the midpoint ${\bf c}_{E_i^0}$ of $E_{i}^0$ and the midpoint ${\bf c}_{E_{i+1}^0}$ of $E_{i+1}^0$ are two vertices of $R_i$, and ${\bf c}_{E_{i+1}^0}=-{\bf c}_{E_i^0}$ for $1\le i\le m$. Naturally, $F_{i+m}^0=F_i^0$ and $R_{i+m}=R_i$.

For convenience, we assume that $E_i$, $M_1^i$, $M_2^i$,..., $M_{w_i}^i$, $E_{i+1}$ in counterclockwise order are edges of $F_i$ and $M_j^i$ be a translate of $S_j^i$ for $1\le i\le m$ and $w_i\in\mathbb{N}^*$, as shown in FIG \ref{Fi}. And it is easy to see that $E_{i+m}$, $M_1^{i+m}$,.., $M_{w_i}^{i+m}$, $E_{i+m+1}$ are edges of $F_{i+m}$ satisfying that $M_k^{i+m}=M_{w_i-k+1}^i+{\bf g}_i$ and $M_k^{i+m}$ is the translate of $S_{w_i-k+1}^i$. Let $\widetilde{\bf g}_i$ be the translation vector of $E_i$ in $F_i$ for $1\le i\le m$, that is $E_{i+1}=E_i+\widetilde{\bf g}_i$, clearly, $\widetilde{\bf g}_{i+m}=-\widetilde{\bf g}_i$. We assume that $\vec{M}_j^i$ is the vector form of $M_j^i$ satisfying that $$\widetilde{\bf g}_i=\sum_{j=1}^{w_i}\vec{M}_j^i.$$ And we have that $${\bf g}_i=\widetilde{\bf g}_{i+1}+\widetilde{\bf g}_{i+2}+\cdots+\widetilde{\bf g}_{i+m-1}$$ and $\vec{M}_k^{i+m}=-\vec{M}_{w_i-k+1}^i$.

\medskip\noindent
{\bf Definition 3.3.} If for any $E'\in A(E)+X$ and any $P'\in P+X$, we have $E'\subset P'$ or $relint(E')\cap P'=\emptyset$, then we call $P+X$ a {\it proper translative tiling relate to $E$}.

\medskip\noindent
{\bf Definition 3.4.} If there exists some $E'\in A(E)+X$ and some $P'\in P+X$ satisfying $E'\not\subset P'$ but $relint(E')\cap P'\not=\emptyset$, we call $P+X$ a {\it non-proper translative tiling relate to $E$}.

\medskip\noindent
{\bf Definition 3.5 (\bf Han, Sriamorn, Yang and Zong \cite{hsyz}).} Assume that $E'\in A(E)+X$, we call a point ${\bf p}\in E'$ {\it a proper point}, if {\bf p} satisfies the following conditions:
\begin{enumerate}
  \item ${\bf p}\notin C(E)+X$.
  \item If ${\bf p}\in relint(F_i)+{\bf x}$ holds for some integer $i$ and some point ${\bf x}\in X$, then there are two points ${\bf p}^*\in(E_i+{\bf x})\setminus\{C(E)+X\}$ and ${\bf p}^\bullet={\bf p}^*+\widetilde{\bf g}_i\in(E_{i+1}+{\bf x})\setminus\{C(E)+X\}$ such that $${\bf p}\in{\bf p}^*+R_i-{\bf c}_{E_i^0}={\bf p}^\bullet+R_i-{\bf c}_{E_{i+1}^0}.$$
  \item If ${\bf p}\in E^*\cap F\subset F^*\cap F$, where $F\in B(E)+X, F^*\in B(E)+X$, and $E^*$ is an edge of $F^*$ which is parallel to $E'$, then the corresponding point of {\bf p} in $F^*\cap F$ is not in $ C(E)+X$.
\end{enumerate}

\medskip\noindent
{\bf Definition 3.6 (\bf Han, Sriamorn, Yang and Zong \cite{hsyz}).} Suppose that ${\bf x}\in X$, ${\bf p}\in B(E)+{\bf x}$ and ${\bf p}\notin C(E)+X$. If {\bf p} belongs to the relative interior of an edge of $P+{\bf x}$, we call $P+{\bf x}$ an {\it E-type translate} at {\bf p}; If {\bf p} belongs to the relative interior of a facet of $P+{\bf x}$, we call $P+{\bf x}$ an {\it F-type translate} at {\bf p}.

Moreover, if ${\bf p}\in int(P)+{\bf x}'$ for ${\bf x}'\in X$, then we call $P+{\bf x}'$ an {\it $I$-type translate} at {\bf p}.

\medskip
Use the definitions and notations above and we introduce several lemmas important in the following discussion.

\medskip\noindent
{\bf Lemma 3.1 (Han, Sriamorn, Yang and Zong \cite{hsyz}).} {\it Assume that ${\bf p}\in E_i\setminus\{C(E)+X\}$, where $i\in\{1,2,...,2m\}$. There are at least $\lceil(m-3)/2\rceil$ different translates $P+{\bf x}_j$ satisfying both $${\bf p}\in B(E)+{\bf x}_j$$ and $$({\bf p}+R_i-{\bf c}_{E_i^0})\setminus\{{\bf p}\}\subset int(P)+{\bf x}_j.$$}

\medskip\noindent
{\bf Lemma 3.2 (Han, Sriamorn, Yang and Zong \cite{hsyz}).} {\it For every ${\bf p}\in E\setminus\{C(E)+X\}$, we have $$\varpi({\bf p})=\kappa({\bf p})\cdot\frac{m-1}{2}+\ell({\bf p})\cdot\frac{1}{2},$$ where $\kappa({\bf p})$ is the covering multiplicity of {\bf p} contributed by the dihedral angles at {\bf p} which are not equal to $\pi$ and $\ell({\bf p})$ is the number of the facets in $B(E)+X$ which take {\bf p} as a relative interior point.}

\medskip\noindent
{\bf Lemma 3.3 (Han, Sriamorn, Yang and Zong \cite{hsyz}).} {\it Assume that $E'\in A(E)+X$. If {\bf p} is a proper point in $E'$, then $$\varphi({\bf p})\ge\left\lceil\frac{m-3}{2}\right\rceil.$$}

\medskip\noindent
{\bf Remark 3.1.} It is easy to see that for any facet $F$ of $P$, there is always an original facet $F^0$ of $P$ corresponding to itself.

\medskip\noindent
{\bf Remark 3.2.} By the conditions in Definition 3.5 and the analysis in \cite{hsyz}, we have known that each $E'\in A(E)$ has infinite numbers of proper points.

\medskip\noindent
{\bf Remark 3.3.} Lemmas 3.1-3.3 also played an important role in \cite{hsyz} and here we only adjust some expressions by the notations in this section.

\vspace{1cm}\noindent
{\LARGE\bf 4. Some properties of zonotopes}

\bigskip\noindent
Let $P$ be a three-dimensional zonotope and $E$ be an edge of $P$, as introduced in Section 3, $\mathbb{S}_P=\{S_0, S_1,....,S_w\}$ is the generator set of $P$ with $\#B(E)=2m$. In this section, we mainly study some properties of zonotopes which play crucial roles in the following discussion. By Definition 3.2 and the analysis of the representation of one zonotope by its generator segments in \cite{mcmu1}, we can easily get Lemma 4.1.

\medskip\noindent
{\bf Lemma 4.1.} {\it Assume that $\mathbb{T}=\{S_{w1}, S_{w2},... S_{wt}\}\subset\mathbb{S}_P$ is a maximal coplanar segments set, that is $S_{w1}$, $S_{w2}$,... $S_{wt}$ are coplanar but for any $S_{w,{t+1}}\in\mathbb{S}_P\setminus\mathbb{T}$, $S_{w1}, S_{w2},... S_{wt}, S_{w,{t+1}}$ are not coplanar, then $$F^0=\sum_{j=1}^t S_{wj}$$ must be an original facet in $P$ for any positive integer $t\ge2$.}

\medskip
Moreover, we have the following results about zonotopes.

\medskip\noindent
{\bf Lemma 4.2.} {\it For $F_i\in B(E)$, $1\le i\le m$, let $F_i^0=S_0+S_1^i+S_2^i+\cdots+S_{w_i}^i$ be the original facet of $P$ corresponding to $F_i$ and $\mathbb{S}_{F_i^0}=\{S_0, S_1^i,..., S_{w_i}^i\}$, then $$\mathbb{S}_P=\bigcup_{i=1}^m\mathbb{S}_{F_i^0}$$ and $$\mathbb{S}_{F_j^0}\bigcap\mathbb{S}_{F_k^0}=\{S_0\}$$ for $j\not=k$.}

\medskip\noindent
{\bf Proof.} It is easy to see that $\bigcup\limits_{i=1}^m\mathbb{S}_{F_i^0}\subseteq \mathbb{S}_P$. Let us show $\mathbb{S}_P=\bigcup\limits_{i=1}^m\mathbb{S}_{F_i^0}$. For the contrary, suppose there is a segment $S^*\in\mathbb{S}_P$ but $S^*\notin\mathbb{S}_{F_i^0}$ for any $i=1,2,...,m$. Since $S^*$ and $S_0$ are coplanar, $S^*+S_0\subseteq F_k^0$ for some $k$, that is $S^*\in\mathbb{S}_{F_k^0}$, contradictorily. Thus we have $$\mathbb{S}_P=\bigcup_{i=1}^m\mathbb{S}_{F_i^0}.$$ Similarly, by Lemma 4.1, we have $$\mathbb{S}_{F_j^0}\bigcap\mathbb{S}_{F_k^0}=\{S_0\}$$ for $j\not=k$.

Thus the proof of this lemma is finished. \hfill{$\Box$}

\medskip
 In fact, let $\mathbb{S}_{B(E)}=\bigcup\limits_{i=1}^m\mathbb{S}_{F_i^0}$ be a generator set of $B(E)$, then we have $$\mathbb{S}_{B(E)}=\mathbb{S}_P.$$ Without loss of generality, let $$\widetilde{\bf g}_i=\vec{S}_1^i+\vec{S}_2^i+\cdots+\vec{S}_{w_i}^i,$$ where $S_j^i\in\mathbb{S}$ and $\vec{S}_j^i=\vec{M}_j^i$ for $1\le i\le m$ and let $\vec{S}_0$ be the vector form of $S_0$ whose endpoint is the joint point of $E_i$ and $M_1^i$ shown in Definition 3.2. Notice $\mathbb{S}_{R_i}=\mathbb{S}_{F_i}^0\setminus\{S_0\}$.

\medskip
Before the following lemma, we have to generalize the definition of the belt in Definition 2.1.

Let $Q$ be a three-dimensional centrally symmetric convex polytope centered at {\bf o} and let $\widetilde{E}$ be one edge of $Q$. Assume that $\{\widetilde{F}_1,\widetilde{F}_2,...,\widetilde{F}_{2l}\}$ is a family of facets of $Q$ for some positive integer $l$. If $\widetilde{E}_i$, $\widetilde{E}_{i+1}\subset\partial(\widetilde{F}_i)$ where $\widetilde{E}_i$ is parallel to $\widetilde{E}$ for $1\le i\le 2l$ and $\widetilde{E}_{2l+1}=\widetilde{E}_1$, then we also call $\{\widetilde{F}_1,\widetilde{F}_2,...,\widetilde{F}_{2l}\}$ a {\it belt} of $Q$ determined by $\widetilde{E}$. To avoid confusion thereinafter, we use $B(Q,\widetilde{E})$ to denote the belt of $Q$ determined by $\widetilde{E}$. Particularly, $B(Q,\widetilde{E})=B(\widetilde{E})$ if $Q$ is a three-dimensional zonotope.

\medskip\noindent
{\bf Lemma 4.3.} {\it Assume that ${\bf g}=\sum\limits_{i=2}^{k-1}\widetilde{\bf g}_i+\sum\limits_{j=1}^{w_k}\lambda_j^k\vec{S}_j^k+\sum\limits_{l=1}^{w_1}\lambda_l^1\vec{S}_l^1+\lambda_0\vec{S}_0$ for $2\le k\le m$, $0\le\lambda_0<1$, $0\le\lambda_l^1<1$ and $0\le\lambda_j^k\le1$ where $1\le l\le w_1$ and $1\le j\le w_k$. Let $P^*$ be a zonotope with $\mathbb{S}_{P^*}=\{(1-\lambda_0)S_0, (1-\lambda_1^1)S_1^1,...,(1-\lambda_{w_1}^1)S_{w_1}^1, (1-\lambda_1^k)S_1^k,...,(1-\lambda_{w_k}^k)S_{w_k}^k\}\bigcup(\bigcup\limits_{i=k+1}^{m}\mathbb{S}_{R_i})$, then we have $$P^*+\frac{1}{2}{\bf g}\subset P\cap(P+{\bf g}),$$ where $P^*=\sum\limits_{S\in\mathbb{S}_{P^*}}S$. Further, let $E^*$ be an edge of $P^*$ which is also a translate of $(1-\lambda_0)E$ and let $B(P^*,E^*)$ be a belt of $P^*$, then $P\cap(P+{\bf g})$ also has a belt $B(P\cap(P+{\bf g}),E^*)$ such that $$\bigcup_{F\in B(P^*,E^*)+\frac{1}{2}{\bf g}}F\subseteq\bigcup_{F\in B(P\cap(P+{\bf g}),E^*)}F.$$ In particular, for $k=m$ and all $\lambda_j^k=1$, we have $$P^*+\frac{1}{2}{\bf g}= P\cap(P+{\bf g})\subseteq F_{m+1}.$$}

\medskip\noindent
{\bf Proof.} By Lemma 4.2, it suffices to show that the boundary of $P\cap(P+{\bf g})$ contains the belt $B(P^*+\frac{1}{2}{\bf g},(1-\lambda_0)E)$. For convenience, let $$\overline{\mathbb{S}}_{R_1}=\{\lambda_1^1S_1^1,(1-\lambda_1^1)S_1^1,...,\lambda_{w_1}^1S_{w_1}^1,(1-\lambda_{w_1}^1)S_{w_1}^1\},$$ $$\overline{\mathbb{S}}_{R_k}=\{\lambda_1^kS_1^k,(1-\lambda_1^k)S_1^k,...,\lambda_{w_k}^kS_{w_k}^k,(1-\lambda_{w_k}^k)S_{w_k}^k\},$$ and $$\overline{\mathbb{S}}_P=\{\lambda_0S_0,(1-\lambda_0)S_0\}\bigcup\overline{\mathbb{S}}_{R_1}\bigcup\overline{\mathbb{S}}_{R_k} \bigcup(\bigcup\limits_{i=k+1}^{m}\mathbb{S}_{R_i}).$$ For convenience, we firstly consider $2\le k<m$.

Since $(1-\lambda_0)S_0\in\mathbb{S}_{P^*}$ which is a translate of $(1-\lambda_0)E$, $P^*$ has an edge $E^*$ that is a translate of $(1-\lambda_0)E$ and naturally has a belt $B(P^*, E^*)$. Moreover, $\mathbb{S}_{B(P^*,E^*)}=\mathbb{S}_{P^*}\subset\overline{\mathbb{S}}_P$. Without loss of generality, we assume that $F^0$ is an original facet of $P^*$ corresponding to a facet in $B(P^*,E^*)$.
\begin{enumerate}
  \item If $\mathbb{S}_{F^0}=\{(1-\lambda_0)S_0,(1-\lambda_1^1)S_1^1,...,(1-\lambda_{w_1}^1)S_{w_1}^1\}$, then there is a vector $${\bf g}_{F^0}=\frac{1}{2}\left[\sum_{j=1}^{w_k}(1-\lambda_j^k)\vec{S}_j^k+\sum_{i=k+1}^m\widetilde{\bf g}_i\right]\in\sum_{T\in \mathbb{S}_{P^*}\setminus\mathbb{S}_{F^0}}T$$ satisfying that $F^0\pm{\bf g}_{F^0}$ are the facets of $P^*$, thus $F^0\pm{\bf g}_{F^0}+\frac{1}{2}{\bf g}$ are the facets of $P^*+\frac{1}{2}{\bf g}$ and $F^0\subseteq F_1^0$.

      Since $$\pm{\bf g}_{F^0}\pm\frac{1}{2}{\bf g}\in\sum_{T\in\mathbb{S}_{P^*}\setminus\mathbb{S}_{F^0}}T+\sum_{T\in\overline{\mathbb{S}}_P\setminus\mathbb{S}_{P^*}}T=\sum_{T\in\overline{\mathbb{S}}_P\setminus\mathbb{S}_{F^0}}T,$$ we have $$F_0\pm{\bf g}_{F^0}\pm\frac{1}{2}{\bf g}\subset\sum_{T\in\mathbb{S}_{F^0}}T+\sum_{T\in\overline{\mathbb{S}}_P\setminus\mathbb{S}_{F^0}}T=\sum_{S\in\mathbb{S}_P}S=P.$$ Note that $${\bf g}_{F^0}+\frac{1}{2}{\bf g}=\frac{1}{2}\sum_{i=2}^m\widetilde{\bf g}_i+\frac{1}{2}\sum_{l=1}^{w_1}\lambda_l^1\vec{S}_l^1+\frac{1}{2}\lambda_0\vec{S}_0=\frac{1}{2}{\bf g}_1+\frac{1}{2}\sum_{l=1}^{w_1}\lambda_l^1\vec{S}_l^1+\frac{1}{2}\lambda_0\vec{S}_0,$$ we have $$F^0+{\bf g}_{F^0}+\frac{1}{2}{\bf g}\subseteq F_1^0+\frac{1}{2}{\bf g}_1=F_{1+m}\subset\partial(P)$$ and $$F^0+{\bf g}_{F^0}+\frac{1}{2}{\bf g}=(F^0+{\bf g}_{F^0}-\frac{1}{2}{\bf g})+{\bf g}\subset P+{\bf g}.$$ Thus $F^0+{\bf g}_{F^0}+\frac{1}{2}{\bf g}\subset\partial[P\cap(P+{\bf g})]$.

      Similarly, since $-{\bf g}_{F^0}+\frac{1}{2}{\bf g}=-({\bf g}_{F^0}+\frac{1}{2}{\bf g})+{\bf g}$, we have $$F^0-{\bf g}_{F^0}+\frac{1}{2}{\bf g}\subseteq F_1^0-\frac{1}{2}{\bf g}_1+{\bf g}=F_1+{\bf g}\subset\partial(P)+{\bf g}$$ and $$F^0-{\bf g}_{F^0}+\frac{1}{2}{\bf g}\subset P,$$ Thus $F^0-{\bf g}_{F^0}+\frac{1}{2}{\bf g}\subset\partial[P\cap(P+{\bf g})]$.
  \item If $\mathbb{S}_{F^0}=\{(1-\lambda_0)S_0, (1-\lambda_1^k)S_1^k,...,(1-\lambda_{w_k}^k)S_{w_k}^k\}$, then there is a vector $${\bf g}_{F^0}=\frac{1}{2}\sum_{i=k+1}^{m}\widetilde{\bf g}_i-\frac{1}{2}\sum_{l=1}^{w_1}(1-\lambda_l^1)\vec{S}_l^1\in\sum_{T\in \mathbb{S}_{P^*}\setminus\mathbb{S}_{F^0}}T$$ satisfying that $F^0\pm{\bf g}_{F^0}$ are the facets of $P^*$, thus $F^0\pm{\bf g}_{F^0}+\frac{1}{2}{\bf g}$ are the facets of $P^*+\frac{1}{2}{\bf g}$ and $F^0\subseteq F_k^0$.

      Similar to (1), we also have $F_0\pm{\bf g}_{F^0}\pm\frac{1}{2}{\bf g}\subset P$. Note that $${\bf g}_{F^0}-\frac{1}{2}{\bf g}=\frac{1}{2}{\bf g}_k-\frac{1}{2}\left(\sum_{j=1}^{w_k}\lambda_j^k\vec{S}_j^k+\lambda_0\vec{S}_0\right),$$ we have $$F^0+{\bf g}_{F^0}+\frac{1}{2}{\bf g}=(F^0+{\bf g}_{F^0}-\frac{1}{2}{\bf g})+{\bf g}\subseteq F_k^0+\frac{1}{2}{\bf g}_k+{\bf g}=F_{k+m}+{\bf g}\subset\partial(P)+{\bf g}$$ and $$F^0+{\bf g}_{F^0}+\frac{1}{2}{\bf g}\subset P.$$ Thus $F^0+{\bf g}_{F^0}+\frac{1}{2}{\bf g}\subset\partial[P\cap(P+{\bf g})]$.

      Similarly, since $-{\bf g}_{F^0}+\frac{1}{2}{\bf g}=-({\bf g}_{F^0}-\frac{1}{2}{\bf g})$, we have $$F^0-{\bf g}_{F^0}+\frac{1}{2}{\bf g}\subseteq F_k^0-\frac{1}{2}{\bf g}_k=F_k\subset\partial(P)$$ and $$F^0-{\bf g}_{F^0}+\frac{1}{2}{\bf g}=F^0-({\bf g}_{F^0}+\frac{1}{2}{\bf g})+{\bf g}\subset P+{\bf g}.$$ Thus $F^0-{\bf g}_{F^0}+\frac{1}{2}{\bf g}\subset\partial[P\cap(P+{\bf g})]$.
  \item If $\mathbb{S}_{F^0}=\mathbb{S}_{R_{k+t}}\bigcup\{(1-\lambda_0)S_0\}$ for some $1\le t\le m-k$, then there is a vector $${\bf g}_{F^0}=\frac{1}{2}\sum_{i=k+t+1}^m\widetilde{\bf g}_i-\frac{1}{2}\sum_{l=1}^{w_1}(1-\lambda_l^1)\vec{S}_l^1-\frac{1}{2}\sum_{j=1}^{w_k}(1-\lambda_j^k)\vec{S}_j^k+\frac{1}{2}\sum_{i=k+m+1}^{k+t+m-1}\widetilde{\bf g}_i\in\sum_{\mathbb{S}_{P^*}\setminus\mathbb{S}_{F^0}}T$$ satisfying that $F^0\pm{\bf g}_{F^0}$ are the facets of $P^*$, thus $F^0\pm{\bf g}_{F^0}+\frac{1}{2}{\bf g}$ are the facets of $P^*+\frac{1}{2}{\bf g}$, and $F^0\subseteq F_{k+t}^0$.

      Similar to (1), we also have $F_0\pm{\bf g}_{F^0}\pm\frac{1}{2}{\bf g}\subset P.$ Note that $${\bf g}_{F^0}-\frac{1}{2}{\bf g}=\frac{1}{2}{\bf g}_{k+t}-\frac{1}{2}\lambda_0\vec{S}_0,$$ we have $$F^0+{\bf g}_{F^0}+\frac{1}{2}{\bf g}=(F^0+{\bf g}_{F^0}-\frac{1}{2}{\bf g})+{\bf g}\subseteq F_{k+t}^0+\frac{1}{2}{\bf g}_{k+t}+{\bf g}=F_{k+t+m}+{\bf g}\subset\partial(P)+{\bf g}$$ and $$F^0+{\bf g}_{F^0}+\frac{1}{2}{\bf g}\subset P.$$ Thus $F^0+{\bf g}_{F^0}+\frac{1}{2}{\bf g}\subset\partial[P\cap(P+{\bf g})]$.

      Similarly, since $-{\bf g}_{F^0}+\frac{1}{2}{\bf g}=-({\bf g}_{F^0}-\frac{1}{2}{\bf g})$, we have $$F^0-{\bf g}_{F^0}+\frac{1}{2}{\bf g}\subseteq F_{k+t}^0-\frac{1}{2}{\bf g}_{k+t}=F_{k+t}\subset\partial(P)$$ and $$F^0-{\bf g}_{F^0}+\frac{1}{2}{\bf g}=F^0-({\bf g}_{F^0}+\frac{1}{2}{\bf g})+{\bf g}\subset P+{\bf g}.$$ Thus $F^0-{\bf g}_{F^0}+\frac{1}{2}{\bf g}\subset\partial[P\cap(P+{\bf g})]$.
\end{enumerate}

Above all, we have $\bigcup\limits_{F\in B(P^*,E^*)+\frac{1}{2}{\bf g}}F\subset\partial[P\cap(P+{\bf g})]$, equivalently, $P\cap(P+{\bf g})$ also has a belt determined by a translate of $(1-\lambda_0)E$. In fact, we can take $E^*=(1-\lambda_0)E_m$. Then by Lemma 4.2 and the Minkowski sum formula of $P$, we get $P^*+\frac{1}{2}{\bf g}\subset P\cap(P+{\bf g})$.

Then for $k=m$, we can also get the analogous result. Notice that if $k=m$ and all $\lambda_j^m=1$, then $$P\cap(P+{\bf g})=\sum_{l=1}^{w_1}(1-\lambda_l^1)S_l^1+(1-\lambda_0)S_0+\frac{1}{2}{\bf g}_1+\frac{1}{2}\sum_{l=1}^{w_1}\lambda_l^1\vec{S}_l^1,$$ which is a centrally symmetric polygon $F_{m+1}^*\subseteq F_{m+1}$. And we have $B(P\cap(P+{\bf g}),E^*)=\{F_{m+1}^*\}$.

Thus this proof is finished. \hfill{$\Box$}

\medskip
In particular, take all $\lambda_j^k=1$, $\lambda_l^1=0$ and $\lambda_0=0$ in the assumption of Lemma 4.3, and we directly obtain Corollary 4.1. Without loss of generality, we proceed to take ${\bf g}=\vec{S}_0$, and we will get Corollary 4.2 following.

\medskip\noindent
{\bf Corollary 4.1.} {\it Assume that ${\bf g}=\sum\limits_{i=2}^{k}\widetilde{\bf g}_i$ for $2\le k\le m$. Let $P^*$ be a zonotope with $\mathbb{S}_{P^*}=\{S_0\}\bigcup\mathbb{S}_{R_1}\bigcup(\bigcup\limits_{i=k+1}^{m}\mathbb{S}_{R_i})$, then we have $$P^*+\frac{1}{2}{\bf g}\subset P\cap(P+{\bf g}),$$ where $P^*=\sum_{S\in\mathbb{S}_{P^*}}S$. Further, $P\cap(P+{\bf g})$ and $P^*+\frac{1}{2}{\bf g}$ have the same belt $B(E)+\frac{1}{2}{\bf g}$.}

\medskip\noindent
{\bf Corollary 4.2.} {\it $P\cap(P+\vec{S}_0)=P^*+\frac{1}{2}\vec{S}_0$ where $P^*$ is a zonotope with $\mathbb{S}_{P^*}=\mathbb{S}_P\setminus\{S_0\}$.}

\medskip\noindent
{\bf Proof.} By Lemma 4.3, we have known that $P^*+\frac{1}{2}\vec{S}_0\subset P\cap(P+\vec{S}_0)$, then it's sufficient to show that $\partial(P^*)+\frac{1}{2}\vec{S}_0=\partial[P\cap(P+\vec{S}_0)]$. Notice that $E$ is a translate of $S_0$.

Let $$F^0=S_{j_1}+S_{j_2}+\cdots+S_{j_l}$$ be an original facet of $P^*$ where $S_{j_{k'}}\in\mathbb{S}_{P^*}$ and $k'=1,2,\cdots,l$. Let $\mathbb{S}_{F^0}=\{S_{j_1}, S_{j_2},\cdots, S_{j_l}\}$. Since $\mathbb{S}_{P^*}=\mathbb{S}_P\setminus\{S_0\}$, it follows by Lemma 4.1 that $F^0$ is (a part of) an original facet of $P$ and $P^*$ has all original facets of $P$ except for those facets parallel to $F'\in B(E)$.
\begin{enumerate}
  \item $\mathbb{S}_{F^0}=\mathbb{S}_{F_k^0}\setminus\{S_0\}$ for some $F_k\in B(E)$ with $1\le k\le m$. Since $F^0$ is an original facet of $P^*$, there exists a vector $${\bf g}_{F^0}=\sum_{i=k+1}^{m+k-1}\widetilde{\bf g}_i={\bf g}_k$$ satisfying that $F^0\pm{\bf g}_{F^0}$ are two facets of $P^*$ and $F^0\pm{\bf g}_{F^0}+\frac{1}{2}\vec{S}_0$ are two facets of $P^*+\frac{1}{2}\vec{S}_0$, and $F^0\subset F_k^0$. Similar to Lemma 4.3, we also have $$F^0\pm{\bf g}_{F^0}\pm\frac{1}{2}\vec{S}_0\subset P.$$ Further, $$F^0\pm\frac{1}{2}\vec{S}_0\pm{\bf g}_{F^0}\subset F_k^0\pm{\bf g}_k\subset\partial(P).$$ Then we have $$F^0\pm{\bf g}_{F^0}+\frac{1}{2}\vec{S}_0\subset\partial[P\cap(P+\vec{S}_0)].$$
  \item $\mathbb{S}_{F^0}\not=\mathbb{S}_{F_k^0}\setminus\{S_0\}$ for any $F_k\in B(E)$. For convenience, let $M$ be another edge of $P$ which is not a translate of $E$ and let $B(M)$ be the belt of $P$ determined by $M$. Apply Lemma 4.3 to $M$, we have that $P\cap(P+\vec{S}_0)$ and $P^*+\frac{1}{2}\vec{S}_0$ have the same belt determined by a translate of $M$. In fact, all the facets in this belt except for two facets parallel to the facets in $B(E)$ are the translates of the facets of $P$.
\end{enumerate}

 Therefore, $$\partial(P^*)+\frac{1}{2}\vec{S}_0\subset\partial[P\cap(P+\vec{S}_0)].$$ Since both $P^*$ and $P\cap(P+\vec{S}_0)$ are polytopes, it is easy to see that $\partial(P^*)+\frac{1}{2}\vec{S}_0=\partial[P\cap(P+\vec{S}_0)]$ and $P^*+\frac{1}{2}\vec{S}_0=P\cap(P+\vec{S}_0).$

 Then the proof of Corollary 4.2 has been finished. \hfill{$\Box$}

\medskip
The results explored above are mainly about the intersection of two translates of a three-dimensional zonotope. In fact, for centrally symmetric polygons in the plane, we also have obtained a similar consequence which is shown in Lemma 3.3 of Han, Sriamorn, Yang and Zong \cite{hsyz}. For convenience, we rephrase this result here. Let $F$ be a centrally symmetric polygon centered at {\bf o} with the generator set $\mathbb{S}_{F}=\{T_1,T_2,..,T_l\}$ for some positive integer $l$, equivalently, $$F=T_1+T_2+\cdots+T_l.$$ Let $M_1$,$M_2$,...,$M_{2l}$ be the edges of $F$ in a circular order where let ${\bf v}_i$ and ${\bf v}_{i+1}$ are two vertices of $M_i$ and ${\bf v}_{2l+1}={\bf v}_1$. Clearly, $M_{l+i}=-M_{l}$ for $1\le i\le l$. Let $\vec{T}_i=\vec{M}_i$ be the translation vector of ${\bf v}_i$ in $M_i$ such that ${\bf v}_{i+1}={\bf v}_i+\vec{M_i}$ and let ${\bf t}_i$ be the translation vector of $M_i$ in $F$ such that $M_{i+l}=M_i+{\bf t}_i$ for $1\le i\le l$. Obviously, $\vec{M}_{i+l}=-\vec{M}_i$. It is easy to see that $${\bf t}_i=\vec{M}_{i+1}+\cdots+\vec{M}_{l+i-1}.$$ Then we can have the following corollary.

\medskip\noindent
{\bf Corollary 4.3 (Han, Sriamorn, Yang and Zong \cite{hsyz}).} Assume that ${\bf t}=\sum_{i=2}^{k-1}\vec{M}_i+\lambda_k\vec{M}_k+\lambda_1\vec{M}_1$ for $2\le k\le l$, $0\le\lambda_k\le1$ and $0\le\lambda_1<1$. Then we have that $$F\cap(F+{\bf t})=F^*+\frac{1}{2}{\bf t},$$ where $F^*=(1-\lambda_1)T_1+(1-\lambda_k)T_k+T_{k+1}+\cdots+T_l$.

\medskip\noindent
{\bf Lemma 4.4.} {\it If $P$ is a three-dimensional non-prism zonotope with $\#B(E)=10$, then there must exist another edge $M$ of $P$ with $\#B(M)\ge8$.}

\medskip\noindent
{\bf Proof.} Firstly show that there exists another edge $M$ of $P$ with $\#B(M)\ge8$ if all the facets in $B(E)$ are parallelograms, then other cases can be obtained by adding generator segments. Notice that $B(E)=\{F_1,F_2,....,F_{10}\}$. For clarity, as introduced in Section 3, let $M_i$ be one edge of $F_i$ not parallel to $E$, and we assume that the generator set $\mathbb{S}_P$ of $P$ satisfies $$\mathbb{S}_P=\{S_0, S_1,..., S_5\},$$ that is $P=S_0+S_1+...+S_5$, where $E$ is the translate of $S_0$ and $M_i$ is the translate of $S_i$. Moreover, let $F_i^0=S_0+S_i$ and $\mathbb{S}_{F_i^0}=\{S_0,S_i\}$ for $i=1,2,3,4,5$. Because of the assumption that $P$ is not a prism, $\#B(M_i)\ge6$ holds for each $M_i$. Without loss of generality, let $i=1$. If $\#B(M_1)\ge8$, there is nothing to prove. Then we need to take into $\#B(M_1)=6$ and have two cases following.

\smallskip\noindent
{\bf Case 1.} There are three segments in $\{S_2,...,S_5\}$ coplanar with $S_1$. Without loss of generality, assume that $S_1$, $S_2$, $S_3$ and $S_4$ are coplanar, then we have $\#B(M_5)=10$ since $S_5$ cannot be coplanar with any two segments in $\{S_2,S_3,S_4\}$, otherwise all $S_i$ except for $S_0$ are coplanar which contradicts our assumption.

\smallskip\noindent
{\bf Case 2.} There are two segments in $\{S_2,...,S_5\}$ coplanar with $S_1$. Without loss of generality, assume that $S_1$, $S_2$ and $S_3$ are coplanar, then we have that $S_1$, $S_4$ and $S_5$ are coplanar by $\#B(M_1)=6$. Thus we have $\#B(M_i)=8$ for $i=2,3,4,5$. Take $S_2$ as an example, and it is easy to see that $S_2$ cannot be coplanar with $S_4$, $S_5$ simultaneously, then $F_2^0=S_2+S_0$, $S_2+S_1+S_3$, $S_2+S_4$ and $S_2+S_5$ are original facets of $P$ and $\#B(M_2)=8$.

Then we consider about adding one generator segment $S_1^*$ in $\mathbb{S}_P$ which is not linear with any $S_i$ for $i=0,1,...,5$, and replace $\mathbb{S}_P\cup\{S_1^*\}$ by $\mathbb{S}_P$. Without loss of generality, we can always assume $4\le\#B(M_1)\le6$ otherwise the proof have to be finished. Let $M_1^*$ be one translate of $S_1^*$ which is also an edge of $P$ and $S_1^*\in\mathbb{S}_{F_j^0}$ for some $j\in\{1,2,3,4,5\}$.
\begin{enumerate}
  \item $\#B(M_1)=4$, that is $S_1$, $S_2$,..., $S_5$ are coplanar, then $S_1^*\in\mathbb{S}_{F_1^0}$. Thus, we have $\#B(M_1^*)=10$.
  \item $\#B(M_1)=6$. For convenience, let $F_1^{*,0}$, $F_2^{*,0}$ and $F_3^{*,0}$ be the original facets of $P$ corresponding to the facets in $B(M_1)$ and $S_0,S_1\in\mathbb{S}_{F_1^{*,0}}=\mathbb{S}_{F_1^0}$.

      If $S_1^*$ are coplanar with $S_0$ and $S_1$, then $\#\mathbb{S}_{F_2^{*,0}}\ge3$ or $\#\mathbb{S}_{F_3^{*,0}}\ge3$ by the principle of pigeon hole. Without loss of generality, let $\#\mathbb{S}_{F_2^{*,0}}\ge\#\mathbb{S}_{F_3^{*,0}}$, $\{S_1,S_2,S_3\}\subseteq\mathbb{S}_{F_2^{*,0}}$ and $S_5\in\mathbb{S}_{F_3^{*,0}}$, then we have $\#B(M_5)\ge8$ since $S_5$ cannot be coplanar with $S_2$ and $S_3$ simultaneously.

      If $S_1^*$ are not coplanar with $S_0$ and $S_1$, then $\#\mathbb{S}_{F_2^{*,0}}\ge4$ or $\#\mathbb{S}_{F_3^{*,0}}\ge4$ by the principle of pigeon hole. Without loss of generality, let $S_1^*$ is coplanar with $S_0$ and $S_2$ and $\{S_1,S_2,S_3,S_4\}\subseteq\mathbb{S}_{F_2^{*,0}}$,  then we have $S_1^*\in\mathbb{S}_{F_3^{*,0}}$ and $\#B(M_1^*)\ge8$.
\end{enumerate}

Then we proceed to add the second generator segment $S_2^*$ which is not linear with $S_1^*$ and any $S_i$, and replace $\mathbb{S}_P\cup\{S_2^*\}$ by $\mathbb{S}_P$, then it suffices to consider the case $\#B(M_1)=6$ by the analysis above.

If both $S_1^*$ and $S_2^*$ are in $\mathbb{S}_{F_j^0}$ for some $j=1,2,..,5$, then we have $\#B(M_i)\ge8$ for each $i\not=j$.

If $S_2^*$ cannot belong to some $\mathbb{S}_{F_j^0}$ which contains $S_1^*$. Without loss of generality, let $S_2^*\in\mathbb{S}_{F_2^0}$. Then, we can have $\#\mathbb{S}_{F_2^{*,0}}\ge4$ or $\#\mathbb{S}_{F_3^{*,0}}\ge4$. Without loss of generality, let $\{S_1,S_2,S_3,S_4\}\subseteq\mathbb{S}_{F_2^{*,0}}$, then $S_2^*\in\mathbb{S}_{F_3^{*,0}}$ and $\#B(M_2^*)\ge8$.

In fact, by the analysis above, we can only consider $2\le\#\mathbb{S}_{F_i^0}\le3$ for each $i$ enough. Now we can successively add three segments $S_3^*$, $S_4^*$ and $S_5^*$ to $\mathbb{S}_P$, respectively. Similarly, it is easy to get that there always exists another edge $M$ of $P$ with $\#B(M)\ge8$.

Then the proof of Lemma 4.4 is finished. \hfill{$\Box$}

\medskip\noindent
{\bf Lemma 4.5.} {\it If $P$ is a three-dimensional non-prism zonotope with $\#B(E)=8$, then there must exist another edge $M$ of $P$ with $\#B(M)\ge8$.}

\medskip\noindent
{\bf Proof.} Similar to Lemma 4.4, firstly show that there exists another edge $M$ of $P$ with $B(M)=8$ if all the facets in $B(E)$ are parallelograms, then other cases can be obtained by adding generator segments. For clarity, as introduced in Section 3, let $M_i$ be one edge of $F_i$ not parallel to $E$, and we assume that the generator set $\mathbb{S}_P$ of $P$ satisfies $$\mathbb{S}_P=\{S_0, S_1,S_2,S_3,S_4\},$$ that is $P=S_0+S_1+...+S_4$, where $E$ is the translate of $S_0$ and $M_i$ is the translate of $S_i$. Moreover, let $F_i^0=S_0+S_i$ and $\mathbb{S}_{F_i^0}=\{S_0,S_i\}$ for $i=1,2,3,4$. Because of the assumption that $P$ is not a prism, $\#B(M_i)\ge6$ holds for each $M_i$. Without loss of generality, let $i=1$. If $\#B(M_1)=8$, there is nothing to prove. Then we need to take into $\#B(M_1)=6$.

For $\#B(M_1)=6$, there need two segments in $\{S_2,S_3,S_4\}$ coplanar with $S_1$ by the principle of pigeon hole. Without loss of generality, assume that $S_1$, $S_2$ and $S_3$ are coplanar, then we have $\#B(M_4)=8$ since $M_4$ cannot be coplanar with $S_2$ and $S_3$ simultaneously, otherwise all $S_i$ except for $S_0$ are coplanar which contradicts our assumption.

Then we consider about adding one generator segment $S_1^*$ in $\mathbb{S}_P$ which is not linear with any $S_i$ for $i=0,1,...,4$, and replace $\mathbb{S}_P\cup\{S_1^*\}$ by $\mathbb{S}_P$. Without loss of generality, we can always assume $4\le\#B(M_1)\le6$ otherwise the proof have to be finished. Let $M_1^*$ be one translate of $S_1^*$ which is also an edge of $P$ and $S_1^*\in\mathbb{S}_{F_j^0}$ for some $j\in\{1,2,3,4\}$.
\begin{enumerate}
  \item $\#B(M_1)=4$, that is $S_1$, $S_2$, $S_3$ and $S_4$ are coplanar, then $S_1^*\in\mathbb{S}_{F_1^0}$. Thus, we have $\#B(M_1^*)=8$.
  \item $\#B(M_1)=6$. For convenience, let $F_1^{*,0}$, $F_2^{*,0}$ and $F_3^{*,0}$ be the original facets of $P$ corresponding to the facets in $B(M_1)$ and $S_0,S_1\in\mathbb{S}_{F_1^{*,0}}=\mathbb{S}_{F_1^0}$.

      If $S_1^*\in\mathbb{S}_{F_1^0}$, then $\#\mathbb{S}_{F_2^{*,0}}\ge3$ or $\#\mathbb{S}_{F_3^{*,0}}\ge3$ by the principle of pigeon hole. Without loss of generality, let $\#\mathbb{S}_{F_2^{*,0}}\ge\#\mathbb{S}_{F_3^{*,0}}$ and $\mathbb{S}_{F_2^{*,0}}=\{S_1,S_2,S_3\}$, then $\mathbb{S}_{F_3^{*,0}}=\{S_1,S_4\}$ and $S_4$ cannot be coplanar with $S_2$ and $S_3$ simultaneously. Thus $\{S_4+S_0, S_4+S_1, S_4+S_2+S_1^*, S_4+S_3\}$ or $\{S_4+S_0, S_4+S_1, S_4+S_2, S_4+S_3+S_1^*\}$ is the set of original facets corresponding to the facets in $B(M_4)$ and $\#B(M_4)=8$.

      If $S_1^*\notin\mathbb{S}_{F_1^0}$, then we also have $\#\mathbb{S}_{F_2^{*,0}}\ge3$ or $\#\mathbb{S}_{F_3^{*,0}}\ge3$. Without loss of generality, let $S_1^*\in\mathbb{S}_{F_2^0}$, $\{S_1,S_2,S_3\}\subseteq\mathbb{S}_{F_2^{*,0}}$ and $S_1^*\in\mathbb{S}_{F_3^{*,0}}$. If $S_4\in\mathbb{S}_{F_2^{*,0}}$, then $\#B(M_1^*)=8$. If $S_4\in\mathbb{S}_{F_3^{*,0}}$, then $\#B(M_4)=8$.
\end{enumerate}

Then we proceed to add the second generator segment $S_2^*$ which is not linear with $S_1^*$ and any $S_i$, and replace $\mathbb{S}_P\cup\{S_2^*\}$ by $\mathbb{S}_P$, then it suffices to consider the case $\#B(M_1)=6$ by the analysis above. Similar to Lemma 4.4, we can also consider $2\le\#\mathbb{S}_{F_i^0}\le3$ for each $i$. Then $S_2^*$ cannot belong to some $\mathbb{S}_{F_{i}^0}$ which contains $S_1^*$. Without loss of generality, let $S_2^*\in\mathbb{S}_{F_2^0}$. Then, we can have $\#\mathbb{S}_{F_2^{*,0}}\ge3$ or $\#\mathbb{S}_{F_3^{*,0}}\ge3$. Without loss of generality, let $\{S_1,S_2,S_3\}\subseteq\mathbb{S}_{F_2^{*,0}}$, then $S_2^*\in\mathbb{S}_{F_3^{*,0}}$. Like (2) above, if $S_4\in\mathbb{S}_{F_2^{*,0}}$, then $\#B(M_1^*)\ge8$. If $S_4\in\mathbb{S}_{F_3^{*,0}}$, then $\#B(M_4)\ge8$.

Now we can successively add two segments $S_3^*$ and $S_4^*$ to $\mathbb{S}_P$, and it is easy to get that there always exists another edge $M$ of $P$ with $\#B(M)\ge8$.

Then the proof of Lemma 4.5 is finished. \hfill{$\Box$}

\medskip
In fact, under the conditional assumptions in Lemma 4.4 and Lemma 4.5, if $P$ is also a fivefold translative tile, then it follows by Lemma 2.4 that there must exist another edge $M$ of $P$ with $8\le\#B(M)\le10$.

\medskip
For any facet $F$ of $P$, let $\#F$ denote the number of edges of $F$. Then we have the following result.

\medskip\noindent
{\bf Lemma 4.6.} For any facet $F$ of the three-dimensional fivefold translative tile $P$, we have $$4\le\#F\le10.$$

\medskip\noindent
{\bf Proof.} Since $P$ is a zonotope with $\mathbb{S}_P=\{S_0,S_1,...,S_w\}$, we have $\#F$ is an even number. For the contrary, suppose that there is  a facet $F'$ of $P$ with $\#F'\ge 12$. Without loss of generality, assume that $F'=F_1$ and $\mathbb{S}_{F_1^0}=\{S_0,S_1,...,S_{k_0}\}$ for $6\le k_0<w$. Take one segment $S\in\mathbb{S}\setminus\mathbb{S}_{F_1^0}$ and let $M$ be one edge of $P$ which is parallel to $S$, then we consider $\#B(M)$. By Lemma 4.1, we have known that $\mathbb{S}_{F_1^0}$ is a maximal coplanar segments set, that is $S$ cannot be simultaneously coplanar with any two different segments in $\mathbb{S}_{F_1^0}$, then the number of the original facets containing $S$ are not less than $2k_0$, consequently, $$\#B(M)\ge12$$ which contradicts the result in Lemma 2.4.

Thus we have finished the proof of Lemma 4.6. \hfill{$\Box$}

\medskip\noindent
{\bf Remark 4.1.} In fact, for the definition of the non-proper fivefold translative tiling $P+X$ relate to $E$ shown in Definition 3.4, we always have $relint(E')\cap\partial(P')\not=\emptyset$ for some $E'\in A(E)+X$ and some $P'\in P+X$, otherwise by Lemma 4.3 we can get that $P+X$ is a proper fivefold translative tiling relate to $E$.

\medskip\noindent
{\bf Remark 4.2.} Generally speaking, for any two points ${\bf x},{\bf y}\in X$, if ${\bf g}={\bf y}-{\bf x}$ in Lemma 4.3 satisfies that $(P+{\bf x})\cap(P+{\bf y})$ has a belt determined by a translate of $(1-\lambda_0)E$ for $\lambda_0\not=1$, then, for convenience in the following discussion, we write such a belt as $B(E,{\bf x},{\bf y})$.

\vspace{1cm}\noindent
{\LARGE\bf 5. Projections of fivefold translative tiles}

\bigskip\noindent
As we have known about the two-, three- and fourfold translative tiles in $\mathbb{E}^3$, it is easy to see that their projections along any edge of the tiles are either parallelograms or centrally symmetric hexagons each of which is also a two-, three- or fourfold translative tile in $\mathbb{E}^2$. Then what about the projection of a three-dimensional fivefold translative tile $P$ along any of its edges? To provide a clear explanation of the proof of this theorem, we need a following result studied by Bolle in 1994.

\medskip\noindent
{\bf Lemma 5.1 (Bolle \cite{boll}).} {\it A convex polygon is a k-fold lattice tile for a lattice $\Lambda$ and some positive integer k if and only if the following conditions are satisfied:
\begin{enumerate}
  \item It is centrally symmetric.
  \item When it is centered at the origin, in the relative interior of each edge $E$ there is a point of $\frac{1}{2}\Lambda$.
  \item If the midpoint of $E$ is not in $\frac{1}{2}\Lambda$ then $E$ is a lattice vector of $\Lambda$.
\end{enumerate}
}

\medskip
Let $P$ be a fivefold translative tile in $\mathbb{E}^3$ and let $E$ be an edge of $P$. By Lemma 2.4, we have known that $\#B(E)$ = 4, 6, 8 or 10. Without loss of generality, we assume that $E$ is parallel to {\it z}-axis and define the positive direction of {\it z}-axis upward. Then we have the following results about the projection of $P$ along $E$.

\medskip\noindent
{\bf Lemma 5.2.} {\it If $P$ is a fivefold translative tile with $\#B(E)=4$ or $6$, then the projection of $P$ along $E$ is a parallelogram or a centrally symmetric hexagon.}

\smallskip
The result in Lemma 5.2 is apparent, then the proof is omitted here.

\medskip\noindent
{\bf Lemma 5.3.} {\it If $P$ is a fivefold translative tile with $\#B(E)=8$, then the projection of $P$ along $E$ must be a fivefold translative tile in $\mathbb{E}^2$.}

\medskip\noindent
{\bf Proof.} Assume that $E'\in A(E)+X$ and let {\bf p} be a proper point in $E'$. By Lemma 3.3, we have $$\varphi({\bf p})\ge\left\lceil\frac{4-3}{2}\right\rceil=1.$$ By Lemma 3.2, we have $$\varpi({\bf p})=\kappa({\bf p})\cdot\frac{3}{2}+\ell({\bf p})\cdot\frac{1}{2},$$ where $\kappa({\bf p})$ is a positive integer and $\ell({\bf p})$ is a nonnegative integer. Notice that $\ell({\bf p})$ is the number of the facets in $B(E)+X$ which take {\bf p} as a relative interior point. Then, we have the following cases:

\smallskip\noindent
{\bf Case 1.} $\varpi({\bf p})\ge 5$ holds for some proper point {\bf p} in $E'$. Since $\varphi({\bf p})\ge 1$, we get $$\varpi({\bf p})+\varphi({\bf p})\ge 5+1=6,$$ which contradicts $\varpi({\bf p})+\varphi({\bf p})=5$.

\smallskip\noindent
{\bf Case 2.} $\varpi({\bf p})=4$ holds for some proper point {\bf p} in $E'$. Then the equation $$\kappa({\bf p})\cdot\frac{3}{2}+\ell({\bf p})\cdot\frac{1}{2}=4$$ has two groups of solutions $\kappa({\bf p})=1,\ell({\bf p})=5$ and $\kappa({\bf p})=2,\ell({\bf p})=2$, so there is an F-type translate $P+{\bf x}$ at {\bf p} for ${\bf x}\in X$. Suppose that ${\bf p}\in relint(F_i)+{\bf x}$. By the definition of a proper point, there exist ${\bf p}^*\in (E_i+{\bf x})\setminus\{C(E)+X\}$ and ${\bf p}^\bullet\in(E_{i+1}+{\bf x})\setminus\{C(E)+X\}$ such that $${\bf p}\in {\bf p}^*+R_i-{\bf c}_{E_i^0}={\bf p}^\bullet+R_i-{\bf c}_{E_{i+1}^0}.$$
By applying Lemma 3.1 to ${\bf p}^*$ and ${\bf p}^\bullet$, respectively, there are ${\bf x}^*, {\bf x}^\bullet\in X$ such that $${\bf p}^*\in B(E)+{\bf x}^*,$$  $${\bf p}^\bullet\in B(E)+{\bf x}^\bullet,$$ $$({\bf p}^*+R_i-{\bf c}_{E_i^0})\setminus\{{\bf p}^*\}\subset int(P)+{\bf x}^*$$ and $$({\bf p}^\bullet+R_i-{\bf c}_{E_{i+1}^0})\setminus\{{\bf p}^\bullet\}\subset int(P)+{\bf x}^\bullet.$$
Since $${\bf p}\in {\bf p}^*+R_i-{\bf c}_{E_i^0}={\bf p}^\bullet+R_i-{\bf c}_{E_{i+1}^0}.$$
we obtain $${\bf p}\in int(P)+{\bf x}^*$$ and $${\bf p}\in int(P)+{\bf x}^\bullet.$$
By the convexity of $P$, since $\#B(E)=8>4$, one can deduce that ${\bf x}^*\not={\bf x}^\bullet$. Consequently, we get $$\varphi({\bf p})\ge 2,$$ then $$\varpi({\bf p})+\varphi({\bf p})\ge 4+2=6,$$
which contradicts $\varpi({\bf p})+\varphi({\bf p})=5$.

\smallskip\noindent
{\bf Case 3.} $\varpi({\bf p})=3$ holds for some proper point {\bf p} in $E'$. Then the equation $$\kappa({\bf p})\cdot\frac{3}{2}+\ell({\bf p})\cdot\frac{1}{2}=3$$ has two groups of solutions $\kappa({\bf p})=1, \ell({\bf p})=3$ and $\kappa({\bf p})=2, \ell({\bf p})=0$.

\smallskip\noindent
{\bf Subcase 3.1.} $\kappa({\bf p})=1, \ell({\bf p})=3$. Then there are at most two dihedral adjacent wheels in $P+X^{\bf p}$ consisting of four E-type translates and three F-type translates. For convenience, let $X^{\bf p}=\{{\bf x}_1, {\bf x}_2,...,{\bf x}_7\}$ with ${\bf x}_1={\bf o}$. Without loss of generality, let $E'=E_6$.

\smallskip\noindent
{\bf Subcase 3.1.1.} There is only one dihedral adjacent wheel $P+X^{\bf p}$ at {\bf p}. For convenience, without loss of generality, we can assume that $P+{\bf x}_1$, $P+{\bf x}_2$,..., $P+{\bf x}_7$ are sequential in $P+X^{\bf p}$. Project $P+X^{\bf p}$ onto the plane {\it x}{\bf o}{\it y}, and we also use the $P+{\bf x}_i$ to denote its projection without confusion, as shown in FIG \ref{P8_7translates}.

\begin{figure}[h!]
 \subfigure[]
 {\begin{minipage}{0.8\linewidth}
 \includegraphics[scale=0.60]{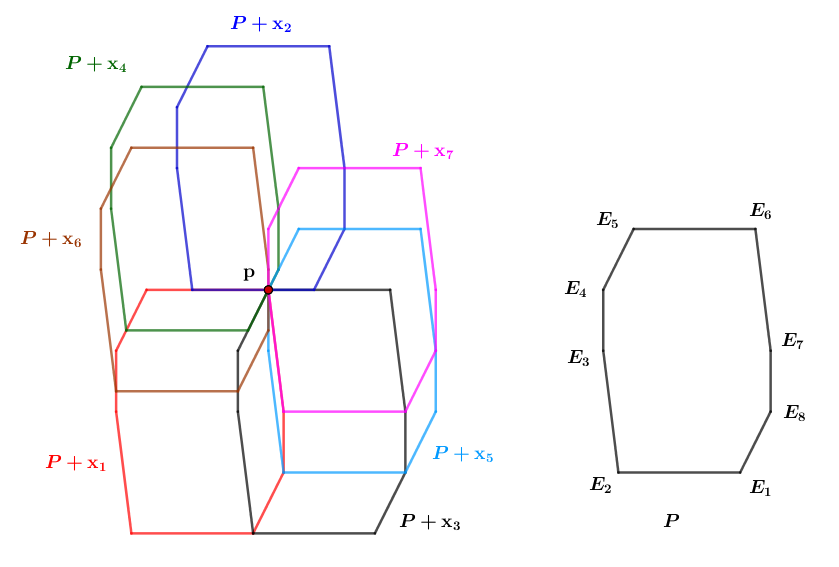}
 \end{minipage}}
 \hspace{1.0cm}
 \subfigure[]
 {\begin{minipage}{0.65\linewidth}
 \includegraphics[scale=0.56]{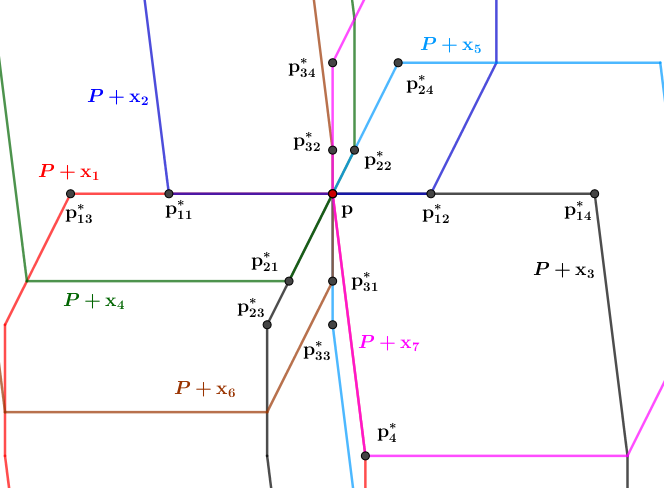}
 \end{minipage}}
 \caption{The only one dihedral adjacent wheel $P+X^{\bf p}$ in Subcase 3.1.1}
 \label {P8_7translates}
\end{figure}

Let ${\bf p}_{11}^*$, ${\bf p}_{12}^*$, ${\bf p}_{13}^*$ and ${\bf p}_{14}^*$ be the corresponding points of {\bf p} in $(F_1+{\bf x}_2)\cap(F_5+{\bf x}_1)$, $(F_1+{\bf x}_2)\cap(F_5+{\bf x}_3)$, $F_5+{\bf x}_1$ and $F_5+{\bf x}_3$, respectively. Let ${\bf p}_{21}^*$, ${\bf p}_{22}^*$, ${\bf p}_{23}^*$ and ${\bf p}_{24}^*$ be the corresponding points of {\bf p} in $(F_8+{\bf x}_4)\cap(F_4+{\bf x}_3)$, $(F_8+{\bf x}_4)\cap(F_4+{\bf x}_5)$, $F_4+{\bf x}_3$ and $F_4+{\bf x}_5$, respectively. Let ${\bf p}_{31}^*$, ${\bf p}_{32}^*$, ${\bf p}_{33}^*$ and ${\bf p}_{34}^*$ be the corresponding points of {\bf p} in $(F_7+{\bf x}_6)\cap(F_3+{\bf x}_5)$, $(F_7+{\bf x}_6)\cap(F_3+{\bf x}_7)$, $F_3+{\bf x}_5$ and $F_3+{\bf x}_7$, respectively. And let ${\bf p}_4^*$ be the corresponding point of {\bf p} in $(F_6+{\bf x}_1)\cap(F_2+{\bf x}_7)$. In fact, ${\bf p}_4^*={\bf p}-\widetilde{\bf g}_2$, we can also say that ${\bf p}_4^*$ is the corresponding point of {\bf p} in $F_2+{\bf x}_7$ or $F_6+{\bf x}_1$. Since {\bf p} is a proper point in $E'$, we have all ${\bf p}_{ij}^*\notin C(E)+X$ and ${\bf p}_4^*\notin C(E)+X$. Applying Lemma 3.1 to each ${\bf p}_{ij}^*$ and ${\bf p}_4^*$, there are points ${\bf y}_{ij}\in X^{{\bf p}_{ij}^*}$ and ${\bf y}_4\in X^{{\bf p}_4^*}$ satisfying that $${\bf p}\in int(P)+{\bf y}_{ij}$$ and $${\bf p}\in int(P)+{\bf y}_4$$ for all $i,j=1,2,3$. Since $\varpi({\bf p})=3$, we need $\varphi({\bf p})=2$, that is there are two I-type translates $P+{\bf z}_1$ and $P+{\bf z}_2$ at {\bf p}. It is easy to see that ${\bf y}_{i1}\not={\bf y}_{i2}$ and ${\bf y}_{i1}\not={\bf y}_{i3}$ by $\#B(E)=8>4$, then we have ${\bf y}_{i1}={\bf y}_{i4}$, ${\bf y}_{i2}={\bf y}_{i3}$ for all $i$. Without loss of generality, let ${\bf y}_{i1}={\bf z}_1$ for all $i$ and let ${\bf y}_4={\bf z}_1$. Then by observing the points ${\bf p}_{11}^*$, ${\bf p}_4^*$ and ${\bf p}_{34}^*$, we can see that $(P+{\bf x}_1)\cap(P+{\bf z}_1)$ has a belt $B(E,{\bf x}_1,{\bf z}_1)$ with $\#B(E,{\bf x}_1,{\bf z}_1)=4$ and $(P+{\bf x}_7)\cap(P+{\bf z}_1)$ has a belt $B(E,{\bf x}_7,{\bf z}_1)$ with $\#B(E,{\bf x}_7,{\bf z}_1)=4$, then it follows that $P+{\bf z}_1$ has a belt $B(E)+{\bf z}_1$ with $\#B(E)=6$ which contradicts $\#B(E)=8$.

\smallskip\noindent
{\bf Subcase 3.1.2.} There are two dihedral adjacent wheels $P+X_1^{\bf p}$ and $P+X_2^{\bf p}$ at {\bf p}. Without loss of generality, let $X_1^{\bf p}=\{{\bf x}_1,{\bf x}_2,...,{\bf x}_4,{\bf x}_5\}$ and $X_2^{\bf p}=\{{\bf x}_6,{\bf x}_7\}$. Clearly, $P+X_1^{\bf p}$ consists of four E-type translates and one F-type translate and $P+X_2^{\bf p}$ consists of two F-type translates. For convenience, without loss of generality, we can assume that $P+{\bf x}_1$, $P+{\bf x}_2$,..., $P+{\bf x}_5$ are sequential in $X_1^{\bf p}$ and $F_1+{\bf x}_2$ has {\bf p} as its relative interior point. Project $P+X_1^{\bf p}$ onto the plane {\it x}{\bf o}{\it y}, and we also use the $P+{\bf x}_i$ to denote its projection without confusion, as shown in FIG \ref{P8_Ftranslates}.

\begin{figure}[h!]
 \includegraphics[scale=0.6]{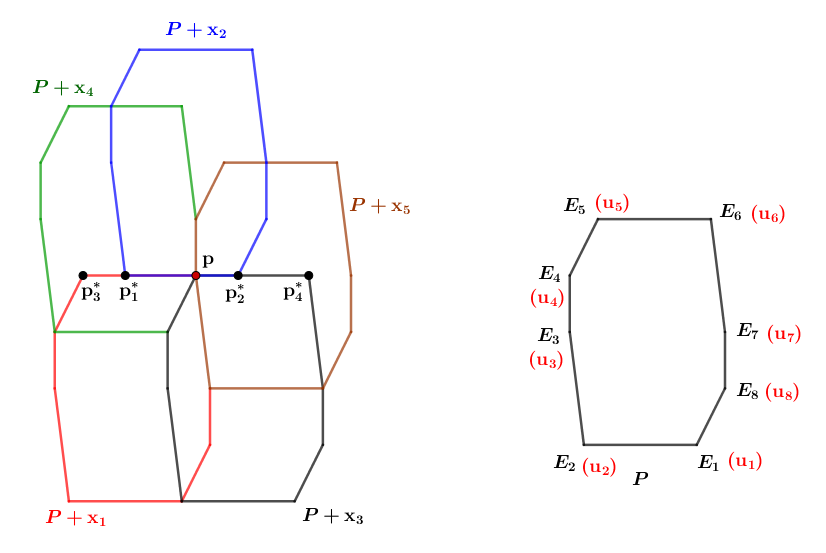}
 \caption{One dihedral adjacent wheel $P+X_1^{\bf p}$ in Subcase 3.1.2}\label{P8_Ftranslates}
\end{figure}

Let ${\bf p}_1^*$, ${\bf p}_2^*$, ${\bf p}_3^*$ and ${\bf p}_4^*$ be the corresponding points of {\bf p} in $(F_1+{\bf x}_2)\cap(F_5+{\bf x}_1)$, $(F_1+{\bf x}_2)\cap(F_5+{\bf x}_3)$, $F_5+{\bf x}_1$ and $F_5+{\bf x}_3$, respectively. Since {\bf p} is a proper point in $E'$, we have ${\bf p}_i^*\notin C(E)+X$ for $1\le i\le 4$. Moreover, by the infiniteness of the number of proper points in $E'$ and the finiteness of the number of translates at {\bf p} and ${\bf p}_i^*$, we can assume that {\bf p} and ${\bf p}_i^*$ are all proper points. Applying Lemma 3.1 to each ${\bf p}_i^*$, we have that there is a point ${\bf y}_i\in X^{{\bf p}_i^*}$ satisfying that $${\bf p}\in int(P)+{\bf y}_i.$$ Clearly, ${\bf y}_1\not={\bf y}_2$ and ${\bf y}_1\not={\bf y}_3$. Since $\varpi({\bf p})=3$, there only need two I-type translates at {\bf p} and it follows that $${\bf y}_1={\bf y}_4$$ and $${\bf y}_2={\bf y}_3.$$ By the definition of the corresponding point and Corollary 4.3, we have ${\bf p}_3^*={\bf p}+\widetilde{\bf g}_1$ and ${\bf p}_1^*\in {\bf p}+R_1-{\bf c}_{E_1^0}={\bf p}_3^*+R_1-{\bf c}_{E_2^0}$, then $${\bf p}_1^*\in relint(F_5)+{\bf x}_1$$ and $${\bf p}_1^*\in int[(P+{\bf x}_4)\cap(P+{\bf y}_2)].$$ Since $\varphi({\bf p}_1^*)\ge 2$, we have $\varpi({\bf p}_1^*)=2$ or 3.
\begin{enumerate}
  \item For $\varpi({\bf p}_1^*)=2$, the structure of the dihedral adjacent wheel at ${\bf p}_1^*$ is determined by $P+{\bf x}_1$ and $P+{\bf x}_2$. One can recall the definition of the structure of the dihedral adjacent wheel at ${\bf p}_1^*$ in Section 2.
  \item For $\varpi({\bf p}_1^*)=3$, it follows by Subcase 3.1.1 above that there are two dihedral adjacent wheels $P+X_1^{{\bf p}_1^*}$ and $P+X_2^{{\bf p}_1^*}$ at ${\bf p}_1^*$ in which $P+X_1^{{\bf p}_1^*}$ consists of four E-type translates and one F-type translate $P+{\bf x}_1$ and $P+X_2^{{\bf p}_1^*}$ consists of two F-type translates. Notice the structure of $P+X_1^{{\bf p}_1^*}$ is determined by $P+{\bf x}_1$ and $P+{\bf x}_2$.
\end{enumerate}
According to (1) and (2), there is always a dihedral adjacent wheel $P+X_1^{{\bf p}_1^*}$ consisting of four E-type translates and one F-type translate $P+{\bf x}_1$ and its structure is determined by $P+{\bf x}_1$ and $P+{\bf x}_2$, then we have $${\bf y}_1\in X_1^{{\bf p}_1^*}$$ and $${\bf p}_1^*\in relint(E_4)+{\bf y}_1.$$ By observation, there is a point ${\bf z}_1\in X_1^{{\bf p}_1^*}$ satisfying that ${\bf p}_1^*\in relint(E_1)+{\bf z}_1$. By Corollary 4.3, we have $${\bf p}_3^*\in relint(F_1)+{\bf z}_1.$$ Then similar to ${\bf p}_1^*$, there is also a dihedral adjacent wheel $P+X_1^{{\bf p}_3^*}$ consisting of four E-type translates and one F-type translate $P+{\bf z}_1$ and its structure is determined by $P+{\bf x}_1$ and $P+{\bf z}_1$. Then we have $${\bf y}_2\in X_1^{{\bf p}_3^*}$$ and $${\bf p}_3^*\in relint(E_3)+{\bf y}_2.$$
Analogously, for ${\bf p}_2^*$ and ${\bf p}_4^*$, there also respectively exist two dihedral adjacent wheels $P+X_1^{{\bf p}_2^*}$ and $P+X_1^{{\bf p}_4^*}$ both consisting of four E-type translates and one F-type translate, then we further have $${\bf y}_2\in X_1^{{\bf p}_2^*}, {\bf p}_2^*\in relint(E_7)+{\bf y}_2$$ and $${\bf y}_4\in X_1^{{\bf p}_4^*}, {\bf p}_4^*\in relint(E_8)+{\bf y}_1.$$ However, in fact, $E_4=E_3+\widetilde{\bf g}_3$ and $E_8=E_7-\widetilde{\bf g}_3$, then ${\bf p}_1^*$, ${\bf p}_2^*$, ${\bf p}_3^*$ and ${\bf p}_4^*$ cannot be coplanar in the plane containing $F_1+{\bf x}_2$.

\smallskip\noindent
{\bf Subcase 3.2.} $\kappa({\bf p})=2, \ell({\bf p})=0$ holds at all proper points ${\bf p}$ in any $E'\in A(E)+X$. Then $P+X^{\bf p}$ is a dihedral adjacent wheel of eight E-type translates $P+{\bf x}_1$,..., $P+{\bf x}_8$ in sequence. Without loss of generality, let $E'=E_6$ and let {\bf p} be a proper point in $E'$. Project $P+X^{\bf p}$ onto the plane {\it x}{\bf o}{\it y}, and we also use the $P+{\bf x}_i$ to denote its projection without confusion, shown in FIG \ref{P8_1k=2,l=0}.

 \begin{figure}[h!]
\includegraphics[scale=0.60]{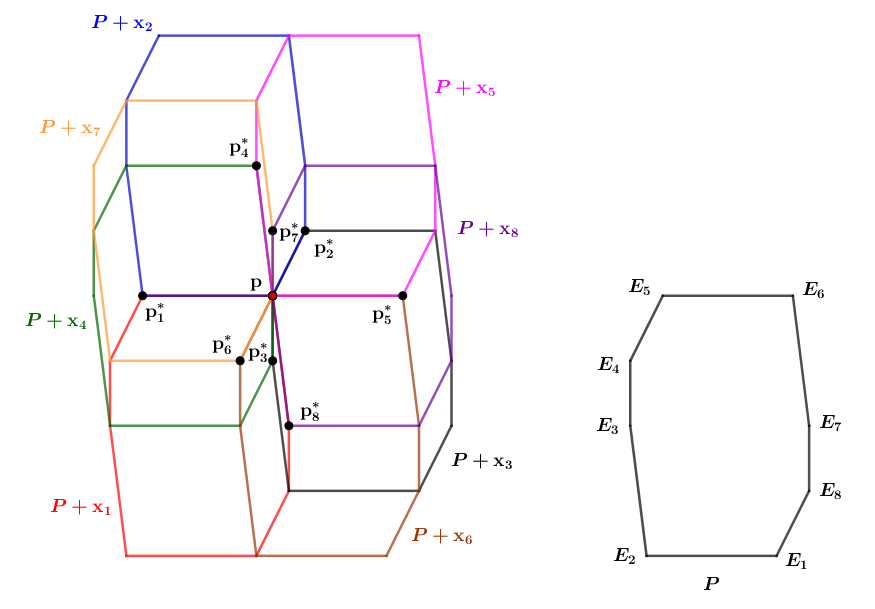}
\caption{The dihedral adjacent wheel $P+X^{\bf p}$ in Subcase 3.2}\label{P8_1k=2,l=0}
\end{figure}

Let ${\bf p}_1^*$, ${\bf p}_2^*$, ${\bf p}_3^*$, ${\bf p}_4^*$, ${\bf p}_5^*$, ${\bf p}_6^*$, ${\bf p}_7^*$ and ${\bf p}_8^*$ be the corresponding points of {\bf p} on $F_5+{\bf x}_1$, $F_8+{\bf x}_2$, $F_3+{\bf x}_3$, $F_6+{\bf x}_4$, $F_1+{\bf x}_5$, $F_4+{\bf x}_6$, $F_7+{\bf x}_7$ and $F_2+{\bf x}_8$, respectively. Similar to Subcase 3.1, we can also assume that {\bf p} and ${\bf p}_i^*$ are all proper points. Since $\varpi({\bf p})=3$, every ${\bf p}_i^*$ is an interior point of exact two of these eight translates. Consequently, for every ${\bf p}_i^*$, there are two different translates $P+{\bf y}_i$ and $P+{\bf y}_i'$ in $P+X^{{\bf p}_i^*}$ both contain {\bf p} as an interior point.

On the other hand, it can be easily deduced that there is only one point ${\bf x}\in X$ such that both ${\bf p}_1^*$ and ${\bf p}_2^*$ belong to $B(E)+{\bf x}$ and ${\bf p}\in int(P)+{\bf x}$ by the structure of the dihedral adjacent wheels at ${\bf p}_i^*$ for i=1,2. Therefore, at least one of the two points ${\bf y}_2$ and ${\bf y}_2'$ is different from both ${\bf y}_1$ and ${\bf y}_1'$. Then we get $$\varphi({\bf p})\ge3$$ and $$\varpi({\bf p})+\varphi({\bf p})\ge6,$$ which contradicts the assumption.

\smallskip\noindent
{\bf Case 4.} $\varpi({\bf p})=2$ holds for some proper point ${\bf p}\in E'$. Then the equation $$\kappa({\bf p})\cdot\frac{3}{2}+\ell({\bf p})\cdot\frac{1}{2}=2$$ has only one group of solutions $$\kappa({\bf p})=1, \ell({\bf p})=1.$$ In other words, $P+X^{\bf p}$ is a dihedral adjacent wheel of five translates. For convenience, we still concentrate on the structure of $P+X^{\bf p}$ with $X^{\bf p}=\{{\bf x}_1,...,{\bf x}_5\}$ as $P+X_1^{\bf p}$ in Subcase 3.1.2 shown in FIG \ref{P8_Ftranslates}, where ${\bf x}_1={\bf o}$ and $E'=E_6$.

Guaranteed by linear transformation, we can assume that the facets $F_1$ and $F_3$ of $P$ are vertical to {\it y}-axis and {\it x}-axis, respectively.

As introduced in Subcase 3.1.2, ${\bf p}_1^*$, ${\bf p}_2^*$, ${\bf p}_3^*$ and ${\bf p}_4^*$ are the corresponding points of {\bf p} on $(F_1+{\bf x}_2)\cap(F_5+{\bf x}_1)$, $(F_1+{\bf x}_2)\cap(F_5+{\bf x}_3)$, $F_5+{\bf x}_1$ and $F_5+{\bf x}_3$ and all ${\bf p}_i^*$ are proper points in their corresponding edges, respectively. And $$\varpi({\bf p}_1^*)=\varpi({\bf p}_2^*)=\varpi({\bf p}_3^*)=\varpi({\bf p}_4^*)=2$$ and ${\bf y}_1={\bf y}_4$ and ${\bf y}_2={\bf y}_3$ cannot holds simultaneously.

Project $P+X^{\bf p}$ and $P+{\bf y}_j$ onto the plane {\it x}{\bf o}{\it y}, and we still use $P+{\bf x}_i$, $P+{\bf y}_j$ and ${\bf p}_k^*$ to denote their projections without confusion. We use ${\bf u}_l$ to denote the projection of $E_l$ for $l=1,2,...,8$. For convenience, we write ${\bf y}_j=(y_1^j, y_2^j, y_3^j)$, ${\bf p}_k^*=(p_1^k, p_2^k, p_3^k)$, ${\bf u}_l=(u_1^l, u_2^l)$. Since ${\bf y}_1={\bf y}_4$ and ${\bf y}_2={\bf y}_3$ are symmetric, it is sufficient to deal with two subcases.

\smallskip\noindent
{\bf Subcase 4.1.} ${\bf y}_2={\bf y}_3$. Let ${\bf q}_1'$ and ${\bf q}_2'$ be the two corresponding points of {\bf p} on $(F_7+{\bf x}_4)\cap(F_3+{\bf x}_5)$ and $(F_8+{\bf x}_4)\cap(F_4+{\bf x}_3)$, as shown in FIG \ref{P8y2y3}. Similar to Subcase 3.1, we can further assume that ${\bf q}_1'$ and ${\bf q}_2'$ are both proper points in their corresponding edges.

\begin{figure}[h!]
\includegraphics[scale=0.75]{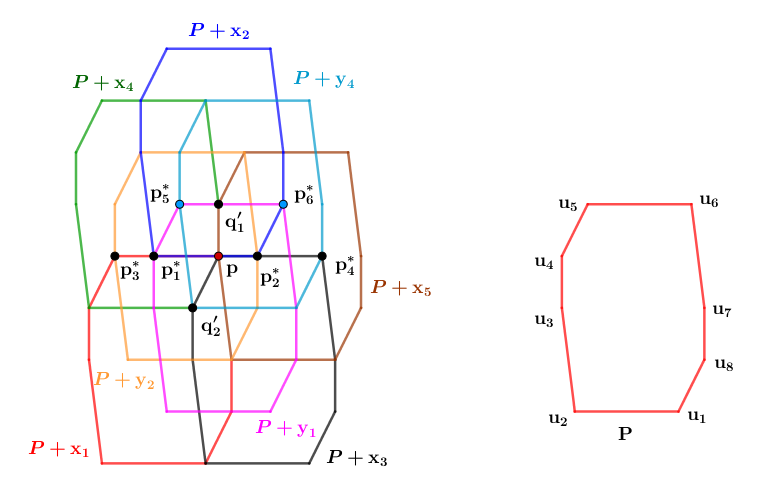}
\caption{${\bf y}_2={\bf y}_3$}\label{P8y2y3}
\end{figure}

By convexity of $P$ and ${\bf q}_1'={\bf p}_4^*+\widetilde{\bf g}_1+\widetilde{\bf g}_3$, we can get that ${\bf q}_1'\in int(P)+{\bf y}_4$. Since ${\bf y}_2={\bf y}_3$, we have $$u_2^3=u_2^7.$$ Similarly, we also have ${\bf q}_1'\in int(P)+{\bf y}_2$. Noticing that ${\bf q}_1'\in int(P)+{\bf x}_2$, we have $\varphi({\bf q}_1')\ge 3$, then we must have $$\varpi({\bf q}_1')=2,$$ otherwise it would contradict the assumption that $P+X$ is a fivefold translative tiling. Since ${\bf q}_1'={\bf p}+\widetilde{\bf g}_3$, apply Lemma 3.1 to ${\bf q}_1'$ and there is a point ${\bf y}_5\in X^{{\bf q}_1'}$ such that ${\bf p}\in int(P)+{\bf y}_5$. Since ${\bf q}_1'\in [int(P)+{\bf y}_4]\cap[int(P)+{\bf y}_2]$, we have $${\bf y}_5\notin\{{\bf y}_2,{\bf y}_4\},$$ then $${\bf y}_5={\bf y}_1$$ by $\varpi({\bf p})=2$ and $\varphi({\bf p})=3$. Next, we will show that $P+{\bf y}_1$ must be an F-type translate at ${\bf q}_1'$ and $F_5+{\bf y}_1$ has ${\bf q}_1'$ as its relative interior point.

\smallskip\noindent
{\bf Subcase 4.1.1.} If ${\bf q}_1'$ is in some edge $E_i+{\bf y}_1$, then we can deduce that ${\bf q}_2'$ is in one edge $E_k+{\bf y}_1$ by ${\bf p}\in\partial(P)+{\bf y}_1$. Specifically, $j=5$, $k=3$. On the other hand, it follows by ${\bf p}_4^*\in relint(E_8)+{\bf y}_4$ and ${\bf p}\in int(P)+{\bf y}_4$ that ${\bf q}_2'={\bf p}+\widetilde{\bf g}_8$ is also a proper point in $E_2+{\bf y}_4$. Together with ${\bf p}_2'\in relint(E_3)+{\bf y}_1$, we have $\varpi({\bf q}_2')\not=2$ by considering the dihedral adjacent wheels at ${\bf q}_2'$ and thus $$\varpi({\bf q}_2')\ge 3.$$
Let ${\bf q}_3'$ be the corresponding point of ${\bf q}_2'$ on $F_2+{\bf y}_1$. Since ${\bf q}_2'$ is a proper point in $E_3+{\bf y}_1$, we have ${\bf q}_3'\notin C(E)+X$.

Apply Lemma 3.1 to ${\bf q}_3'$, and there is a point ${\bf z}\in X^{{\bf q}_3'}$ such that $${\bf q}_2'\in int(P)+{\bf z}.$$ And it can also be deduced that $${\bf q}_3'\in int(P)+{\bf x}_1$$ by ${\bf p}_1^*\in relint(F_5)+{\bf x}_1$ and ${\bf q}_3'={\bf p}_1^*-\widetilde{\bf g}_3-\widetilde{\bf g}_2$. Thus $${\bf z}\not={\bf x}_1.$$ Since $u_2^3=u_2^7$, it can be shown that ${\bf q}_3'={\bf q}_2'-\widetilde{\bf g}_2\notin P+{\bf y}_2$ and therefore ${\bf z}\not={\bf y}_2$. Since ${\bf q}_2'={\bf p}+\widetilde{\bf g}_8\in int(P)+{\bf y}_2$ by ${\bf p}$ is the corresponding point of ${\bf p}_2^*$ in $(F_1+{\bf x}_2)\cap(F_5+{\bf x}_3)$, we have $$\varphi({\bf q}_2')\ge3$$ and consequently $$\varpi({\bf q}_2')+\varphi({\bf q}_2')\ge 6,$$ which contradicts the assumption.

\smallskip\noindent
{\bf Subcase 4.1.2.} If ${\bf q}_1'\in relint(F_i)+{\bf y}_1$, then we can deduce that $${\bf q}_1'\in relint(F_5)+{\bf y}_1,$$ since $\varpi({\bf q}_1')=2$ and $P+{\bf y}_1$ has only two facets $F_4+{\bf y}_1$ and $F_5+{\bf y}_1$ containing interior points of $P+{\bf x}_2$. And we have $$u_2^5-u_2^4=u_2^4-u_2^3.$$

By Lemma 4.3, we have that $(P+{\bf y}_1)\cap(P+{\bf x}_2)$ contains a zonotope with a part or the whole of $E_1+{\bf x}_2$ as one of its edges. Notice that ${\bf q}_1'\in int(P)+{\bf x}_2$ and we can choose a proper point ${\bf p}_6^*\in relint\left[(E_6+{\bf y}_1)\cap(E_8+{\bf x}_2)\right]$ and a point ${\bf p}_5^*={\bf p}_6^*+\widetilde{\bf g}_1$ such that $${\bf q}_1'\in {\bf p}_6^*+R_1-{\bf c}_{E_1^0}={\bf p}_5^*+R_1-{\bf c}_{E_2^0}.$$

Apparently, ${\bf p}_5^*\in int(P)+{\bf x}_2$. Applying Lemma 3.1 to ${\bf p}_5^*$, there is a point ${\bf y}_6\in X^{{\bf p}_5^*}$ such that $${\bf q}_1'\in int(P)+{\bf y}_6.$$ Then, ${\bf y}_6\not={\bf x}_2$. It is easy to see that ${\bf y}_2\notin X^{{\bf p}_5^*}$, then ${\bf y}_6\not={\bf y}_2$. If ${\bf p}_5^*\notin B(E)+{\bf y}_4$, then we have ${\bf y}_4\not={\bf y}_6$. Consequently, all ${\bf y}_2, {\bf y}_4, {\bf y}_6$ and ${\bf x}_2$ are pairwise distinct. Thus we have $$\varphi({\bf q}_1')\ge 4$$ and $$\varpi({\bf q}_1')+\varphi({\bf q}_1')\ge6,$$ which contradicts the assumption. Then we must have $${\bf p}_5^*\in B(E)+{\bf y}_4.$$ Furthermore, since the {\it y}-coordinate of ${\bf p}_5^*$ is equal to the {\it y}-coordinate of both ${\bf q}_1'$ and ${\bf u}_3+{\bf y}_4$, we have ${\bf p}_5^*\in E_3+{\bf y}_4$. By ${\bf q}_2'\in relint(E_2)+{\bf y}_4$ and ${\bf p}_5^*\in relint(E_3)+{\bf y}_4$, we have $$u_2^3-u_2^2=2(u_2^4-u_2^3).$$

For convenience, let ${\bf p}=(p_1, p_2, p_3)$, and let $w_1$, $w_2$, $w_3$ denote the {\it x}-coordinates of ${\bf u}_3+{\bf y}_4$, ${\bf p}_1^*$, ${\bf p}_5^*$, respectively. First, by computing the {\it x}-coordinate of ${\bf p}_4^*$ in two ways we get $$w_1+(u_1^7-u_1^6)+(u_1^6-u_1^5)+(u_1^5-u_1^4)=p_1+(u_1^6-u_1^5)$$ and thus $$w_1=p_1-(u_1^7-u_1^6)-(u_1^5-u_1^4).$$
On the other hand, since ${\bf y}_2={\bf y}_3$, by computing the difference of the {\it x}-coordinates between ${\bf p}_3^*$ and ${\bf p}_4^*$ in two ways we get $$(u_1^5-u_1^4)+(u_1^7-u_1^6)+(p_1-w_2)+(u_1^6-u_1^5)=2(u_1^6-u_1^5)$$ and thus $$w_2=p_1+(u_1^7-u_1^6)-(u_1^6-u_1^5)+(u_1^5-u_1^4).$$
Since ${\bf p}_5^*\in E_5+{\bf y}_1$ by ${\bf p}_5^*={\bf p}_6^*+\widetilde{\bf g}_1$, we get $$w_3=w_2+(u_1^5-u_1^4)=p_1+(u_1^7-u_1^6)-(u_1^6-u_1^5)+2(u_1^5-u_1^4).$$
Then, since $w_1=w_3$, we get $$2(u_1^7-u_1^6)+3(u_1^5-u_1^4)=u_1^6-u_1^5.$$

Thus, a zonotope $P$ with $\#B(E)=8$, $F_1$ vertical to {\it y}-axis, $F_3$ vertical to {\it x}-axis and ${\bf y}_2={\bf y}_3$ is a fivefold translative tile only if
\begin{equation*}
  \begin{cases}
  u_2^3=u_2^7,\\
  u_2^5-u_2^4=u_2^4-u_2^3,\\
  u_2^3-u_2^2=2(u_2^4-u_2^3,)\\
  u_1^6-u_1^5=2(u_1^7-u_1^6)+3(u_1^5-u_1^4).
  \end{cases}
\end{equation*}

Guaranteed by linear transformations, we can choose $u_2^4-u_2^3=1$, $u_1^6-u_1^5=2$ and $u_1^5-u_1^4=\alpha$ for $0<\alpha<\frac{2}{3}$. Noticing these conditions are exactly suitable to one of fivefold translative tiles in $\mathbb{E}^2$, see Subcase 4.1.2 of Lemma 3.8 of Yang and Zong \cite{yz2}. Thus in this case the projection of $P$ onto the plane {\it x}{\bf o}{\it y} must be a fivefold translative tile in $\mathbb{E}^2$.

\smallskip\noindent
{\bf Subcase 4.2.} ${\bf y}_3={\bf y}_4$. Let ${\bf q}_1'$ be the corresponding point of {\bf p} on $(F_7+{\bf x}_4)\cap(F_3+{\bf x}_5)$, as shown in FIG \ref{P8y3y4}. Similarly above, we can also assume that ${\bf q}_1'$ is a proper point in its corresponding edge.

\begin{figure}[h!]
\includegraphics[scale=0.7]{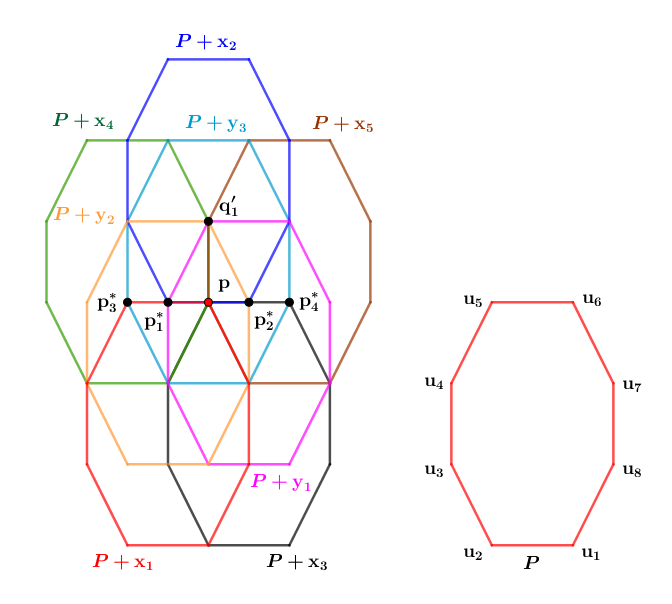}
\caption{${\bf y}_3={\bf y}_4$}\label{P8y3y4}
\end{figure}

By ${\bf p}_3^*\in relint(E_3)+{\bf y}_3$ and ${\bf p}_4^*\in relint(E_8)+{\bf y}_4$, we have $$u_2^3=u_2^8$$ and $$u_1^8-u_1^3=2(u_1^1-u_1^2).$$ Applying Lemma 3.1 to ${\bf q}_1'$, there is a point ${\bf y}_5\in X^{{\bf q}_1'}$ such that ${\bf p}\in int(P)+{\bf y}_5$. Since ${\bf q}_1'={\bf p}_3^*-\widetilde{\bf g}_1+\widetilde{\bf g}_3\in int(P)+{\bf y}_3$, we have ${\bf y}_5\not={\bf y}_3$.

If $\varpi({\bf q}_1')=3$ and ${\bf q}_1'\in int(F)$ holds for some $F\in B(E)+X$, by Subcase 3.1 we have $$\varpi({\bf q}_1')+\varphi({\bf q}_1')\ge6,$$ which contradicts the assumption. Thus we have either $\varpi({\bf q}_1')=2$ or $\varpi({\bf q}_1')=3$ and ${\bf q}_1'$ is on one edge $E^*\in A(E)+{\bf x}$ of P+{\bf x} for all ${\bf x}\in X^{{\bf q}_1'}$.

If both ${\bf y}_5\not={\bf y}_1$ and ${\bf y}_5\not={\bf y}_2$ holds simultaneously, then $$\varphi({\bf p})\ge 4$$ and therefore $$\varpi({\bf p})+\varphi({\bf p})\ge6,$$ which contradicts the assumption. So we must have either ${\bf y}_5={\bf y}_1$ or ${\bf y}_5={\bf y}_2$.

Suppose that ${\bf y}_5={\bf y}_1$. If ${\bf q}_1'$ is in one edge $E^*\in A(E)+{\bf y}_1$, then we have $$u_2^5-u_2^4=u_2^4-u_2^3.$$ If ${\bf q}_1'$ is a relative interior point in one facet $F\in B(E)+{\bf y}_1$, then by Subcase 3.1 we have $\varpi({\bf q}_1')=2$. By the structure of the dihedral adjacent wheel at ${\bf q}_1'$, one can deduce that ${\bf q}_1'$ must be an interior point of $F_5+{\bf y}_1$. Since $F_5+{\bf y}_1$ is vertical to the {\it y}-axis, we also obtain $$u_2^5-u_2^4=u_2^4-u_2^3.$$ If ${\bf y}_5={\bf y}_2$, we can also obtain these two equations above.

In conclusion, a zonotope $P$ with $\#B(E)=8$, $F_1$ vertical to {\it y}-axis, $F_3$ vertical to {\it x}-axis and ${\bf y}_3={\bf y}_4$ is a fivefold translative tile only if
\begin{equation*}
  \begin{cases}
  u_2^3=u_2^8,\\
  u_2^5-u_2^4=u_2^4-u_2^3,\\
  u_1^8-u_1^3=2(u_1^1-u_1^2).
  \end{cases}
\end{equation*}

Guaranteed by linear transformations, we can choose $u_2^4-u_2^3=2$, $u_1^1-u_1^2=2$ and $u_1^6=\beta$ for $0<\beta\le1$. Noticing these conditions are exactly suitable to one of fivefold translative tiles in $\mathbb{E}^2$, see Subcase 4.2 of Lemma 3.8 of Yang and Zong \cite{yz2}. Thus in this case the projection of $P$ onto the plane {\it x}{\bf o}{\it y} must be a fivefold translative tile in $\mathbb{E}^2$.

Therefore, this lemma is proved. \hfill{$\Box$}

\medskip\noindent
{\bf Lemma 5.4.} {\it If $P$ is a fivefold translative tile with $\#B(E)=10$, then the projection of $P$ along $E$ must be a fivefold translative tile in $\mathbb{E}^2$.}

\medskip\noindent
{\bf Proof.} Assume that $E'\in A(E)+X$ and let {\bf p} be a proper point in $E'$. By Lemma 3.3, we have $$\varphi({\bf p})\ge\left\lceil\frac{5-3}{2}\right\rceil=1.$$ By Lemma 3.2, we have $$\varpi({\bf p})=\kappa({\bf p})\cdot2+\ell({\bf p})\cdot\frac{1}{2},$$ where $\kappa({\bf p})$ is a positive integer and $\ell({\bf p})$ is a nonnegative integer. Notice that $\ell({\bf p})$ is the number of the facets in $B(E)+X$ which take {\bf p} as a relative interior point. Then we have the following two cases to discuss.

\smallskip\noindent
{\bf Case 1.} $\ell({\bf p})\not=0$ at one proper point ${\bf p}\in E'$ for some $E'\in A(E)+X$, so there  is one F-type translate $P+{\bf x}$ at {\bf p} for ${\bf x}\in X$. Suppose that ${\bf p}\in relint(F_i)+{\bf x}$ for some $i$, then there exist two points ${\bf p}_1^*\in(E_i+{\bf x})\setminus\{C(E)+X\}$ and ${\bf p}_2^*\in(E_{i+1}+{\bf x})\setminus\{C(E)+X\}$ such that $${\bf p}\in {\bf p}_1^*+R_i-{\bf c}_{E_i^0}={\bf p}_2^*+R_i-{\bf c}_{E_{i+1}^0}.$$ Applying Lemma 3.1 to these two points ${\bf p}_1^*, {\bf p}_2^*$, respectively, there are two different points ${\bf z}_1\in X^{{\bf p}_1^*}$ and ${\bf z}_2\in X^{{\bf p}_2^*}$ such that $${\bf p}\in \left[int(P)+{\bf z}_1\right]\cap\left[int(P)+{\bf z}_2\right].$$ Then we have $\varphi({\bf p})\ge2$. If $\varpi({\bf p})\ge4$, we can deduce that $$\varpi({\bf p})+\varphi({\bf p})\ge6,$$ which contradicts the assumption.

If $\varpi({\bf p})=3$, then $\varpi({\bf p})=\kappa({\bf p})\cdot2+\ell({\bf p})\cdot\frac{1}{2}$ has one group of solution $\kappa({\bf p})=2$, $\ell({\bf p})=2$ and it follows that $P+X^{\bf p}$ consists of five E-type translates and two E-type translates and there are at most two dihedral adjacent wheels in $P+X^{\bf p}$. Without loss of generality, assume that $E'=E_1$. For convenience, let $X^{\bf p}=\{{\bf x}_1, {\bf x}_2,...,{\bf x}_7\}$ with ${\bf x}_1={\bf o}$ and we still use the projections of translates at {\bf p} onto the plane {\it x}{\bf o}{\it y} to discuss the following cases.

\smallskip\noindent
{\bf Subcase 1.1.} There is only one dihedral adjacent wheel $P+X^{\bf p}$ at {\bf p} consisting of five E-type translates and two F-type translates, and $P+{\bf x}_1$,..., $P+{\bf x}_7$ are sequential in $P+X^{\bf p}$. Without loss of generality, this wheel has the following two kinds of structure.

\smallskip\noindent
{\bf Subcase 1.1.1.} As FIG \ref{P10_7translates1} shows, we assume that $P+{\bf x}_2$ and $P+{\bf x}_4$ are F-type translates and others are E-type translates.

\begin{figure}[h!]
 \subfigure[]
 {\begin{minipage}{0.85\linewidth}
 \includegraphics[scale=0.55]{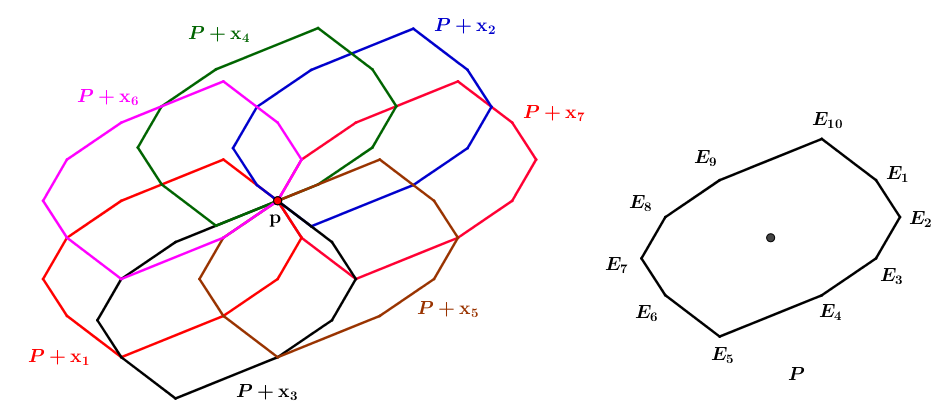}
 \end{minipage}}
 \subfigure[]
 {\begin{minipage}{0.55\linewidth}
  \includegraphics[scale=0.7]{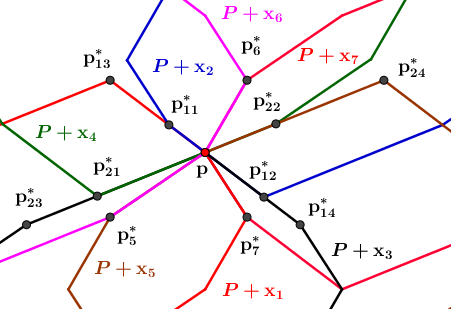}
  \end{minipage}}
  \caption{The dihedral adjacent wheel $P+X^{\bf p}$ in Subcase 1.1.1.}
  \label {P10_7translates1}
\end{figure}

Moreover, let ${\bf p}_{11}^*$, ${\bf p}_{12}^*$, ${\bf p}_{13}^*$ and ${\bf p}_{14}^*$ be the corresponding points of {\bf p} in $(F_5+{\bf x}_2)\cap(F_{10}+{\bf x}_1)$, $(F_5+{\bf x}_2)\cap(F_{10}+{\bf x}_3)$, $F_{10}+{\bf x}_1$ and $F_{10}+{\bf x}_3$, let ${\bf p}_{21}^*$, ${\bf p}_{22}^*$, ${\bf p}_{23}^*$ and ${\bf p}_{24}^*$ be the corresponding points of {\bf p} in $(F_4+{\bf x}_4)\cap(F_9+{\bf x}_3)$, $(F_4+{\bf x}_4)\cap(F_9+{\bf x}_5)$, $F_9+{\bf x}_3$ and $F_9+{\bf x}_5$, and let ${\bf p}_5^*$, ${\bf p}_6^*$ and ${\bf p}_7^*$ be the corresponding points of {\bf p} in $F_8+{\bf x}_5$, $F_2+{\bf x}_6$ and $F_6+{\bf x}_7$. Since {\bf p} is a proper point, we have that all these corresponding points are not in $C(E)+X$. Applying Lemma 3.1 to all ${\bf p}_{ij}^*$ and ${\bf p}_k^*$, there are ${\bf y}_{ij}\in X^{{\bf p}_{ij}^*}$ and ${\bf y}_k\in X^{{\bf p}_k^*}$ satisfying that $${\bf p}\in int(P)+{\bf y}_{ij}$$ and $${\bf p}\in int(P)+{\bf y}_k$$ for all $i=1,2$, $j=1,2,3,4$ and $k=5,6,7$. Since $\varpi({\bf p})=3$, we need $\varphi({\bf p})=2$, in other words, there are two I-type translates $P+{\bf z}_1$ and $P+{\bf z}_2$ at {\bf p}, then ${\bf y}_{i1}={\bf y}_{i4}$ and ${\bf y}_{i2}={\bf y}_{i3}$ by $\#B(E)=10$ and Lemma 3.1 for $i=1,2$. Notice that if one of ${\bf y}_5={\bf y}_6$ and ${\bf y}_7={\bf y}_6$ holds, it is easy to see that $P+{\bf y}_6$ has a belt $B(E)+{\bf y}_6$ with $\#B(E)=6$ contradicting $\#B(E)=10$. Without loss of generality, let ${\bf y}_{i1}={\bf z}_1$ for all $i$ and let ${\bf y}_k={\bf z}_1$ for $k=5,7$ and ${\bf y}_6={\bf z}_2$. By observing the two pairwise points ${\bf p}_{11}^*$ and ${\bf p}_7^*$, ${\bf p}_5^*$ and ${\bf p}_{24}^*$, we can see that $(P+{\bf x}_1)\cap(P+{\bf z}_1)$ has a blet $B(E,{\bf x}_1,{\bf z}_1)$ with $\#B(E,{\bf x}_1,{\bf z}_1)=4$, $(P+{\bf x}_5)\cap(P+{\bf z}_1)$ has a blet $B(E,{\bf x}_5,{\bf z}_1)$ with $\#B(E,{\bf x}_5,{\bf z}_1)=4$, then we need ${\bf p}_5^*\in F_4+{\bf z}_1$ which is impossible since $B(E,{\bf x}_1,{\bf z}_1)\cap(P+{\bf x}_1)\subseteq(F_5+{\bf z}_1)\cup(F_6+{\bf z}_1)$ and $[int(P)+{\bf x}_1]\cap(F_4+{\bf z}_1)=\emptyset$.

\smallskip\noindent
{\bf Subcase 1.1.2.} As FIG \ref{P10_7translates2} shows, we assume that $P+{\bf x}_2$ and $P+{\bf x}_6$ are F-type translates and others are E-type translates.

\begin{figure}[h!]
 \subfigure[]
 {\begin{minipage}{0.95\linewidth}
 \includegraphics[scale=0.6]{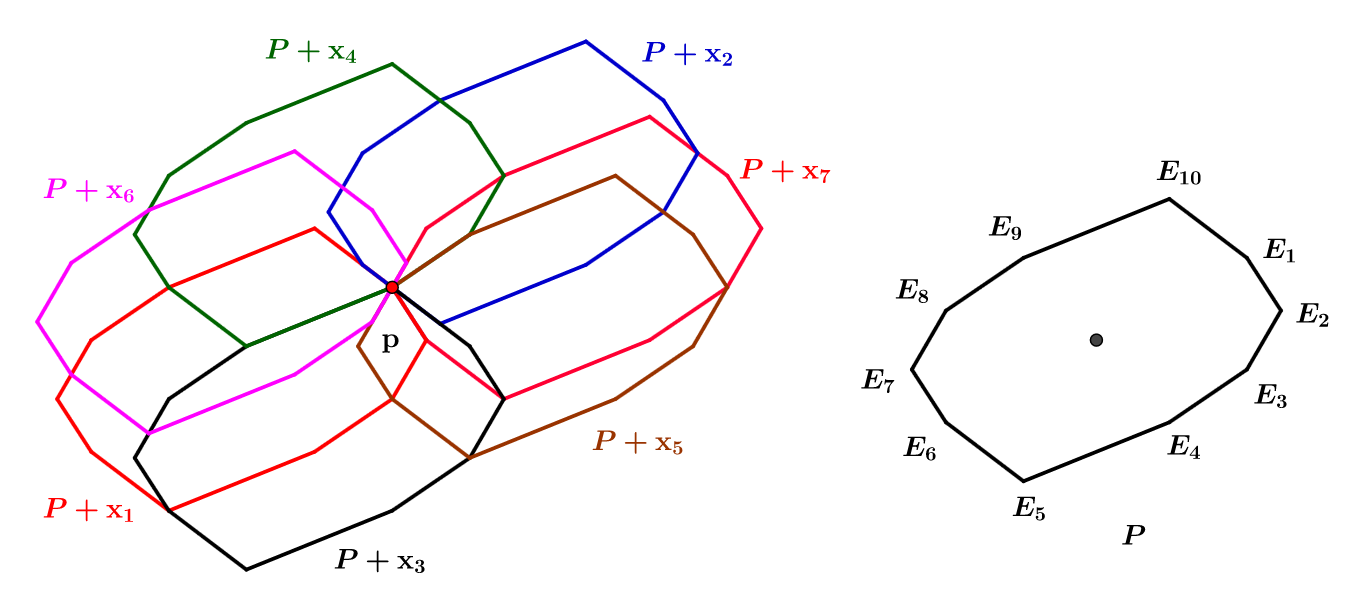}
 \end{minipage}}
 \subfigure[]
 {\begin{minipage}{0.6\linewidth}
 \vspace{0.5cm}
 \includegraphics[scale=0.7]{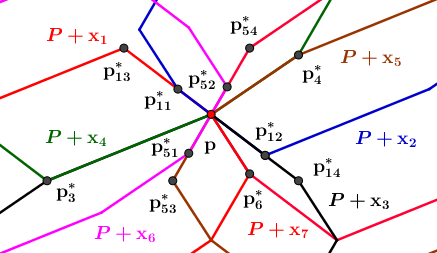}
 \end{minipage}}
 \caption{The dihedral adjacent wheel $P+X^{\bf p}$ in Subcase 1.1.2.}
 \label {P10_7translates2}
\end{figure}

Moreover, let ${\bf p}_{11}^*$, ${\bf p}_{12}^*$, ${\bf p}_{13}^*$ and ${\bf p}_{14}^*$ be the corresponding points of {\bf p} in $(F_5+{\bf x}_2)\cap(F_{10}+{\bf x}_1)$, $(F_5+{\bf x}_2)\cap(F_{10}+{\bf x}_3)$, $F_{10}+{\bf x}_1$ and $F_{10}+{\bf x}_3$, let ${\bf p}_3^*$, ${\bf p}_4^*$ and ${\bf p}_6^*$ be the corresponding points of {\bf p} in $F_9+{\bf x}_3$, $F_3+{\bf x}_4$ and $F_6+{\bf x}_7$ and let ${\bf p}_{51}^*$, ${\bf p}_{52}^*$, ${\bf p}_{53}^*$ and ${\bf p}_{54}^*$ be the corresponding points of {\bf p} in $(F_2+{\bf x}_6)\cap(F_7+{\bf x}_5)$, $(F_2+{\bf x}_6)\cap(F_7+{\bf x}_7)$, $F_7+{\bf x}_5$ and $F_7+{\bf x}_7$. Since {\bf p} is a proper point, we have that all these corresponding points are not in $C(E)+X$. Applying Lemma 3.1 to all ${\bf p}_{ij}^*$ and ${\bf p}_k^*$, there are ${\bf y}_{ij}\in X^{{\bf p}_{ij}^*}$ and ${\bf y}_k\in X^{{\bf p}_k^*}$ satisfying that $${\bf p}\in int(P)+{\bf y}_{ij}$$ and $${\bf p}\in int(P)+{\bf y}_k$$ for all $i=1,5$, $j=1,2,3,4$ and $k=3,4,6$. Since $\varpi({\bf p})=3$, we need $\varphi({\bf p})=2$, in other words, there are two I-type translates $P+{\bf z}_1$ and $P+{\bf z}_2$ at {\bf p}, then ${\bf y}_{i1}={\bf y}_{i4}$ and ${\bf y}_{i2}={\bf y}_{i3}$ for $i=1,2$. Without loss of generality, let ${\bf y}_{i1}={\bf z}_1$ for $i=1,5$ and let ${\bf y}_6={\bf z}_1$. By observing the two pairwise points ${\bf p}_6^*$ and ${\bf p}_{11}^*$, ${\bf p}_6^*$ and ${\bf p}_{54}^*$, we can see that $(P+{\bf x}_1)\cap(P+{\bf z}_1)$ has a blet $B(E,{\bf x}_1,{\bf z}_1)$ with $\#B(E,{\bf x}_1,{\bf z}_1)=4$, $(P+{\bf x}_7)\cap(P+{\bf z}_1)$ has a blet $B(E,{\bf x}_7,{\bf z}_1)$ with $\#B(E,{\bf x}_7,{\bf z}_1)=4$, then it follows that $P+{\bf z}_1$ has a belt $B(E)+{\bf z}_1$ with $\#B(E)=6$ contradicting $\#B(E)=10$.

\smallskip\noindent
{\bf Subcase 1.2.} There are two dihedral adjacent wheels $P+X_1^{\bf p}$ and $P+X_2^{\bf p}$ at {\bf p} with $X_1^{\bf p}=\{{\bf x}_1, {\bf x}_2,...,{\bf x}_5\}$, $X_2^{\bf p}=\{{\bf x}_6,{\bf x}_7\}$ where $E'=E_1$ and ${\bf x}_1={\bf o}$. Without loss of generality, we can assume that both $F_4+{\bf x}_6$ and $F_9+{\bf x}_7$ contain {\bf p} as their relative interior point, as shown in FIG \ref{P10_2wheels}.

\begin{figure}[h!]
 \subfigure[]
 {\begin{minipage}{0.95\linewidth}
 \includegraphics[scale=0.65]{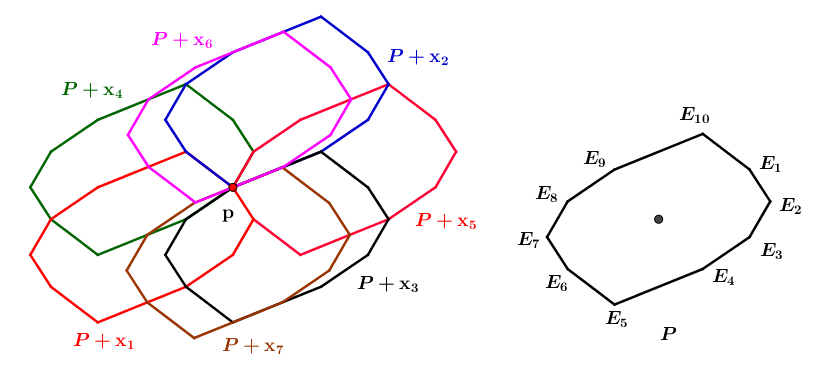}
 \end{minipage}}
 \subfigure[]
 {\begin{minipage}{0.6\linewidth}
 \vspace{0.5cm}
 \includegraphics[scale=0.7]{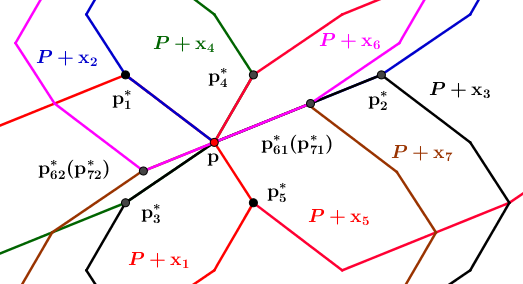}
 \end{minipage}}
 \caption{The two dihedral adjacent wheels at {\bf p} in Subcase 1.2.}
 \label {P10_2wheels}
\end{figure}

Clearly, $P+X_1^{\bf p}$ consists of five E-type translates with $P+{\bf x}_1$,..., $P+{\bf x}_5$ in sequence and $P+X_2^{\bf p}$ consists of two F-type translates.
Let ${\bf p}_1^*$, ${\bf p}_2^*$, ${\bf p}_3^*$, ${\bf p}_4^*$ and ${\bf p}_5^*$ be the corresponding points of {\bf p} in $F_{10}+{\bf x}_1$, $F_4+{\bf x}_2$, $F_8+{\bf x}_3$, $F_2+{\bf x}_4$ and $F_6+{\bf x}_5$. For convenience, write $(F_4+{\bf x}_2)\cap(F_9+{\bf x}_3)$ as $F_4^*$. Let ${\bf p}_{61}^*$ be the corresponding point of {\bf p} in $(F_4+{\bf x}_6)\cap F_4^*$ and ${\bf p}_{62}^*={\bf p}_{61}+\widetilde{\bf g}_4$, and let ${\bf p}_{71}^*$ be the corresponding point of {\bf p} in $(F_9+{\bf x}_7)\cap F_4^*$ and ${\bf p}_{72}^*={\bf p}_{71}+\widetilde{\bf g}_4$. Moreover, we can aslo assume that all ${\bf p}_k^*$ and ${\bf p}_{ij}^*$ are proper points.
Applying Lemma 3.1 to all ${\bf p}_{ij}^*$ and ${\bf p}_k^*$, there are points ${\bf y}_{ij}\in X^{{\bf p}_{ij}^*}$ and ${\bf y}_k\in X^{{\bf p}_k^*}$ satisfying that $${\bf p}\in int(P)+{\bf y}_{ij}$$ and $${\bf p}\in int(P)+{\bf y}_k$$ for $i=6,7$, $j=1,2$ and $k=1,2,...,5$. Clearly, ${\bf y}_{i1}\not={\bf y}_2$. Since $\varpi({\bf p})=3$, we need $\varphi({\bf p})=2$, in other words, there are two I-type translates $P+{\bf z}_1$ and $P+{\bf z}_2$ at {\bf p} and $${\bf y}_{61}={\bf y}_{71}.$$ Since ${\bf p}_{61}^*$ is the corresponding point of {\bf p} in $(F_4+{\bf x}_6)\cap F_4^*$, we have $${\bf p}_{61}^*\in [relint(F_4)+{\bf x}_2]\cap[relint(F_9)+{\bf x}_3],$$ and both $P+{\bf x}_2$ and $P+{\bf x}_3$ are the F-type translates at ${\bf p}_{61}^*$. According to the discussion above, we obtain that $\varpi({\bf p}_{61}^*)=3$ and $P+X^{{\bf p}_{61}^*}$ can be divided into two dihedral adjacent wheels $P+X_1^{{\bf p}_{61}^*}$ and $P+X_2^{{\bf p}_{61}^*}$. Without loss of generality, let $X_2^{{\bf p}_{61}^*}=\{{\bf x}_2,{\bf x}_3\}$. Apparently, $P+X_1^{{\bf p}_{61}^*}$ consists of five E-type translates and ${\bf z}_1\in X_1^{{\bf p}_{61}^*}$. Since the structure of $P+X_1^{{\bf p}_{61}^*}$ can be determined by $P+{\bf x}_6$, we can see that $(P+{\bf z}_1)\cap(P+{\bf x}_6)$ has a belt $B(E,{\bf z}_1, {\bf x}_6)$ with $\#B(E,{\bf z}_1, {\bf x}_6)=6$ and $(P+{\bf z}_1)\cap(P+{\bf x}_7)$ has a belt $B(E,{\bf z}_1, {\bf x}_7)$ with $\#B(E,{\bf z}_1, {\bf x}_7)=6$. It is easy to get that ${\bf p}_1^*\in int(P)+{\bf z}_1$ and ${\bf p}_3^*\in int(P)+{\bf z}_1$, then we have $${\bf y}_i={\bf z}_2$$ for $i=1,2,3$. By the structure of $(P+{\bf z}_2)\cap(P+{\bf x}_2)$  and $(P+{\bf z}_2)\cap(P+{\bf x}_2)$, we get $\#B(E)=6$ contradicting $\#B(E)=10$.

\smallskip\noindent
{\bf Case 2.} $\ell({\bf p})=0$ for all proper points ${\bf p}$ in any $E'\in A(E)+X$. Then it follows by $\varpi({\bf p})=\kappa({\bf p})\cdot2+\ell({\bf p})\cdot\frac{1}{2}$ that $\varpi({\bf p})=4$ or $\varpi({\bf p})=2$.

\smallskip\noindent
{\bf Subcase 2.1.} $\varpi({\bf p})=4$ at one proper point ${\bf p}\in E'$ for some $E'\in A(E)+X$, then the translates in $P+X^{\bf p}$ can be divided into two dihedral adjacent wheels each of which contains five E-type translates. Without loss of generality, assume that $E'=E_1+{\bf x}_1$ with ${\bf x}_1={\bf o}$ and $P+{\bf x}_1$, $P+{\bf x_2}$,..., $P+{\bf x}_5$ in sequence generate such a wheel at {\bf p}. Project $P+{\bf x}_i$ onto the plane {\it x}{\bf o}{\it y}, and we also use the $P+{\bf x}_i$ to denote its projection without confusion, as shown in FIG \ref{P10_l=0}. Let ${\bf p}_1^*$ and ${\bf p}_2^*$ be the corresponding points of {\bf p} on $(F_1+{\bf x}_1)\cap(F_6+{\bf x}_5)$ and $(F_{10}+{\bf x}_1)\cap(F_5+{\bf x}_2)$, respectively. Then, applying Lemma 3.1 to each ${\bf p}_i^*$, there is ${\bf y}_i\in X^{{\bf p}_i^*}$ such that $${\bf p}\in int(P)+{\bf y}_i,$$ for $i=1,2$.

\begin{figure}[h!]
\includegraphics[scale=0.62]{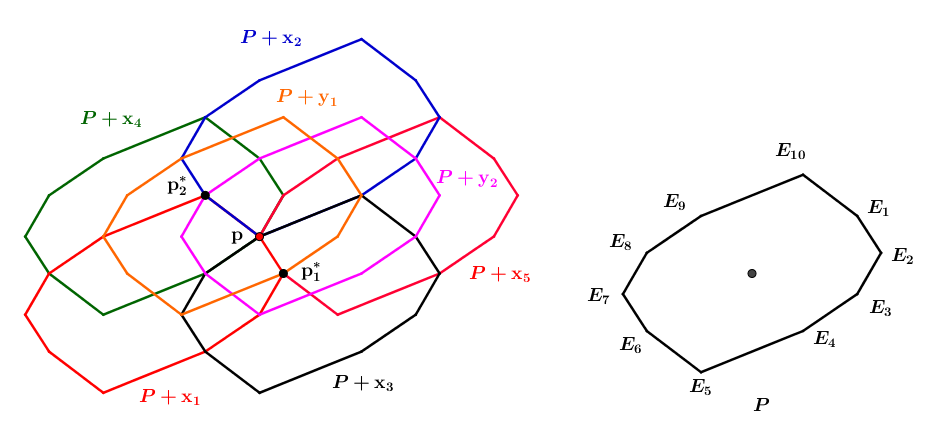}
\caption{One dihedral adjacent wheel at {\bf p} in Subcase 2.1}\label{P10_l=0}
\end{figure}

In fact, by the previous analysis, we have $${\bf y}_1=-\widetilde{\bf g}_3-\widetilde{\bf g}_2+\lambda\vec{E}$$ for $|\lambda|<1$ and $(P+{\bf x}_1)\cap(P+{\bf y}_1)$ has a belt $B(E,{\bf x}_1,{\bf y}_1)$ with $\#B(E,{\bf x}_1,{\bf y}_1)=6$. Since ${\bf p}_2^*={\bf p}+\widetilde{\bf g}_5={\bf p}_1^*-\widetilde{\bf g}_1+\widetilde{\bf g}_5$, we have $${\bf p}_2^*\in int(P)+{\bf y}_1.$$ Then ${\bf y}_1\not={\bf y}_2$ and we have $$\varphi({\bf p})\ge 2$$ and $$\varpi({\bf p})+\varphi({\bf p})\ge6,$$ which contradicts the assumption.

\smallskip\noindent
{\bf Subcase 2.2.} $\varpi({\bf p})=2$ holds for all proper points ${\bf p}$ in each $E'\in A(E)+X$.

\smallskip
For convenience, we still concentrate on the structure of $P+X^{\bf p}$ as $P+X_1^{\bf p}$ and the corresponding points ${\bf p}_1^*$,..., ${\bf p}_5^*$ of {\bf p} in Subcase 1.2 above, shown in FIG \ref{P10_2wheels}. Since $\varpi({\bf p})=2$, there must only three I-type translates at {\bf p}, denoted by $P+{\bf z}_1$, $P+{\bf z}_2$, $P+{\bf z}_3$, respectively. Notice that the intersection of any two E-type translates $P+{\bf x}$ and $P+{\bf y}$ at {\bf p} is either a centrally symmetric polygon or a centrally symmetric polytope having a belt $B(E,{\bf x},{\bf y})$ with $\#B(E,{\bf x},{\bf y})=6$. Then by $\varpi({\bf p}_i^*)=2$ and the structure of the dihedral adjacent wheel at ${\bf p}_i^*$ determined by $P+{\bf x}_i$, we have ${\bf y}_i\not={\bf y}_{i+1}$ for ${\bf y}_6={\bf y}_1$ and the following five cases.

\begin{enumerate}
  \item ${\bf y}_1={\bf y}_3$ and ${\bf y}_2={\bf y}_4$;
  \item ${\bf y}_1={\bf y}_3$ and ${\bf y}_2={\bf y}_5$;
  \item ${\bf y}_1={\bf y}_4$ and ${\bf y}_2={\bf y}_5$;
  \item ${\bf y}_1={\bf y}_4$ and ${\bf y}_3={\bf y}_5$ and
  \item ${\bf y}_2={\bf y}_4$ and ${\bf y}_3={\bf y}_5$.
\end{enumerate}

 Without loss of generality, we take (2) above as an example for convenience and the projections of translates at {\bf p} onto the plane {\it x}{\bf o}{\it y} is shown in FIG \ref{P10translates}.

\begin{figure}[h!]
\includegraphics[scale=0.65]{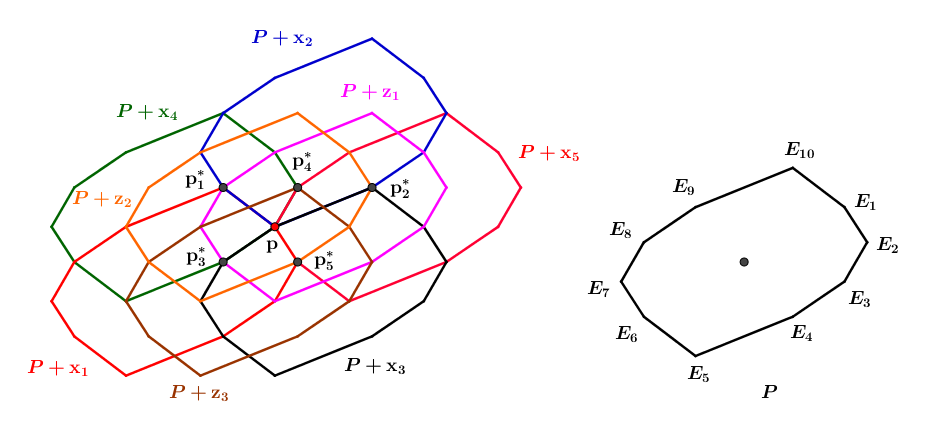}
\caption{The translates at {\bf p} for $\#B(E)=10$ in (2)}\label{P10translates}
\end{figure}

To show that the projection of $P$ onto the plane {\it x}{\bf o}{\it y} is the fivefold translative tile in $\mathbb{E}^2$, we have the following two cases.

\smallskip\noindent
{\bf Subcase 2.2.1.} $P+X$ is a non-proper fivefold translative tiling relate to $E$. For convenience, without loss of generality, we can assume that $E_1\not\subset \partial(P)+{\bf x}_2$ and $relint(E_1)\cap(\partial(P)+{\bf x}_2)\not=\emptyset$, specifically, that is $E_1\not=E_5+{\bf x}_2$, then we have the following claims. Notice that ${\bf x}_1={\bf o}$.

\medskip\noindent
{\bf Claim 5.1.} For the proper point {\bf q}$\in relint(E_1)\setminus(E_5+{\bf x}_2)$, we have ${\bf x}_2+\vec{E}\in X$ or ${\bf x}_2-\vec{E}\in X$.

\medskip\noindent
{\bf Claim 5.2.} ${\bf x}_i+k\vec{E}\in X, {\bf z}_j+k\vec{E}\in X$ for $i=1,2,3,4,5$, $j=1,2,3$ and $k\in\mathbb{Z}$.

\medskip
The proof of Claim 5.1 will be shown at the end of this section. Claim 5.2 can be easily deduced following Claim 5.1 by $E_1+k\vec{E}\not=E_5+{\bf x}_2+k\vec{E}$ for $k\in\mathbb{Z}$. Then we continue the discussion about Subcase 2.2.1.

If $P$ is a prism with $\#B(E)=10$, we can conclude by Claims 5.1-5.2 that the translates of $\bigcup_{k\in\mathbb{Z}}(P+k\vec{E})$ consisting of all translates in $P+X$ generate a fivefold translative tiling of $\mathbb{E}^3$. Apparently, the projection of $\bigcup_{k\in\mathbb{Z}}(P+k\vec{E})$ onto the plane {\it x}{\bf o}{\it y}, identical to the projection of $P$ onto the plane {\it x}{\bf o}{\it y}, is the fivefold translative tile in $\mathbb{E}^2$.

If $P$ is not a prism with $\#B(E)=10$, we have that $int[P\cap(P+\vec{E})]\not=\emptyset$ by Corollary 4.2. Without loss of generality, we can assume that the upper vertex {\bf v} of $E_1+{\bf x}_1$ is in $E_5+{\bf x}_2$. By Claim 5.1, we have ${\bf x}_i+\vec{E}\in X$ and ${\bf z}_j+\vec{E}\in X$ for $i=1,2,...,5$ and $j=1,2,3$ where the direction of $\vec{E}$ is the positive direction of {\it z}-axis. Let $N({\bf v},\epsilon)$ be a suitable neighborhood of {\bf v} for some small positive number $\epsilon$, then we have following steps to analyze the non-proper fivefold translative tiling relate to $E$ for $\#B(E)=10$.

\smallskip\noindent
{\bf Step 1.} $int[(P+{\bf x}_1)\cap(P+{\bf x}_1+\vec{E})]\cap N({\bf v}, \epsilon)\not=\emptyset$ by $int[(P+{\bf x}_1)\cap(P+{\bf x}_1+\vec{E})]\not=\emptyset$. Let $$C({\bf v},\epsilon)=int[(P+{\bf x}_1)\cap(P+{\bf x}_1+\vec{E})]\cap N({\bf v}, \epsilon)$$ and proceed to the next step.

\smallskip\noindent
{\bf Step 2.} $int[(P+{\bf x}_4)\cup(P+{\bf x}_4+\vec{E})]\cap C({\bf v},\epsilon)\not=\emptyset$ by ${\bf v}\in relint[(E_3+{\bf x}_4)\cup(E_3+{\bf x}_4+\vec{E})]$. Replace $$int[(P+{\bf x}_4)\cup(P+{\bf x}_4+\vec{E})]\cap C({\bf v}, \epsilon)$$ with $C({\bf v},\epsilon)$ and proceed to the next step.

\smallskip\noindent
{\bf Step 3.} $C({\bf v},\epsilon)\subset int[(P+{\bf z}_j)\cup(P+{\bf z}_j+\vec{E})]$ by ${\bf v}\in int[(P+{\bf z}_j)\cup(P+{\bf z}_j+\vec{E})]$ for small enough positive number $\epsilon$ and $j=1,2,3$.

\smallskip
According to these steps above, we can conclude that there exists a point ${\bf p}'\in C({\bf v},\epsilon)$ satisfying $$\varphi({\bf p}')\ge 6,$$ which contradicts that $P+X$ is a fivefold translative tiling of $\mathbb{E}^3$. Thus, $P$ must be a prism with $\#B(E)=10$ for the non-proper fivefold translative tiling $P+X$ relate to $E$. As proved above, we can conclude the projection of $P$ onto the plane {\it x}{\bf o}{\it y} is a fivefold translative tile in $\mathbb{E}^2$.

\smallskip\noindent
{\bf Subcase 2.2.2.} $P+X$ is a proper fivefold translative tiling relate to $E$. We mainly see the claims following. And for convenience, proofs of these claims are also shown at the end of this section.

\medskip\noindent
{\bf Claim 5.3.} The points in any $E'\in A(E)+X$ except for two vertices of $E'$ are all proper points.

\medskip\noindent
{\bf Claim 5.4.} $\sum_{i=1}^5 k_i\cdot {\bf g}_i\in X$, where ${\bf g}_i$ is the translation vector of $F_i$ in $P$ and $k_i\in\mathbb{Z}$.

\medskip\noindent
{\bf Claim 5.5.} Let $X'=\{{\bf x}: {\bf x}=\sum_{i=1}^5 k_i\cdot {\bf g}_i\in X, k_i\in\mathbb{Z}\}$, then the projection of $P+X'$ onto the plane {\it x}{\bf o}{\it y} is a fivefold translative tiling in $\mathbb{E}^2$.

\medskip
By Claims 5.3-5.5, we also have that the projection of $P$ onto the plane {\it x}{\bf o}{\it y} is a fivefold translative tile in $\mathbb{E}^2$.

Thus, the proof of this lemma has been finished. \hfill{$\Box$}

\medskip\noindent
{\bf Proof of Claim 5.1.} Since $\varpi({\bf q})=2$, the structure of the dihedral adjacent wheel at {\bf q} is also determined by $P+{\bf x}_1$. Without loss of generality, denote the translates by $P+{\bf x}_1', P+{\bf x}_2', P+{\bf x}_3', P+{\bf x}_4', P+{\bf x}_5'$ in sequence. Certainly, they are all E-type translates at {\bf q} and  ${\bf x}_1={\bf x}_1'$. Naturally, $${\bf x}_2'={\bf x}_2+\lambda\vec{E}$$ for a real number $\lambda$, then we conclude $\lambda=\pm1$.
\begin{enumerate}
  \item If $|\lambda|<1$, then $$relint[(E_5+{\bf x}_2)\cap(E_5+{\bf x}_2')]\not=\emptyset,$$
 thus there exists a proper point {\bf q}$'\in relint[(E_5+{\bf x}_2)\cap(E_5+{\bf x}_2')]$. Then both $P+{\bf x}_2$ and $P+{\bf x}_2'$ are the E-type translates at ${\bf q}'$ and we have $\varpi({\bf q}')\ge3$ which contradicts $\varpi({\bf q}')=2$.
  \item If $|\lambda|>1$, then there exists a proper point {\bf q}$''\in relint\left\{E_1\setminus[(E_5+{\bf x}_2)\cup(E_5+{\bf x}_2')]\right\}$. Since $\varpi({\bf q}'')=2$, the structure of the dihedral adjacent wheel at {\bf q}$''$ is determined by $P+{\bf x}_1$. Without loss of generality, denoted by $P+{\bf x}_1'', P+{\bf x}_2'', P+{\bf x}_3'', P+{\bf x}_4'', P+{\bf x}_5''$ in sequence and ${\bf x}_1''={\bf x}_1$, we must have $$relint[(E_5+{\bf x}_2)\cap(E_5+{\bf x}_2'')]\not=\emptyset$$ by the finiteness of the length of $E_1$, then $\mid\lambda\mid>1$ is also impossible as in (1) above.
\end{enumerate}

Consequently, $\lambda=\pm1$. Thus the proof of Claim 5.1 is finished.  \hfill{$\Box$}

\medskip\noindent
{\bf Proof of Claim 5.3.} Since $P+X$ is a proper fivefold translative tiling relate to $E$, for any $E'\in A(E)+X$ and any $P'\in P+X$, we have that $E'\subset P'$ or $relint(E')\cap P'=\emptyset$, then for all points ${\bf p}\in relint(E')$, we have $${\bf p}\notin C(E)+X.$$ Then by the definition of the proper point, we need to show that every corresponding point ${\bf p}^*$ of {\bf p} in its corresponding facet $F\in B(E)+X$ is not in $C(E)+X$. Without loss of generality, we can take $E'=E_1$ and ${\bf p}\in relint(E_1)$. And let the structure of the dihedral adjacent wheel $P+X^{\bf p}$ at {\bf p} and the corresponding points ${\bf p}_1^*$,...,${\bf p}_5^*$ of {\bf p} be shown in FIG \ref{P10translates}, where $X^{\bf p}=\{{\bf x}_1,...,{\bf x}_5\}$ and ${\bf x}_1={\bf o}$. Since ${\bf p}_1^*={\bf p}+\widetilde{\bf g}_5$, ${\bf p}_2^*={\bf p}-\widetilde{\bf g}_4$, ${\bf p}_3^*={\bf p}+\widetilde{\bf g}_3$, ${\bf p}_4^*={\bf p}-\widetilde{\bf g}_2$ and ${\bf p}_5^*={\bf p}+\widetilde{\bf g}_1$, we have ${\bf p}_1^*\in relint(E_{10})+{\bf x}_1$, ${\bf p}_2^*\in relint(E_4)+{\bf x}_2$, ${\bf p}_3^*\in relint(E_8)+{\bf x}_3$, ${\bf p}_4^*\in relint(E_2)+{\bf x}_4$, and ${\bf p}_5^*\in relint(E_6)+{\bf x}_5$, then $${\bf p}_i^*\notin C(E)+X.$$ Notice that there are only these five corresponding points of {\bf p}, and we conclude that {\bf p} is a proper point in $relint(E_1)$, further this proof is finished.  \hfill{$\Box$}

\medskip\noindent
{\bf Proof of Claim 5.4.} Since $\varpi({\bf p})=2$ and $\varphi({\bf p})=3$ holds for all relative interior points {\bf p} in $E'\in A(E)+X$, there is only one dihedral adjacent wheel at {\bf p} consisting of five E-type translates. For convenience, we call the translates at {\bf p} also the translates at $E$. Notice that ${\bf x}_1={\bf o}\in X$. Without loss of generality, let us observe the $P+{\bf x}_1$ as shown in FIG \ref{P10translates}. It is easy to see that $P-{\bf g}_i$ and $P+{\bf g}_i$ are the E-type translates at $E_i-{\bf g}_i$ and $E_i+{\bf g}_i$, respectively, for $i=1,2,3,4,5$, and we have $$\{\pm{\bf g}_1,\pm{\bf g}_2,...,\pm{\bf g}_5\}\subset X.$$ Similarly, $P+k_i{\bf g}_i$ is an E-type translate at $E_i+k_i{\bf g}_i$ for all $i$, then we have $$k_i{\bf g}_i\in X$$ for $k_i\in\mathbb{Z}$ and $i=1,2,3,4,5$. It follows by recursion that $P+\sum_{i=1}^5 k_i\cdot {\bf g}_i$ must be an E-type translate at some $E'\in A(E)+X$. In fact, if $k_j$ is the last non-zero coefficient, equivalently, $k_j\not=0$ but $k_{j+1}=\cdots=k_5=0$ for some $1\le j\le5$, then $P+\sum_{i=1}^5 k_i\cdot {\bf g}_i$ is an E-type translate at $E_j+\sum_{i=1}^5 k_i\cdot {\bf g}_i$ and we have $$\sum_{i=1}^5 k_i\cdot{\bf g}_i\in X$$ for for $k_i\in\mathbb{Z}$. And this claim is proved. \hfill{$\Box$}

\medskip\noindent
{\bf Proof of Claim 5.5.} We firstly show that for every point {\bf p} in each $E'\in A(E)+X'$, $\kappa({\bf p})=5$ holds in $P+X'$ with $\varpi({\bf p})=2$ and $\varphi({\bf p})=3$. Without loss of generality, we can take $E'=E_1$. By the structure of the dihedral adjacent wheel at $E_1$ as shown in FIG \ref{P10translates}, we have
\begin{equation*}
  \begin{cases}
    {\bf x}_1 = {\bf o},\\
    {\bf x}_2 = {\bf g}_5,\\
    {\bf x}_3 = {\bf g}_5-{\bf g}_4,\\
    {\bf x}_4 = {\bf g}_5-{\bf g}_4+{\bf g}_3,\\
    {\bf x}_5 = {\bf g}_5-{\bf g}_4+{\bf g}_3-{\bf g}_2.
  \end{cases}
\end{equation*}
and $${\bf g}_5-{\bf g}_4+{\bf g}_3-{\bf g}_2+{\bf g}_1={\bf o}.$$ Clearly, all ${\bf x}_i\in X'$. Notice that $\varpi({\bf p}_i^*)=2$ and $\varphi({\bf p}_i^*)=3$, and it follows by the structure of dihedral adjacent wheel at each ${\bf p}_i^*$ for $i=1,2,3,4,5$ that
\begin{equation*}
  \begin{cases}
    {\bf y}_1 = {\bf g}_4-{\bf g}_3,\\
    {\bf y}_2 = {\bf g}_5-{\bf g}_3+{\bf g}_2,\\
    {\bf y}_3 = {\bf g}_5-{\bf g}_4+{\bf g}_2-{\bf g}_1,\\
    {\bf y}_4 = -{\bf g}_1+{\bf g}_3-{\bf g}_4,\\
    {\bf y}_5 = -{\bf g}_1-{\bf g}_5+{\bf g}_4.
  \end{cases}
\end{equation*}
Clearly, each ${\bf y}_i\in X'$ for $i=1,2,...,5$. By $\varphi({\bf p})=3$, we also have five subcases (1)-(5) shown in Subcase 2.2, then for each point {\bf p} in $E_1$, $\kappa({\bf p})=5$ holds in $P+X'$.

Now let $P^*$ be the projection of $P$ on the plane {\it x}{\bf o}{\it y}, and let ${\bf u}_i$ and ${\bf g}_i^*$ be the projection of $E_i$ and ${\bf g}_i$ on the plane {\it x}{\bf o}{\it y}, respectively, for $i=1,2,...,5$. Apparently, the projection of the dihedral adjacent wheel at $E_i$ is an adjacent wheel at ${\bf u}_i$ and the projection of each I-type translate at $E_i$ also contains ${\bf u}_i$ in its relative interior. In other words, the projections of the translates at $E'$ generate a fivefold tiling at the suitable neighbor of ${\bf u}_i$. Let $X^*=\{{\bf x}^*: {\bf x}^*=\sum_{i=1}^5 k_i\cdot{\bf g}_i^*\}$. For convenience, we assume that ${\bf c}_i$ is the midpoint between ${\bf u}_i$ and ${\bf u}_{i+1}$. By the argument above, it is easy to see that ${\bf c}_i\in\frac{1}{2}X^*$. According to Lemma 5.1 above and Subcase 2.2 of Lemma 3.5 of Yang and Zong \cite{yz2}, we can deduce that $X^*$ is a two-dimensional lattice and $P^*+X^*$ is a fivefold tiling in $\mathbb{E}^2$.

Thus, this claim is proved. \hfill{$\Box$}

\medskip
Since both a parallelogram and a centrally symmetric hexagon are also fivefold translative tiles, it follows by Lemmas 5.2-5.4 that we have the following Theorem 5.1.

\medskip\noindent
{\bf Theorem 5.1.} {\it Let $P$ be a fivefold translative tile in $\mathbb{E}^3$ and let $E$ be an edge of $P$, then the projection of $P$ along $E$ must be a fivefold translative tile in $\mathbb{E}^2$.}

\medskip\noindent
{\bf Remark 5.1.} In Subcase 4.1 for $\#B(E)=8$ and ${\bf y}_2={\bf y}_3$ of Lemma 5.3, according to the analysis, we can conclude that for any proper point ${\bf p}'\in E'\setminus\{C(E)+X\}$, $\varpi({\bf p}')=2$ and the structure of the dihedral adjacent wheel at ${\bf p}'$ is determined under the assumption that $F_1+{\bf x}_2$ has {\bf p} as its relative interior point at the beginning of Case 4.

\medskip\noindent
{\bf Remark 5.2.} In fact, Claim 5.3 are not only suitable for $\#B(E)=10$, but also available for the proper fivefold translative tiling $P+X$ relate to $E$ with $\#B(E)=8$.

\vspace{1.0cm}\noindent
{\LARGE\bf 6. Non-proper fivefold translative tilings }

\bigskip\noindent
In this section, we assume that $P+X$ is a non-proper fivefold translative tiling relate to $E$, see Definition 3.4. Without loss of generality, we also assume that $E$ is parallel to {\it z}-axis and define the positive direction of {\it z}-axis upward. Then we show that $P$ is a prism for $\#B(E)=10$ and $\#B(E)=8$, respectively.

\bigskip\noindent
{\Large\bf 6.1. $P+X$ for $\#B(E)=10$}

\bigskip\noindent
By Subcase 2.2.1 in Lemma 5.4 for $P$ with $\#B(E)=10$, we can directly obtain the following result.

\medskip\noindent
{\bf Lemma 6.1.} {\it If $P+X$ is a non-proper fivefold translative tiling relate to $E$ for $\#B(E)=10$, then $P$ must be a cylinder over a two-dimensional fivefold decagonal tile.}

\bigskip\noindent
{\Large\bf 6.2. $P+X$ for $\#B(E)=8$}

\bigskip\noindent
See Lemma 2.5, and we have known that there are two kinds of centrally symmetric octagonal translative tiles in $\mathbb{E}^2$.

\medskip
In $\mathrm{E}^3$, use the notations above. By Theorem 5.1, we have two kinds of $P$ whose projection along $E$ is the fivefold translative tile in $\mathbb{E}^2$ for $\#B(E)=8$. Then denote these two kinds of $P$ by $P_1$ and $P_2$, respectively. We assume that the projections of $P_1$ and $P_2$ are shown in FIG \ref{P8_E1} and FIG \ref{P8_E2}, respectively. According to the analysis in Lemma 5.3, we have known that there must be an F-type translate at some proper point {\bf p} in one $E'\in A(E)+X$ and $\varpi({\bf p})=2$. Without loss of generality, let $E'=E_8$. Then we get a dihedral adjacent wheel at {\bf p} consisting of $P+{\bf x}_1$, $P+{\bf x}_2$, $P+{\bf x}_3$, $P+{\bf x}_4$, $P+{\bf x}_5$ in sequence with ${\bf x}_1=${\bf o}, where $P+{\bf x}_4$ be the F-type translate at {\bf p} and ${\bf p}\in relint(F_1)+{\bf x}_4$. Besides, let $P+{\bf y}_1$, $P+{\bf y}_2$, $P+{\bf y}_3$ be the I-type translates at {\bf p}. For convenience, we still use the projections of translates at {\bf p} onto the plane {\it x}{\bf o}{\it y} to discuss the following cases and keep their respective notations unchanged.

\begin{figure}[h]
  \centering
  \includegraphics[scale=0.65]{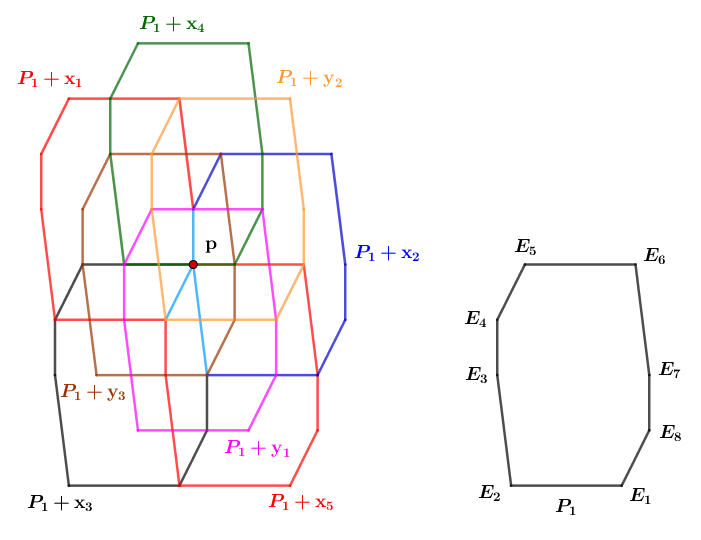}\\
  \caption{The translates at {\bf p} in $P_1+X$ for $\#B(E)=8$}\label{P8_E1}
\end{figure}

\smallskip
Recall the definition of the non-proper fivefold translative tiling relate to $E$ in Section 3. Now let us consider $P=P_1$ and $P=P_2$, respectively.

\bigskip\noindent
{\Large\bf 6.2.1. $P=P_1$}

\bigskip\noindent
According to Remark 5.1, we have that $\varpi({\bf p}')=2$ holds for all proper points ${\bf p}'$ in any $E'\in A(E)+X$ and the structure of the dihedral adjacent wheel at ${\bf p}'$ is determined. Since $P_1+X$ is a non-proper translative tiling relate to $E$, without loss of generality, by Remark 4.1, we can assume that $E_6+{\bf x}_3\not\subset\partial(P_1)+{\bf x}_j$ but $[relint(E_6)+{\bf x}_3]\cap[\partial(P_1)+{\bf x}_j]\not=\emptyset$ for some $j=1,2,4,5$. Then we get the following claim.

\medskip\noindent
{\bf Claim 6.1.} If $P_1+X$ is a non-proper fivefold translative tiling relate to $E$ for $\#B(E)=8$, then we have $${\bf x}_i+k\vec{E}\in X$$ and $${\bf y}_j+k\vec{E}\in X$$ for $i=1,2,3,4,5, y=1,2,3$ and $k\in\mathbb{Z}$.

\medskip\noindent
{\bf Proof.} By the introduction above, without loss of generality, we can consider the following two cases.

\smallskip\noindent
{\bf Case 1.} $E_6+{\bf x}_3\not=E_5+{\bf x}_5$. Since the structure of the dihedral adjacent wheel at {\bf p} determined by the E-type translate $P_1+{\bf x}_1$ and the F-type translate $P_1+{\bf x}_4$, similar to the analysis in Claim 5.1, we directly have $${\bf x}_i\pm\vec{E}\in X$$ for $i=1,2,3,5$. Since ${\bf x}_3\pm\vec{E}\in X$, we have $${\bf y}_j\pm\vec{E}\in X$$ and $${\bf x}_4\pm\vec{E}\in X.$$ Since $E_6+{\bf x}_3+k\vec{E}\not=E_5+{\bf x}_5+k\vec{E}$, we have $${\bf x}_i+k\vec{E}\in X$$ and $${\bf y}_j+k\vec{E}\in X$$ for $i=1,2,...,5$, $j=1,2,3$ and $k\in\mathbb{Z}$.

\smallskip\noindent
{\bf Case 2.} $E_6+{\bf x}_3=E_5+{\bf x}_5$. Based on (1) above, we can only discuss the case in which $E_8+{\bf x}_1=E_3+{\bf x}_2=E_6+{\bf x}_3$ and $E_6+{\bf x}_3\not\subset F_1+{\bf x}_4$. By the intersection of $F_1+{\bf x}_4$ and $F_5+{\bf x}_5$, we have $$E_1+{\bf x}_4\not\subset F_5+{\bf x}_5.$$ Since we have known that $relint[(E_7+{\bf y}_3)\cap(E_1+{\bf x}_4)]\not=\emptyset$, if $E_7+{\bf y}_3\not=E_1+{\bf x}_4$, then we can get the result as that in Case 1; if $E_7+{\bf y}_3=E_1+{\bf x}_4$, then we have the following discussion.

\smallskip\noindent
{\bf Subcase 2.1.} $\#F_1=4$, that is $F_1$ is a parallelogram. Since the structure of the dihedral adjacent wheel at each proper point ${\bf p}'\in E'$ has been determined for each $E'\in A(E)+X$, we can easily get $${\bf x}_i+k\vec{E}\in X$$ and $${\bf y}_j+k\vec{E}\in X$$ for $i=1,2,...,5$, $j=1,2,3$ and $k\in\mathbb{Z}$.

\smallskip\noindent
{\bf Subcase 2.2.} $\#F_1\ge6$. By Lemma 4.6, we have $\#F_1\le10$. Without loss of generality, we can assume that $F_1$ is a centrally symmetric decagon, as shown in FIG \ref{6_F1=6Y}. Other cases can be obtained by the similar method.

\begin{figure}[h!]
 \subfigure[]
 {\begin{minipage}{0.45\linewidth}
 \includegraphics[scale=0.55]{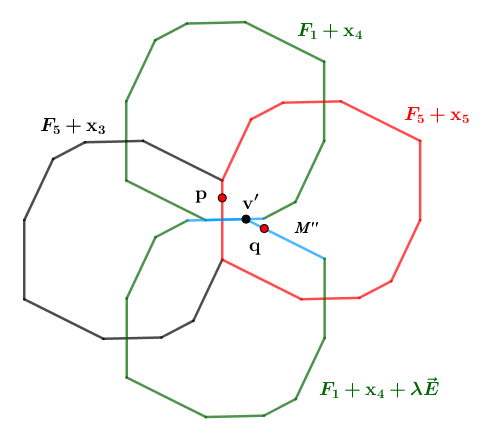}
 \end{minipage}}
 \subfigure[]
 {\begin{minipage}{0.5\linewidth}
 \includegraphics[scale=0.53]{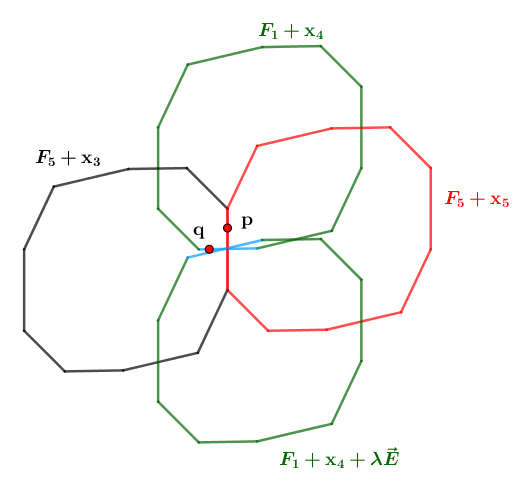}
 \end{minipage}}
 \caption{$E_6+{\bf x}_3=E_5+{\bf x}_5$ for $\#F_1=10$ in $P_1+X$}
 \label {6_F1=6Y}
\end{figure}

Since $\varphi({\bf p}')=2$ holds for all proper points ${\bf p}'$ in $E_6+{\bf x}_3$, we have ${\bf x}_4+\lambda\vec{E}\in X$ for some real number $\lambda$. For convenience, let $M$ be the edge of $F_1+{\bf x}_4$ intersecting with $E_6+{\bf x}_3$ and let $M'$ be the edge of $F_1+{\bf x}_4+\lambda\vec{E}$ intersecting with $E_6+{\bf x}_3$, then we consider the relation between $M$ and $M'$.

\smallskip\noindent
{\bf Subcase 2.2.1.} If $M\parallel M'$, then we take one edge $M''$ of $F_1+{\bf x}_4+\lambda\vec{E}$ adjacent to $M'$ but non-parallel to $E_6+{\bf x}_3$. As shown in FIG \ref{6_F1=6Y}(a), without loss of generality, let $M''$ be located between $M'$ and $E_1+{\bf x}_4+\lambda\vec{E}$ and let ${\bf v}'$ be the intersection point of $M'$ and $M''$. Take one point ${\bf q}\in M''\setminus\{C(M'')+X\}$ very close to ${\bf v}'$, and observe the translates in $P_1+X^{\bf q}$. Assume that $P_1+X_1^{\bf q}$ is one dihedral adjacent wheel containing $P_1+{\bf x}_4+\lambda\vec{E}$ and let ${\bf z}\in X_1^{\bf q}$ satisfying ${\bf q}\in F_5+{\bf z}$. Since {\bf q} is very close to ${\bf v}'$, we can assume that ${\bf v}'\in F_5+{\bf z}$.
\begin{enumerate}
  \item If ${\bf q}\in\partial(F_5)+{\bf z}$, that is {\bf q} is in one edge of $F_5+{\bf z}$ parallel to $M''$, by ${\bf v}'\in(F_1+{\bf x}_4)\cap(F_5+{\bf z})$ and Corollary 4.3, there must exist a proper point ${\bf q}'\in E_6+{\bf x}_3\setminus\{C(E)+X\}$ satisfying $${\bf q}'\in relint[(F_1+{\bf x}_4+\lambda\vec{E})\cap(F_1+{\bf z})]$$ which contradicts $\varpi({\bf q}')=2$.
  \item If ${\bf q}\in relint(F_5)+{\bf z}$, for convenience, let ${\bf z}={\bf x}_5$, then by the structure of $P_1+X_1^{\bf q}$, there is a point ${\bf z}'\in X_1^{\bf q}$ such that $${\bf q}\in\partial(F_1)+{\bf z}'.$$ Since {\bf q} is very close to ${\bf v}'$, we can assume that ${\bf v}'\in\partial(F_1)+{\bf z}'$. Similarly, by ${\bf v}'\in M''$ and Corollary 4.3, there must exist a proper point ${\bf q}'\in E_1+{\bf x}_4\setminus\{C(E)+X\}$ satisfying $${\bf q}'\in relint[(F_5+{\bf x}_5)\cap(F_1+{\bf z}')]$$ contradicting $\varpi({\bf q}')=2$.
\end{enumerate}

\smallskip\noindent
{\bf Subcase 2.2.2.} If $M\nparallel M'$, as shown in FIG \ref{6_F1=6Y}(b), we take one point ${\bf q}\in M\setminus(F_1+{\bf x}_4+\lambda\vec{E})$ and ${\bf q}\notin\{C(M)+X\}$ and consider the translates in $P+X^{\bf q}$. Similarly above, there exists a point ${\bf q}'\in(E_2+{\bf x}_4)\cap(F_5+{\bf x}_3)$ or ${\bf q}'\in(E_2+{\bf x}_4+\lambda\vec{E})\cap(F_5+{\bf x}_3)$ satisfying $\varpi({\bf q}')\ge3$ which contradicts $\varpi({\bf q}')=2$.

Therefore, Subcase 2.2 cannot exist, and we obtain the result of Claim 6.1.   \hfill{$\Box$}

\medskip
By Claim 6.1, we have the following result.

\medskip\noindent
{\bf Lemma 6.2.} {\it If $P_1+X$ is a non-proper fivefold translative tiling in $\mathbb{E}^3$ relate to $E$ for $\#B(E)=8$, then $P_1$ must be a cylinder over a two-dimensional fivefold octagonal tile.}

\medskip\noindent
{\bf Proof.} For the contrary, we suppose that $P_1$ is not a prism. Since $P_1+X$ is a non-proper fivefold translative tiling in $\mathbb{E}^3$, by the discussion in Claim 6.1, for convenience, we can assume that the upper vertex {\bf v} of $E_5+{\bf x}_5$ is in $E_6+{\bf x}_3$ and $N({\bf v},\epsilon)$ is a suitable neighborhood of {\bf v} for some small positive number $\epsilon$. Take one proper point ${\bf q}\in(E_6+{\bf x}_3)\setminus(E_5+{\bf x}_5)$ such that ${\bf q}\in N({\bf v},\epsilon)\cap[relint(F_1)+{\bf x}_4]$ and consider the translates at {\bf q}. By Claim 6.1 and $\varpi({\bf p})=\varpi({\bf q})=2$, we have ${\bf x}_5+\vec{E}\in X^{\bf q}$ and it is easy to see that ${\bf p}$ and {\bf q} have the same structure of the dihedral adjacent wheel. Let $P_1+{\bf x}_1'$,...,$P_1+{\bf x}_5'$ in sequence generate the wheel at {\bf q} where ${\bf x}_3'={\bf x}_3$, ${\bf x}_4'={\bf x}_4$ and ${\bf x}_5'={\bf x}_5+\vec{E}$, and let $P_1+{\bf y}_j'$ be the I-type translate at {\bf q} for $j=1,2,3$. It follows by Claim 6.1 that $${\bf x}_i'\in\{{\bf x}_i, {\bf x}_i+\vec{E}\}$$ and $${\bf y}_j'\in\{{\bf y}_j,{\bf y}_j+\vec{E}\}$$ for $i=1,2$ and $j=1,2,3$. In other words, if ${\bf x}_i'\not={\bf x}_i$, then ${\bf x}_i'={\bf x}_i+\vec{E}$, likewise ${\bf y}_j'$. Then we have followings steps similar to those in Subcase 2.2.1 of Lemma 5.4.

\smallskip\noindent
{\bf Step 1.} $int[(P_1+{\bf x}_5)\cap(P_1+{\bf x}_5+\vec{E})\cap N({\bf v},\epsilon)]\not=\emptyset$ by $int[(P_1+{\bf x}_5)\cap(P_1+{\bf x}_5+\vec{E})]\not=\emptyset$. Let $$C({\bf v},\epsilon)=int[(P_1+{\bf x}_5)\cap(P_1+{\bf x}_5+\vec{E})\cap N({\bf v}, \epsilon)]$$ and proceed to the next step.

\smallskip\noindent
{\bf Step 2.} $int[(P_1+{\bf x}_3)\cap C({\bf v}, \epsilon)]\not=\emptyset$ by ${\bf v}\in relint(E_6)+{\bf x}_3$. Replace $$int[(P_1+{\bf x}_6)\cap C({\bf v}, \epsilon)]$$ with $C({\bf v},\epsilon)$ and proceed to the next step.

\smallskip\noindent
{\bf Step 3.} $C({\bf v},\epsilon)\subset int[(P_1+{\bf y}_j)\cup(P_1+{\bf y}_j+\vec{E})]$ by ${\bf v}\in int[(P_1+{\bf y}_j)\cup(P_1+{\bf y}_j+\vec{E})]$ for small enough positive number $\epsilon$ and $j=1,2,3$.

According to the Steps 1-3 above, we can conclude that there exists a point ${\bf p}'\in C({\bf v},\epsilon)$ satisfying $$\varphi({\bf p}')\ge6,$$ which contradicts the assumption $P_1+X$ is a fivefold translative tiling in $\mathbb{E}^3$, thus $P_1$ must be a prism. And the proof of Lemma 6.2 is finished. \hfill{$\Box$}

\bigskip\noindent
{\Large\bf 6.2.2. $P=P_2$}

\bigskip\noindent
We now consider the non-proper fivefold translative tiling $P_2+X$. Notice that $\varpi({\bf p}+\widetilde{\bf g}_3)=2$ or $3$ by Subcase 4.2 in Lemma 5.3. Without loss of generality, correspondingly, let $P_2+{\bf y}_1$ be the F-type translate at ${\bf p}+\widetilde{\bf g}_3$ or let $P_2+{\bf y}_1$ be the E-type translate at ${\bf p}+\widetilde{\bf g}_3$, as shown in FIG \ref{P8_E2}(a) and (b), respectively.

\begin{figure}[h!]
\hspace{-3.5cm}
\subfigure[]
{\begin{minipage}{0.65\linewidth}
  \includegraphics[scale=0.65]{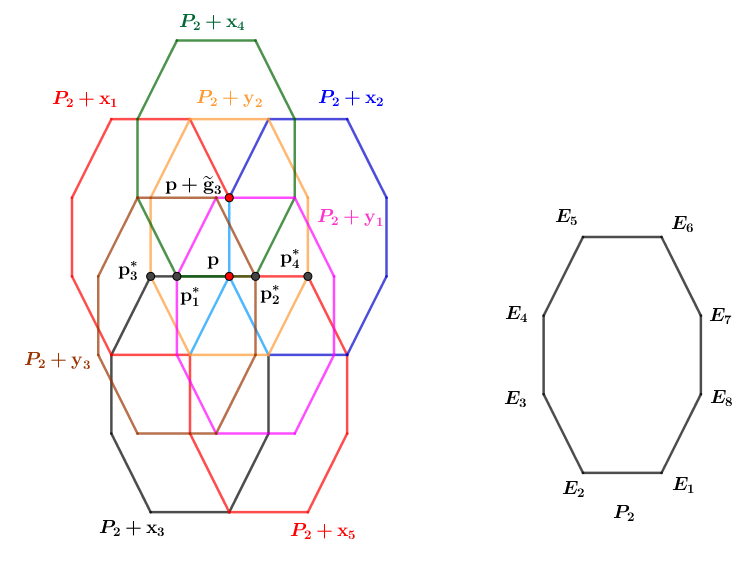}
  \end{minipage}}
\subfigure[]
{\begin{minipage}{0.8\linewidth}
  \includegraphics[scale=0.65]{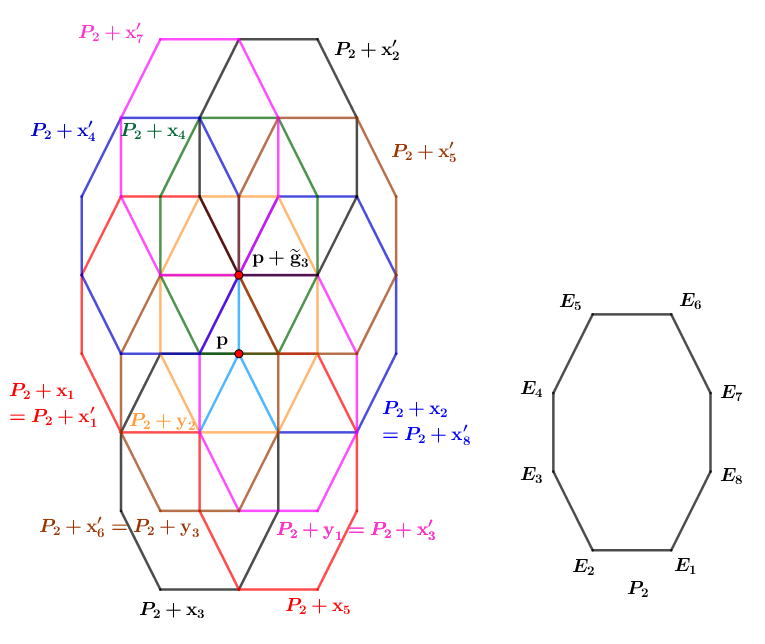}
  \end{minipage}}
\caption{The translates at {\bf p} in $P_2+X$ for $\#B(E)=8$}\label{P8_E2}
\end{figure}

\medskip\noindent
{\bf Lemma 6.3.} {\it If $P_2+X$ is a non-proper fivefold translative tiling relate to $E$ for $\#B(E)=8$, then $P_2$ must be a cylinder over a two-dimensional fivefold octagonal tile.}

\medskip\noindent
{\bf Proof.} Since $P_2+X$ is a non-proper translative tiling relate to $E$, there exists one edge $E'\in A(E)+X$ and one translate $P'\in P_2+X$ satisfying $E'\not\subset P'$ but $relint{E'}\cap \partial(P')\not=\emptyset$ by Remark 4.1. Notice that for each proper point ${\bf p}'\in E'$, we also have $\varphi({\bf p}')=3$ or $2$ by Subcase 4.2 in Lemma 5.3. Let ${\bf p}_1^*$, ${\bf p}_2^*$, ${\bf p}_3^*$ and ${\bf p}_4^*$ be the corresponding points of {\bf p} on $(F_1+{\bf x}_4)\cap(F_5+{\bf x}_3)$, $(F_1+{\bf x}_4)\cap(F_5+{\bf x}_5)$, $F_5+{\bf x}_3$ and $F_5+{\bf x}_5$, respectively.

\smallskip\noindent
{\bf Case 1.} $\varpi({\bf p}')=3$ and there is no F-type translate at ${\bf p}'$. Clearly, there are eight E-type translates $P_2+{\bf x}_1'$, $P_2+{\bf x}_2'$,..., $P_2+{\bf x}_8'$ in sequence generating the dihedral adjacent wheel and two I-type translates $P_2+{\bf y}_1'$ and $P+{\bf y}_2'$ at ${\bf p}'$. Without loss of generality, let ${\bf p}'={\bf p}+\widetilde{\bf g}_3$, ${\bf x}_1'={\bf x}_1$, ${\bf x}_8'={\bf x}_2$ and ${\bf y}_1'={\bf x}_4$, then ${\bf x}_3'={\bf y}_1$, ${\bf x}_6'={\bf y}_3$ and ${\bf y}_2'={\bf y}_2$, as shown in FIG \ref{P8_E2}(b).

\smallskip\noindent
{\bf Subcase 1.1.} $E_7+{\bf x}_1'\not=E_4+{\bf x}_8'$, in other words, if $E'=E_7+{\bf x}_1'$ and $P'=P_2+{\bf x}_8'$, then it follows that $E_8+{\bf x}_1\not=E_3+{\bf x}_2$ by ${\bf x}_1'={\bf x}_1$ and ${\bf x}_8'={\bf x}_2$. Without loss of generality, we can assume that the lower vertex of $E_8+{\bf x}_1$ is in $E_3+{\bf x}_2$ and the intersection of $(E_8+{\bf x}_1)\setminus(E_3+{\bf x}_2))$ and $F_1+{\bf x}_4$ has non-empty relative interior. Similar to Case 1 in Claim 6.1, we have $${\bf x}_2+\vec{E}\in X.$$

Let {\bf v} be the upper vertex of $E_3+{\bf x}_2$ and take a proper point ${\bf q}\in(E_8+{\bf x}_1)\setminus(E_3+{\bf x}_2))$ very close to ${\bf v}$ and having $P+{\bf x}_4$ as its F-type translate, then the structure of the dihedral adjacent wheel at {\bf q} is determined. According to the translates at {\bf p}, we correspondingly assume that $P_2+{\bf x}_1^*$,..., $P+{\bf x}_5^*$ in sequence generate the dihedral adjacent wheel at {\bf q} and $P_2+{\bf y}_j^*$ is the I-type translates at {\bf q} for $j=1,2,3$, where ${\bf x}_1^*={\bf x}_1$, ${\bf x}_2^*={\bf x}_2+\vec{E}$ and ${\bf x}_4^*={\bf x}_4$, then it follows that $${\bf x}_i^*=\{{\bf x}_i,{\bf x}_i+\vec{E}\}$$ and $${\bf y}_j^*=\{{\bf y}_j,{\bf y}_j+\vec{E}\}$$ for $i=3,5$ and $j=1,2,3$.

Let $N({\bf v},\epsilon)$ be a suitable neighborhood of {\bf v} for some small positive number $\epsilon$ and let $P^*=(P_2+{\bf x}_2)\cap(P_2+{\bf x}_2+\vec{E})$. Since $P_2$ isn't a prism, there is an open domain $C({\bf v},\epsilon)=N({\bf v},\epsilon)\cap P^*$. Similar to the steps in Subcase 2.2.1 of Lemma 5.4, we can get a new open domain $C({\bf v},\epsilon)$ satisfying that there is a point ${\bf q}'\in C({\bf v},\epsilon)$ with $\varphi({\bf q}')\ge6$ which contradicts $P_2+X$ is a fivefold translative tiling in $\mathbb{E}^3$.

\smallskip\noindent
{\bf Subcase 1.2.} $E_7+{\bf x}_1'=E_4+{\bf x}_8'$, by the discussion above and observing the translates at the points ${\bf p}'\pm\widetilde{\bf g}_2$ and ${\bf p}'\pm\widetilde{\bf g}_4$, respectively, we need $$E_1+{\bf x}_7'=E_4+{\bf x}_8'=E_7+{\bf x}_1'=E_2+{\bf x}_2'$$ and $$E_5+{\bf x}_3'=E_8+{\bf x}_4'=E_3+{\bf x}_5'=E_6+{\bf x}_6',$$ moreover, $$E_6+{\bf x}_1'=E_5+{\bf y}_2$$ and $$E_1+{\bf x}_4'=E_2+{\bf x}_4.$$

Take $E'=E_7+{\bf x}_1'$ and $P'=P_2+{\bf x}_6'$, then we need $$E_7+{\bf x}_1'\not=E_6+{\bf x}_6',$$ that is $$E_7+{\bf x}_1\not=E_6+{\bf y}_3,$$ and without loss of generality, let the upper vertex of $E_7+{\bf x}_1$ is in $E_6+{\bf y}_3$. Choose one proper point ${\bf q}\in(E_7+{\bf x}_1)\setminus(E_6+{\bf y}_3)$ very close to the lower vertex {\bf v} of $E_6+{\bf y}_3$.

\smallskip\noindent
{\bf Subcase 1.2.1.} If $\varpi({\bf q})=3$, then it follows by the structure of the dihedral adjacent wheel at ${\bf p}'$ that $${\bf x}_i'-\vec{E}\in X$$ for $i=3,4,5,6$. Let ${\bf y}_1''$ and ${\bf y}_2''$ be the two I-type translates at ${\bf p}'$. Since $E_6+{\bf x}_1'=E_5+{\bf y}_2$, we have $${\bf y}_2''={\bf y}_2'={\bf y}_2.$$ And we can also have $${\bf y}_1''\in\{{\bf y}_1', {\bf y}_1'-\vec{E}\}.$$ Notice that we have assume that $P_2$ is not a prism with $\#B(E)=8$. Take a suitable neighborhood $N({\bf v},\epsilon)$ of {\bf v} with ${\bf q}$ in its interior for some small positive number $\epsilon$, and similar to Subcase 2.2.1 in Lemma 5.4, we can conclude that there exists a point ${\bf q}'\in N({\bf v},\epsilon)$ by the translates at ${\bf p}'$ and ${\bf q}$ satisfying that $$\varphi({\bf q}')\ge6,$$ which contradicts $P_2+X$ is a fivefold translative tiling in $\mathbb{E}^3$.

\smallskip\noindent
{\bf Subcase 1.2.2.} If $\varpi({\bf q})=2$, then it follows that there is an F-type translate $P_2+{\bf y}$ at {\bf q} for ${\bf y}\in X$ and ${\bf q}\in relint(F_5)+{\bf y}$. In facet, since $P_2$ is a convex non-prism polytope, we can get that $P^*=(P_2+{\bf y})\cap[(P_2+{\bf y}_1)\cup(P_2+{\bf y}_3)]$ has nonempty interior by $(F_5+{\bf y})\cap[(F_5+{\bf y}_3)\cup(F_5+{\bf y}_1)]\not=\emptyset$ and the following cases. Since {\bf q} is very close to {\bf v}, for convenience, we can assume that ${\bf v}\in F_5+{\bf y}$.
\begin{enumerate}
 \item $\#F_1=8$ or $10$. It is easy to see that $(P_2+{\bf y})\cap[(P_2+{\bf y}_1)\cup(P_2+{\bf y}_3)]$ has nonempty interior by $relint\{(F_5+{\bf y})\cap[(F_5+{\bf y}_3)\cup(F_5+{\bf y}_1)]\}\not=\emptyset$ and the convexity of $P_2$.
 \item $\#F_1=6$. Let $M$ be one edge of $F_1$ adjacent to $E_1$ at the upper vertex of $E_1$, let $M^*$ be another edge of $F_1$ adjacent to $E_1$. If $relint\{(F_5+{\bf y})\cap[(F_5+{\bf y}_3)\cup(F_5+{\bf y}_1)]\}\not=\emptyset$, then $(P_2+{\bf y})\cap[(P_2+{\bf y}_1)\cup(P_2+{\bf y}_3)]$ has nonempty interior. Then we consider the case for $relint\{(F_5+{\bf y})\cap[(F_5+{\bf y}_3)\cup(F_5+{\bf y}_1)]\}=\emptyset$ and ${\bf q}\in relint(F_5)+{\bf y}$. Since $P_2$ is not a prism, we must have $\#B(M)\ge6$ or $\#B(M^*)\ge6$. Without loss of generality, we assume that $\#B(M^*)\ge6$. Then we have that $(P_2+{\bf y}_1)\cap(P_2+{\bf y})$ has non-empty interior by Corollary 4.1 and ${\bf y}={\bf y}_1+\vec{E}+\vec{M}$, where $\vec{M}$ is the vector form of $M$ and its starting point is the upper vertex of $E_1$.
 \item $\#F_1=4$. Let $M$ be one edge of $F_1$ adjacent to $E_1$ at the upper vertex of $E_1$. Since $P_2$ is not a prism, we must have $\#B(M)\ge6$. Then we also have that $(P_2+{\bf y})\cap[(P_2+{\bf y}_1)\cup(P_2+{\bf y}_3)]$ has nonempty interior.
\end{enumerate}

We have known that ${\bf q}$ and ${\bf p}'$ has one identical I-type translate $P_2+{\bf y}_2$ by $E_5+{\bf y}_2=E_6+{\bf x}_1$. By $\varpi({\bf q})=2$, let $P_2+{\bf y}_j''$ be the I-type translate at {\bf q} for $j=1,2,3$ and ${\bf y}_2''={\bf y}_2$. Then by the structure of the dihedral adjacent wheel at {\bf q} and the analysis above, we have ${\bf v}\in int[(P_2+{\bf y}_1'')\cup(P_2+{\bf y}_3'')\cup(P_2+{\bf x}_4)]$ has non-empty interior.

Let $N({\bf v},\epsilon)$ be a suitable neighborhood of {\bf v} for some small positive number $\epsilon$ and let $P^*=(P_2+{\bf y})\cap[(P_2+{\bf y}_1)\cup(P_2+{\bf y}_3)]$. Since $P_2$ isn't a prism, there is an open domain $C({\bf v},\epsilon)=N({\bf v},\epsilon)\cap P^*$. Similar to the steps in Subcase 2.2.1 of Lemma 5.4, we can get a new open domain $C({\bf v},\epsilon)$ satisfying that there is a point ${\bf q}'\in C({\bf v},\epsilon)$ with $\varphi({\bf q}')\ge6$ which contradicts $P_2+X$ is a fivefold translative tiling in $\mathbb{E}^3$.

Then by the translates at {\bf q} and Lemma 4.3, we obtain that there is a point ${\bf q}'\in N({\bf v},\epsilon)$ satisfying that $$\varphi({\bf q}')\ge6.$$

Therefore, Case 1 cannot exist.

\smallskip\noindent
{\bf Case 2.}  $\varpi({\bf p}')=2$ and there is an F-type translate at ${\bf p}'$. Without loss of generality, let ${\bf p}'={\bf p}$, $E'=E_6+{\bf x}_3$, and $P'=P_2+{\bf x}$ for some ${\bf x}\in\{{\bf x}_2, {\bf x}_3, {\bf x}_4, {\bf x}_5\}$.

If ${\bf x}\not={\bf x}_4$, without loss of generality, let ${\bf x}={\bf x}_5$, that is $$E_6+{\bf x}_3\not=E_5+{\bf x}_5.$$ Naturally, we have ${\bf p}\in relint[(E_6+{\bf x}_3)\cap(E_5+{\bf x}_5)]$. For convenience, let the upper vertex {\bf v} of $E_5+{\bf x}_5$ is in $E_6+{\bf x}_3$ and $relint\{[(E_6+{\bf x}_3)\setminus(E_5+{\bf x}_5)]\cap[relint(F_1)+{\bf x}_4]\}\not=\emptyset$ and let $N({\bf v},\epsilon)$ be a suitable neighborhood of {\bf v} for some small positive number $\epsilon$. Take one proper point ${\bf q}\in(E_6+{\bf x}_3)\setminus(E_5+{\bf x}_5)$ such that ${\bf q}\in N({\bf v},\epsilon)\cap[relint(F_1)+{\bf x}_4]$.
By Claim 6.1 and $\varpi({\bf p})=\varpi({\bf q})=2$, we have ${\bf x}_5+\vec{E}\in X^{\bf q}$ and it is easy to see that ${\bf p}$ and {\bf q} have the same structure of the dihedral adjacent wheel. Let $P_2+{\bf x}_1'$,...,$P_2+{\bf x}_5'$ in sequence generate the wheel at {\bf q} where ${\bf x}_2'={\bf x}_2+\vec{E}$, ${\bf x}_3'={\bf x}_3$ and ${\bf x}_4'={\bf x}_4$, and let $P_2+{\bf y}_j'$ be the I-type translate at {\bf q} for $j=1,2,3$. It follows by Claim 6.1 that $${\bf x}_i'\in\{{\bf x}_i, {\bf x}_i+\vec{E}\}$$ and $${\bf y}_j'\in\{{\bf y}_j,{\bf y}_j+\vec{E}\}$$ for $i=1,5$ and $j=1,2,3$. Then we have the steps similar to those in Lemma 6.2, and finally we can also get a new open domain $C({\bf v},\epsilon)$ of {\bf v} satisfying that there is a point ${\bf p}'\in C({\bf v},\epsilon)$ satisfying $\varphi({\bf p}')\ge6$, which contradicts the assumption $P_2+X$ is a fivefold translative tiling in $\mathbb{E}^3$. For ${\bf x}={\bf x}_2$ or ${\bf x}={\bf x}_3$, similarly, we can get the same contradiction above.

Then we can only consider the case ${\bf x}={\bf x}_4$, that is $E_6+{\bf x}_3\not\subset F_1+{\bf x}_4$. Moreover, according to the discussion above, we can assume that $$E_6+{\bf x}_3=E_5+{\bf x}_5=E_8+{\bf x}_1=E_3+{\bf x}_2.$$ Since {\bf p} is a proper point, ${\bf p}_i^*\notin C(E)+X$ for $i=1,2,3,4$. By FIG \ref{P8_E2}, we have that $relint\left[(E_7+{\bf y}_3)\cap(E_1+{\bf x}_4)\right]\not=\emptyset$. By the analysis for ${\bf x}\not={\bf x}_4$ above, we can conclude that $$E_7+{\bf y}_3=E_1+{\bf x}_4,$$ since the F-type translate $P_2+{\bf x}_5$ at ${\bf p}_2^*$ has been determined. Similarly, we also have that $$E_4+{\bf y}_1=E_2+{\bf x}_4.$$ Since $4\le\#F_1\le10$ by Lemma 4.6, we have the following three cases.

\smallskip\noindent
{\bf Subcase 2.1.} $\#F_1=8$ or 10. Let $M$ be one edge of $F_1$ adjacent to $E_1$ at the upper vertex of $E_1$ and let {\bf v} be the upper vertex of $E_6+{\bf x}_3$. Then $M+{\bf g}_1+{\bf x}_3$ is the edge of $F_5+{\bf x}_3$. Since $relint[(F_5+{\bf x}_3)\cap(F_1+{\bf x}_4)]\not=\emptyset$ by ${\bf p}\in relint(F_1+{\bf x}_4)$, $(F_5+{\bf x}_3)\cap(F_1+{\bf x}_4)$ has at least four edges. We choose a point ${\bf q}\in (M+{\bf g}_1+{\bf x}_3)\setminus\{C(M)+X\}$ very close to {\bf v} such that ${\bf q}\in relint(F_1+{\bf x}_4)$. And consider one dihedral adjacent wheel $P_2+X_1^{\bf q}$ which contains $P_2+{\bf x}_3$ and let ${\bf z}\in X_1^{\bf q}$ satisfying that ${\bf q}\in F_1+{\bf z}$, as shown in FIG \ref{F1=10Y}.

\begin{figure}[h!]
\subfigure[]
{\begin{minipage}{0.5\linewidth}
  \includegraphics[scale=0.55]{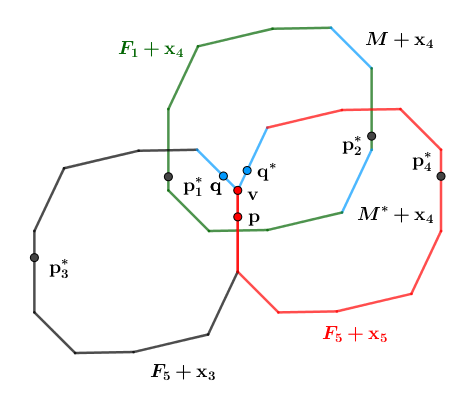}
  \end{minipage}}
\subfigure[]
{\begin{minipage}{0.45\linewidth}
  \includegraphics[scale=0.55]{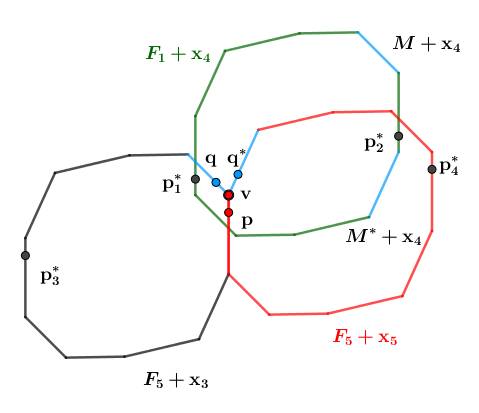}
  \end{minipage}}
 \caption{$E_6+{\bf x}_3=E_5+{\bf x}_5$ for $\#F_1=10$ in $P_2+X$}
 \label {F1=10Y}
\end{figure}

\smallskip\noindent
{\bf Subcase 2.1.1.} If ${\bf q}\in\partial(F_1)+{\bf z}$, since ${\bf q}$ is very close to {\bf v}, without loss of generality, we can assume that ${\bf v}\in\partial(F_1)+{\bf z}$, then there is always a point ${\bf p}'\in relint[(F_1+{\bf z})\cap(E_2+{\bf x}_4)]$ or ${\bf p}'\in relint[(F_1+{\bf z})\cap(E_5+{\bf x}_3)]$ by Corollary 4.3 and there would be two F-type translates at ${\bf p}'$ contradicting $\varpi({\bf p}')=2$.

\smallskip\noindent
{\bf Subcase 2.1.2.} If ${\bf q}\in relint(F_1)+{\bf z}$, without loss of generality, we can assume that ${\bf z}={\bf x}_4$, then there is ${\bf z}'\in X_1^{\bf q}$ satisfying that ${\bf q}\in\partial(F_5)+{\bf z}'$. Since ${\bf q}$ is very close to {\bf v}, without loss of generality, we can assume that ${\bf v}\in\partial(F_5)+{\bf z}'$. By observation, ${\bf z}'={\bf x}_5+(\vec{E}+\lambda\vec{M})$ for $0<\lambda\le1$ where $\vec{M}$ is the vector form of $M$ and its starting point is the upper vertex of $E$.
\begin{enumerate}
 \item If the edge $N+{\bf x}_4$ of $F_1+{\bf x}_4$ intersecting $E_5+{\bf x}_5$ is not a translate of $M$ as shown in FIG \ref{F1=10Y}(a), then $(F_5+{\bf x}_5)\cap(F_5+{\bf z}')$ contains $-N+{\bf x}_5$ by Corollary 4.3. Let $F'=(F_5+{\bf x}_5)\cap(F_5+{\bf z}')$. Similar to Subcase 1.1, there is a point ${\bf p}''\in relint[(E_1+{\bf x}_4)\cap F']$, which contradicts $\varpi({\bf p}'')=2$.
 \item If $N+{\bf x}_4$ is a translate of $M$ as shown in FIG \ref{F1=10Y}(b), there is no way to get such ${\bf p}''$ in (1), then we have to take into another edge $M^*$ of $F_1$ which is adjacent to $E_1$ at the lower vertex of $E_1$. It is easy to see that $-M^*+{\bf x}_5$ is associated with {\bf v}. We choose a point ${\bf q}^*\in relint(-M^*+{\bf x}_5)$ very close to {\bf v}. In the same way, we consider the dihedral adjacent $P_2+X_1^{{\bf q}^*}$ which contains $P_2+{\bf x}_5$. Based on the analysis above, we can only consider the case in which both ${\bf x}_4$ and ${\bf x}_3+\vec{E}-\lambda\vec{M^*}$ belong to $X_1^{{\bf q}^*}$ for $0<\lambda\le1$ where $\vec{M^*}$ is the vector form of $M^*$ and its starting point is the lower vertex of $E_1$. Since $M^*$ is not a translate of $M$, $(F_5+{\bf x}_3)\cap(F_5+{\bf x}_3+\vec{E}-\lambda\vec{M^*})$ contains $M+{\bf g}_1+{\bf x}_3$ as its edge. Let $F^*=(F_5+{\bf x}_3)\cap(F_5+{\bf x}_3+\vec{E}-\lambda\vec{M^*})$. Since $N+{\bf x}_4$ intersects $E_6+{\bf x}_3$, we have that $M+{\bf g}_1+{\bf x}_3$ intersects $E_2+{\bf x}_4$ and $relint[(E_2+{\bf x}_4)\cap F^*]\not=\emptyset$, which is also impossible.
\end{enumerate}

\smallskip\noindent
{\bf Subcase 2.2.} $\#F_1=6$. Let $M$ be one edge of $F_1$ adjacent to $E_1$ at the upper vertex of $E_1$, let $M^*$ be another edge of $F_1$ adjacent to $E_1$, and let {\bf v} be the upper vertex of $E_6+{\bf x}_3$. Similar to the analysis in Subcase 2.1, we can get that there is a translate $P_2+{\bf x}_5+\vec{E}+\vec{M}$ for $\#F_1=6$, as shown in FIG \ref{F1=6Y}, where $\vec{M}$ is the vector form of $M$ and its starting point is the upper vertex of $E_1$.

\begin{figure}[h!]
\subfigure[]
{\begin{minipage}{0.5\linewidth}
  \includegraphics[scale=0.64]{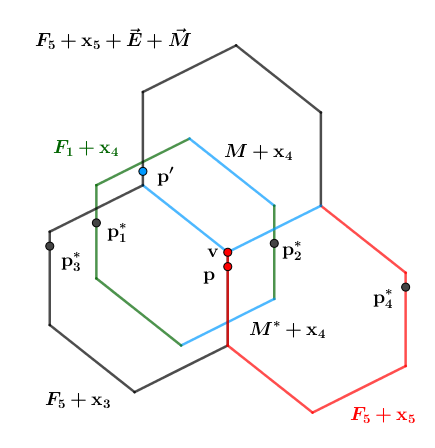}
  \end{minipage}}
\subfigure[]
{\begin{minipage}{0.45\linewidth}
  \includegraphics[scale=0.66]{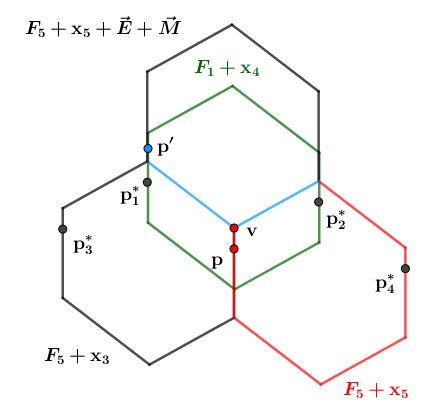}
  \end{minipage}}
 \caption{$E_6+{\bf x}_3=E_5+{\bf x}_5$ for $\#F_1=6$ in $P_2+X$}
 \label {F1=6Y}
\end{figure}

Let $N({\bf v},\epsilon)$ be a neighborhood of {\bf v} for some small positive number $\epsilon$. Notice that since $P_2$ is not a prism, we must have $\#B(M)\ge6$ or $\#B(M^*)\ge6$. Without loss of generality, we assume that $\#B(M^*)\ge6$. Then we have that $(P_2+{\bf x}_5)\cap(P_2+{\bf x}_5+\vec{E}+\vec{M})$ has non-empty interior by Corollary 4.1, written as $P^*$. Then we have the following two subcases.

\smallskip\noindent
{\bf Subcase 2.2.1.} For FIG \ref{F1=6Y}(a), $(E_6+{\bf x}_3)\cap(F_1+{\bf x}_4)$ doesn't contain the vertex of $F_1+{\bf x}_4$.

Let ${\bf p}'$ be a proper point on $(E_5+{\bf x}_5+\vec{E}+\vec{M})\cap(F_1+{\bf x}_4)$, then clearly ${\bf p}'$ and {\bf p} have the same F-type translate $P_2+{\bf x}_4$, and the dihedral adjacent wheels $P_2+X^{\bf p}$ and $P_2+X^{{\bf p}'}$ have the same structure, thus $P_2+{\bf x}_2+\vec{E}+\vec{M}$ is an E-type translate at ${\bf p}'$. And we have the following steps.

\smallskip\noindent
{\bf Step 1.} $int[N({\bf v},\epsilon)\cap P^*]\not=\emptyset$. Then let $$C({\bf v},\epsilon)=int[N({\bf v},\epsilon)\cap P^*]$$ and proceed to the next step.

\smallskip\noindent
{\bf Step 2.} $int[(P_2+{\bf x}_2+\vec{E}+\vec{M})\cap C({\bf v},\epsilon)]\not=\emptyset$ by $relint(F_5+{\bf x}_5+\vec{E}+\vec{M})\subset int(P_2+{\bf x}_2+\vec{E}+\vec{M})$ and ${\bf v}\in P_2+{\bf x}_2+\vec{E}+\vec{M}$. Then replace $$int[(P_2+{\bf x}_2+\vec{E}+\vec{M})\cap C({\bf v},\epsilon)]$$ with $C({\bf v},\epsilon)$, and proceed to the next step.

\smallskip\noindent
{\bf Step 3.} $C({\bf v},\epsilon)\subset int(P_2)+{\bf y}_j$ by ${\bf v}\in int(P_2)+{\bf y}_j$ for small enough positive number $\epsilon$ and $j=1,2,3$.

Based on the steps above, we can conclude that there is a point ${\bf q}\in C({\bf v},\epsilon)$ such that $$\varphi({\bf q})\ge6,$$ which contradicts the assumption that $P_2+X$ is a fivefold translative tiling. Thus this case is impossible.

\smallskip\noindent
{\bf Subcase 2.2.2.} For FIG \ref{F1=6Y}(b), $(E_6+{\bf x}_3)\cap(F_1+{\bf x}_4)$ contains a vertex of $F_1+{\bf x}_4$. Let ${\bf p}'$ be a proper point on $(E_5+{\bf x}_5+\vec{E}+\vec{M})\cap(F_1+{\bf x}_4)$. Since $P_2+{\bf x}_4$ and $P_2+{\bf x}_5+\vec{E}+\vec{M}$ are two E-type translates at ${\bf p}'$, we have $$\varpi({\bf p}')=3$$ and $$E_5+{\bf x}_5+\vec{E}+\vec{M}\not=E_2+{\bf x}_4.$$ It follows by Subcase 1.1 above that this case cannot exist.

\smallskip\noindent
{\bf Subcase 2.3.} $\#F_1=4$. Let ${\bf p}'$ be one proper point on $(E_6+{\bf x}_3)\setminus(F_1+{\bf x}_4)$ and let $M$ be an edge of $F_1$ adjacent to $E_1$ at the upper vertex of $E_1$, then we consider the translates at ${\bf p}'$. Since $F_1$ is a parallelogram and $P_2$ is not a prism, $\#B(M)\ge6$. Let {\bf v} be the intersection point of $E_6+{\bf x}_3$ and $M-\vec{E}$, and let $N({\bf v},\epsilon)$ be a suitable neighborhood at {\bf v} for a small positive number $\epsilon$. For convenience, we assume that ${\bf p}'$ is very close to {\bf v}.

\begin{figure}[h!]
\subfigure[]
{\begin{minipage}{0.4\linewidth}
  \includegraphics[scale=0.68]{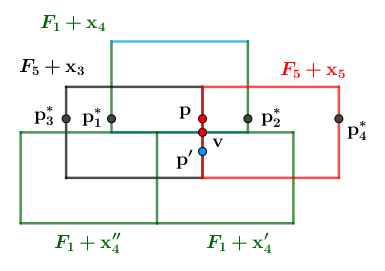}
  \end{minipage}}
 \subfigure[]
{\begin{minipage}{0.5\linewidth}
\includegraphics[scale=0.68]{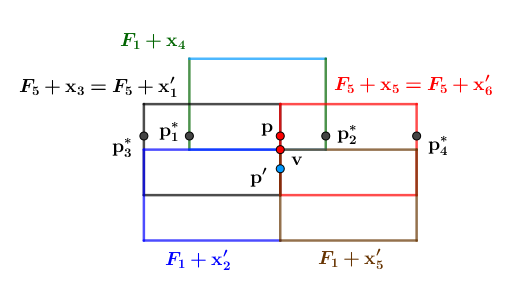}
  \end{minipage}}
 \caption{$E_6+{\bf x}_3=E_5+{\bf x}_5$ for $\#F_1=4$ in $P_2+X$}\label{F1=4Y}
\end{figure}

\smallskip\noindent
{\bf Subcase 2.3.1.} There is an F-type translate $P_2+{\bf x}_4'$ at ${\bf p}'$. Then we can assume that $P_2+{\bf x}_i'$ in sequence and $P_2+{\bf y}_j'$ are respectively the E-type translates and I-type translates at ${\bf p}'$ for $i=1,2,3,4,5$ and $j=1,2,3$. Without loss of generality, let $P+{\bf x}_4'$ is the F-type translate at ${\bf p}'$ and it is easy to see that ${\bf y}_2'={\bf y}_2$ and ${\bf x}_i'={\bf x}_i$ for all $i\not=4$, and $${\bf x}_4'={\bf x}_4-\vec{E}\pm\lambda\vec{M}$$ for $0\le\lambda<1$ where $\vec{M}$ is the vector form of $M$ and its starting point is the upper vertex of $E_1$. Without loss of generality, take ${\bf x}_4'={\bf x}_4-\vec{E}-\lambda\vec{M}$ as shown in FIG \ref{F1=4Y}(a), then we can assume that ${\bf y}_j'={\bf y}_j-\vec{E}-\lambda\vec{M}$ for $j=1,3$. By Lemma 4.3, we have that $$int[(P_2+{\bf x}_4)\cap(P_2+{\bf x}_4')]\not=\emptyset$$ and $$int[(P_2+{\bf y}_j)\cap(P_2+{\bf y}_j')]\not=\emptyset$$ for $j=1,3$. For convenience, let $P^*=(P_2+{\bf x}_4)\cap(P_2+{\bf x}_4')$ Then we have the following steps.

\smallskip\noindent
{\bf Step 1.} $int[N({\bf v},\epsilon)\cap P^*]\not=\emptyset$ by ${\bf v}\in P^*$. Then replace $int[N({\bf v},\epsilon)\cap P^*]$ with $C({\bf v},\epsilon)$ and proceed to the next step.

\smallskip\noindent
{\bf Step 2.} $int[(P_2+{\bf x}_1)\cap C({\bf v},\epsilon)]\not=\emptyset$ by ${\bf v}\in relint(E_8)+{\bf x}_1$. Replace $int[(P_2+{\bf x}_1)\cap C({\bf v},\epsilon)]$ with $C({\bf v},\epsilon)$, then we proceed to the next step.

\smallskip\noindent
{\bf Step 3.} $C({\bf v},\epsilon)\subset int(P_2)+{\bf y}_2$ by ${\bf v}\in int(P_2)+{\bf y}_2$ for small enough positive number $\epsilon$.

\smallskip\noindent
{\bf Step 4.} $C({\bf v},\epsilon)\subset int[(P_2+{\bf y}_j)\cup(P_2+{\bf y}_j')]$ by ${\bf v}\in int[(P_2+{\bf y}_j)\cup(P_2+{\bf y}_j')]$ for $j=1,3$ and small enough positive number $\epsilon$.

Based on the steps above, we can conclude that there is a point ${\bf q}\in C({\bf v},\epsilon)$ satisfying that $$\varphi({\bf q})\ge 6,$$ which contradicts that $P_2+X$ is a fivefold translative tiling.

\smallskip\noindent
{\bf Subcase 2.3.2.} There is no F-type translate at ${\bf p}'$, as in FIG \ref{F1=4Y}(b), then we have $$\varpi({\bf p}')=3$$ and $$E_1+{\bf x}_2'\not=E_6+{\bf x}_1'.$$ It follows by Subcase 1.2 above that this case cannot exist.

Therefore, we have finished the proof of Lemma 6.3. \hfill{$\Box$}

\vspace{1cm}\noindent
{\LARGE\bf 7. Proper fivefold translative tilings}

\bigskip\noindent
{We assume that $P+X$ is a proper fivefold translative tiling relate to $E$ in this section, see Definition 3.3. Then we show that $P$ is a prism for $\#B(E)=10$ and $\#B(E)=8$, respectively.

Let $F\in B(E)$ and let $M$ be another edge of $F$ which is not parallel to $E$. For clarity, let $M_1^*, M_2^*,...,M_{2m'}^*$ in circular order be the edges of $P$ parallel to $M$ where $M_{2m'+1}^*=M_1^*$ and let $B(M)=\{F_1^*, F_2^*,...,F_{2m'}^*\}$ where $M_i^*$ and $M_{i+1}^*$ are two edges of $F_i^*$ and $m'\in\{2,3,4,5\}$. Moreover, let ${\bf g}_i^*$ be the translate vector of $F_i^*$ in $P$ and let $\widetilde{\bf g}_i^*$ be the translate vector of $M_i^*$ in $F_i^*$, in other words, $F_{i+m'}^*=-F_i^*=F_i^*+{\bf g}_i^*$ and $M_{i+1}^*=M_i^*+\widetilde{\bf g}_i^*$ for $1\le i\le m'$. By Lemma 4.2, it is easy to get $\pm F\in B(E)\cap B(M)$.

\bigskip\noindent
{\Large\bf 7.1. $P+X$ for $\#B(E)=10$}

\bigskip\noindent
Recall that {\bf p} is a proper point on $E'\in A(E)+X$, $P+{\bf x}_1$, $P+{\bf x}_2$,..., $P+{\bf x}_5$ in clock order generate a dihedral adjacent wheel at {\bf p} all of which are E-type translates and $P+{\bf y}_1$, $P+{\bf y}_2$, $P+{\bf y}_3$ are the I-type translates at {\bf p}. Since $P+X$ is a proper fivefold translative tiling relate to $E$, we also call $P+{\bf x}_i$ and $P+{\bf y}_j$ {\it E-type translates} and {\it I-type translates} at $E'$, respectively, by Claim 5.3. Without loss of generality, let $E'=E_1$. By observing the dihedral adjacent wheel at $E_1$ in FIG \ref{P10_E}, we can get that $F'\in\{F_1+{\bf x}_1, F_5+{\bf x}_2, F_4+{\bf x}_2, F_3+{\bf x}_4, F_2+{\bf x}_4\}$ is a translate of $F$. Let $M^*$ be the edge of $F'$ that is a translate of $M$. Then we mainly analyze two cases following in this subsection.

\begin{figure}[h!]
\includegraphics[scale=0.6]{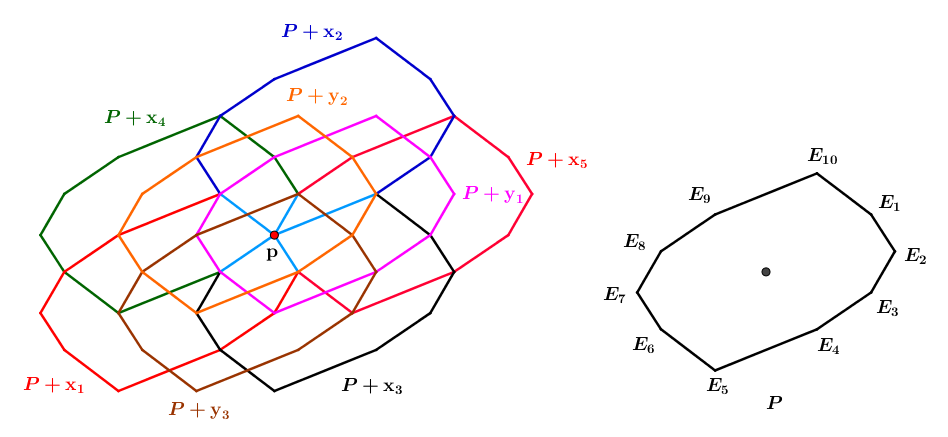}
\caption{The dihedral adjacent wheel at $E$ for $\#B(E)=10$}\label{P10_E}
\end{figure}

\smallskip\noindent
{\bf Case 1.} $F'\in\{F_1+{\bf x}_1, F_3+{\bf x}_4, F_2+{\bf x}_4\}$. By Lemma 3.1 and Corollary 4.1, we have that there are four translates of $P$ containing $relint(F')$ in their interior in this case. And these four translates can be divided into two groups by symmetry about the center of $F'$. Without loss of generality, we only analyze the case for $F'=F_1+{\bf x}_1$ and another two cases can be obtained by the similar method.

If $F'=F_1+{\bf x}_1$, we have $F'\subset P+{\bf x}$ for ${\bf x}\in\{{\bf x}_1, {\bf x}_5, {\bf x}_3, {\bf y}_2,{\bf y}_3, {\bf y}_1\}$ and $relint(F')\subset int(P)+{\bf x}$ for ${\bf x}\in\{{\bf x}_3, {\bf y}_2,{\bf y}_3, {\bf y}_1\}$. Then by simple observation and analysis, we get that $(P+{\bf x}_3)\cap(P+{\bf y}_2)$ and $(P+{\bf y}_3)\cap(P+{\bf y}_1)$ are both centrally symmetric about the center of $F'$. Since $P+X$ is a proper fivefold translative tiling relate to $E$, by Corollary 4.1, we can conclude that $(P+{\bf x}_3)\cap(P+{\bf y}_2)$ has a belt $B(E,{\bf x}_3,{\bf y}_2)$ and $(P+{\bf y}_3)\cap(P+{\bf y}_1)$ has a belt $B(E,{\bf y}_3,{\bf y}_1)$, seen in Remark 4.2.

\smallskip\noindent
{\bf Case 2.} $F'\in\{F_5+{\bf x}_2, F_4+{\bf x}_2\}$. By Lemma 3.1 and Corollary 4.1, we have that there are three translates of $P$ containing $relint(F')$ in their interior. And these three translates can be divided into two groups by symmetry about the center of $F'$. Without loss of generality, we only analyze the case for $F'=F_4+{\bf x}_2$ and the other case can be obtained by the similar method.

If $F'=F_4+{\bf x}_2$, we have $F'\subset P+{\bf x}$ for ${\bf x}\in\{{\bf x}_2,{\bf x}_3, {\bf x}_5,{\bf y}_2, {\bf y}_1\}$ and $relint(F')\subset int(P)+{\bf x}$ for ${\bf x}\in\{{\bf x}_5,{\bf y}_2, {\bf y}_1\}$. By simple observation and analysis, we can get that $(P+{\bf x}_5)\cap(P+{\bf y}_2)$ and $P+{\bf y}_1$ are both centrally symmetric about the center of $F'$. Similar to Case 1, we can conclude that $(P+{\bf x}_5)\cap(P+{\bf y}_2)$ has a belt $B(E,{\bf x}_5,{\bf y}_2)$.

Then we have the following lemmas.


\medskip\noindent
{\bf Lemma 7.1.} {\it If $P$ is a fivefold translative tile in $\mathbb{E}^3$ with $\#B(E)=10$, then $\#B(M)\not=10$.}

\medskip\noindent
{\bf Proof.} For the contrary, suppose that $P$ is a fivefold translative tile with $\#B(E)=\#B(M)=10$, then $P+X$ is both a proper fivefold translative tiling relate to $E$ and a proper fivefold translative tiling relate to $M$ based on the results above. Let $P+{\bf x}_1^*$,..., $P+{\bf x}_5^*$ in sequence generate the dihedral adjacent wheel at $M^*$ and let $P+{\bf y}_j^*$ be the I-type translates at $M^*$ for $j=1,2,3$. And we can also have the projection of the dihedral adjacent wheel at $M^*$ along $M$ shown in FIG \ref{P10_M}, where ${\bf q}$ is a relative interior point on $M^*$.

\begin{figure}[h!]
\includegraphics[scale=0.75]{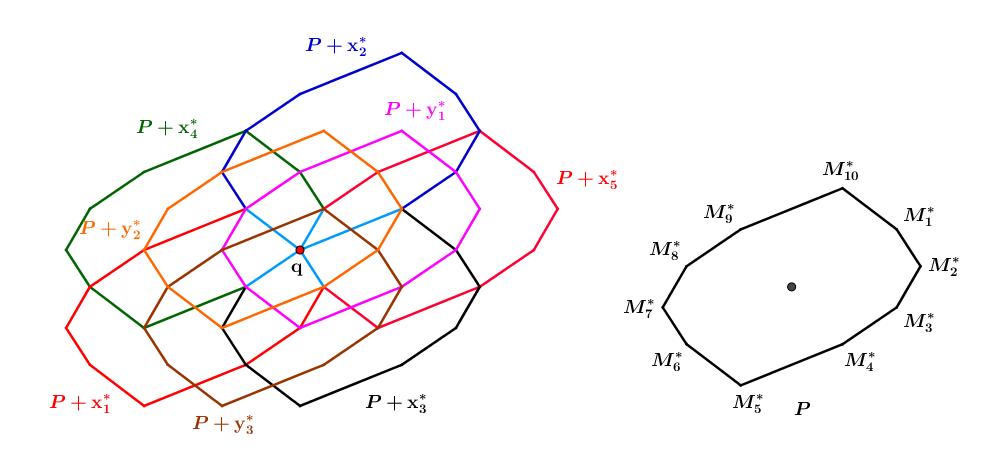}
\caption{The dihedral adjacent wheel at $M^*$ for $\#B(M)=10$}\label{P10_M}
\end{figure}

\smallskip\noindent
{\bf Case 1.} $F'\in\{F_1+{\bf x}_1, F_3+{\bf x}_4, F_2+{\bf x}_4\}$ and $F'=F_1+{\bf x}_1$. Let us observe the dihedral adjacent wheel at $M^*$, and we have $F'\in\{F_1^*+{\bf x}_1^*, F_3^*+{\bf x}_4^*, F_2^*+{\bf x}_4^*\}$, since there are four translates of $int(P)$ containing $relint(F')$.

\smallskip\noindent
{\bf Subcase 1.1.} If $F'=F_1^*+{\bf x}_1^*$, we have $relint(F')\subset int(P)+{\bf x}$ for ${\bf x}\in\{{\bf x}_3^*, {\bf y}_2^*,{\bf y}_3^*, {\bf y}_1^*\}$, and both $(P+{\bf x}_3^*)\cap(P+{\bf y}_2^*)$ and $(P+{\bf y}_3^*)\cap(P+{\bf y}_1^*)$ are centrally symmetric about the center of $F'$. Similarly above, $(P+{\bf x}_3^*)\cap(P+{\bf y}_2^*)$ has a belt $B(M,{\bf x}_3^*,{\bf y}_2^*)$ and $(P+{\bf y}_3^*)\cap(P+{\bf y}_1^*)$ has a belt $B(M,{\bf x}_3^*,{\bf y}_2^*)$. And we have $$\{{\bf x}_3, {\bf y}_2\}=\{{\bf x}_3^*, {\bf y}_2^*\}$$ or $$\{{\bf x}_3, {\bf y}_2\}=\{{\bf y}_3^*, {\bf y}_1^*\},$$ that is $(P+{\bf x}_3)\cap(P+{\bf y}_2)$ also has a belt $B(M,{\bf x}_3^*,{\bf y}_2^*)$ or $B(M,{\bf y}_3^*,{\bf y}_1^*)$. Then $(P+{\bf x}_3)\cap(P+{\bf y}_2)$ should have two facets parallel to $F'$, which is impossible since $B(E,{\bf x}_3,{\bf y}_2)$ doesn't have any facet parallel to $F'$ by observation and Corollary 4.1.

\smallskip\noindent
{\bf Subcase 1.2.} If $F'=F_3^*+{\bf x}_4^*$, by similar analysis to Subcase 1.1 above, this case is also impossible.

\smallskip\noindent
{\bf Subcase 1.3.} If $F'=F_2^*+{\bf x}_4^*$, by similar analysis to Subcase 1.1 above and Lemma 4.3, this case is also impossible.

\smallskip\noindent
{\bf Case 2.} $F'\in\{F_5+{\bf x}_2, F_4+{\bf x}_2\}$ and $F'=F_4+{\bf x}_2$. Let us observe the dihedral adjacent wheel at $M^*$, and we have $F'\in\{F_4^*+{\bf x}_2^*, F_5^*+{\bf x}_2^*\}$.

\smallskip\noindent
{\bf Subcase 2.1.} If $F'=F_4^*+{\bf x}_2^*$, we have $relint(F')\subset int(P)+{\bf x}$ for ${\bf x}\in\{{\bf x}_5^*,{\bf y}_2^*,{\bf y}_1^*\}$, and $(P+{\bf x}_5^*)\cap(P+{\bf y}_2^*)$ is centrally symmetric about the center of $F'$. Similarly, we can conclude that $(P+{\bf x}_5^*)\cap(P+{\bf y}_2^*)$ has a belt $B(M,{\bf x}_5^*,{\bf y}_2^*)$. And we have $$\{{\bf x}_5,{\bf y}_2\}=\{{\bf x}_5^*,{\bf y}_2^*\},$$ that is $(P+{\bf x}_5)\cap(P+{\bf y}_2)$ also has a belt $B(M,{\bf x}_5^*,{\bf y}_2^*)$. Then $(P+{\bf x}_5)\cap(P+{\bf y}_2)$ has two facets parallel to $F'$, which is impossible since $B(E,{\bf x}_5,{\bf y}_2)$ doesn't have any facet parallel to $F'$.

\smallskip\noindent
{\bf Subcase 2.2.} If $F'=F_5^*+{\bf x}_2^*$, by similar analysis to Subcase 2.1 above, this case is also impossible.

Thus, this lemma is proved.\hfill{$\Box$}


\medskip\noindent
{\bf Lemma 7.2.} {\it If $P$ is a fivefold translative tile in $\mathbb{E}^3$ with $\#B(E)=10$, then $\#B(M)\not=8$.}

\medskip\noindent
{\bf Proof.} For the contrary, suppose that $P$ is a fivefold translative tile with $\#B(E)=10$ and $\#B(M)=8$, then we have that $P+X$ is both a proper fivefold translative tiling relate to $E$ and a proper fivefold translative tiling relate to $M$ based on the results above. Since there are two kinds of projection of $P$ for $\#B(M)=8$, without loss of generality, we have the projection of the dihedral adjacent wheel at $M^*$ along $M$ shown in FIG \ref{P8_M1} and FIG \ref{P8_M2} where {\bf q} is a relative interior point on $M^*$ and we assume that $P+{\bf x}_2^*$ is the F-type translate at $M^*$. Notice that $F'\notin\{F_1^*+{\bf x}_2^*, F_5^*+{\bf x}_1^*,F_5^*+{\bf x}_3^*\}$ since $(P+{\bf x}_2^*)\cap(P+{\bf x}_1^*)\not= F_1^*+{\bf x}_2^*$ and $(P+{\bf x}_2^*)\cap(P+{\bf x}_3^*)\not= F_1^*+{\bf x}_2^*$. For convenience, we write $P$ whose projection along $M$ is shown in FIG \ref{P8_M1} as $P_1$ and write $P$ whose projection along $M$ is shown in FIG \ref{P8_M2} as $P_2$.

For $P=P_1$, we have the following discussions.

\smallskip\noindent
{\bf Case 1.} $F'\in\{F_1+{\bf x}_1, F_3+{\bf x}_4, F_2+{\bf x}_4\}$ and $F'=F_1+{\bf x}_1$. Observe the dihedral adjacent wheel at $M^*$, and we have $F'\in\{F_4^*+{\bf x}_3^*, F_3^*+{\bf x}_5^*\}$ since there need four translates containing $F'$ in their interior.

\begin{figure}[ht!]
  \centering
  \includegraphics[scale=0.65]{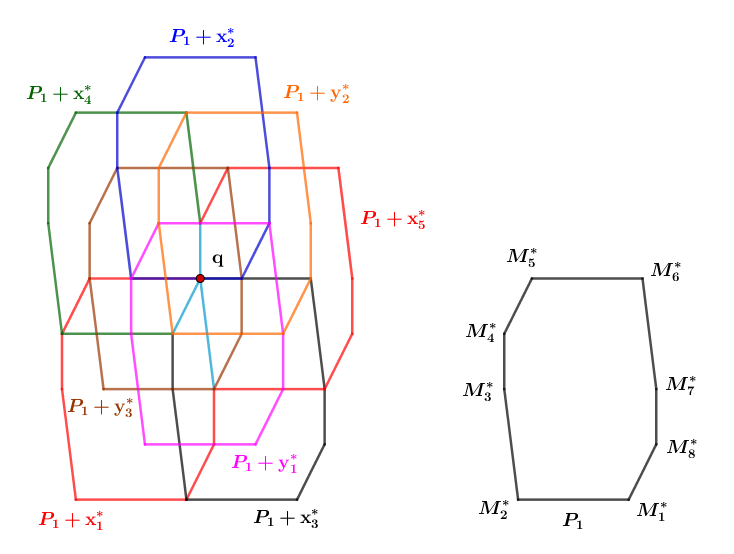}
  \caption{The dihedral adjacent wheel at $M^*$ in $P_1+X$ for $\#B(M)=8$}
 \label {P8_M1}
\end{figure}

\begin{figure}[h!]
\hspace{-1cm}
\subfigure[]
{\begin{minipage}{0.5\linewidth}
  \includegraphics[scale=0.64]{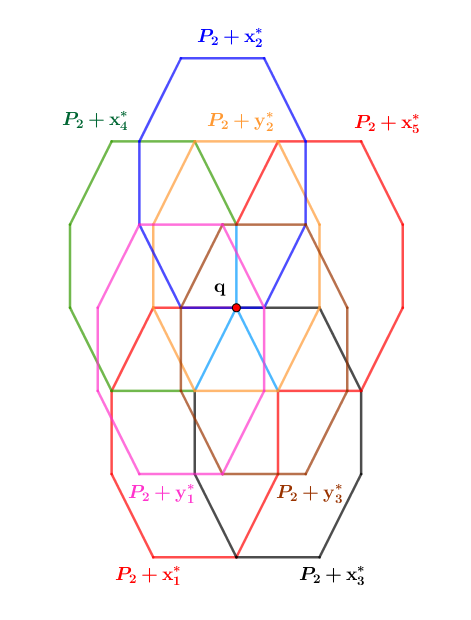}
  \end{minipage}}
\subfigure[]
{\begin{minipage}{0.4\linewidth}
  \includegraphics[scale=0.64]{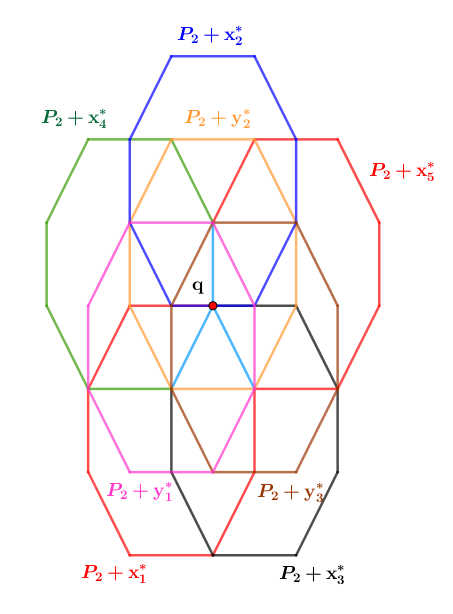}
  \end{minipage}}
  \subfigure[]
{\begin{minipage}{0.9\linewidth}
\vspace{1cm}
\includegraphics[scale=0.59]{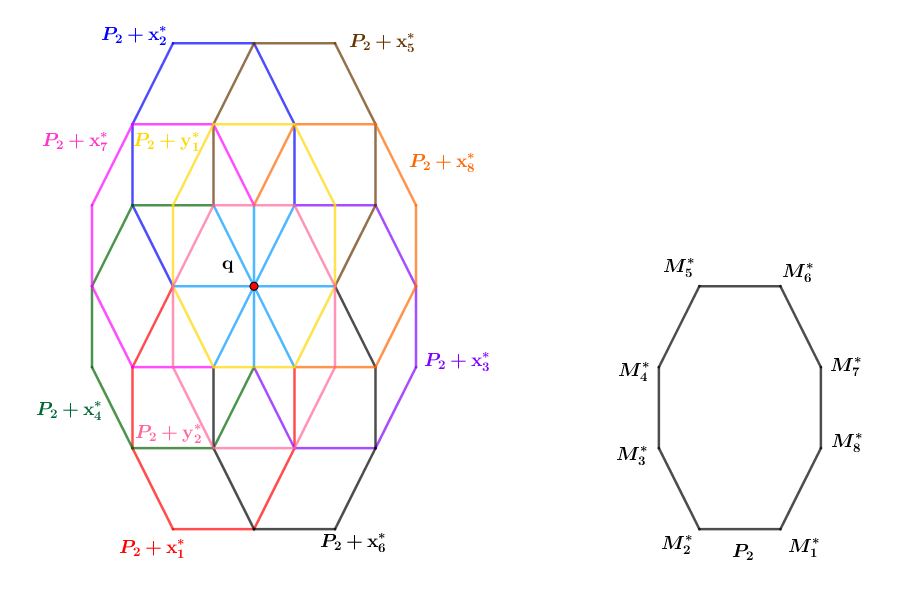}
  \end{minipage}}
 \caption{The dihedral adjacent wheel at $M^*$ in $P_2+X$ for $\#B(M)=8$}
 \label {P8_M2}
\end{figure}

\smallskip\noindent
{\bf Subcase 1.1.} $F'=F_4^*+{\bf x}_3^*$, we have $relint(F')\subset int(P)+{\bf x}$ for ${\bf x}\in\{{\bf x}_1^*,{\bf y}_1^*,{\bf y}_2^*,{\bf y}_3^*\}$, and $(P+{\bf x}_1^*)\cap(P+{\bf y}_2^*)$ and $(P+{\bf y}_1^*)\cap(P+{\bf y}_3^*)$ are both centrally symmetric about the center of $F'$. Then we have $$\{{\bf x}_3,{\bf y}_2\}=\{{\bf x}_1^*, {\bf y}_2^*\}$$ or $$\{{\bf x}_3,{\bf y}_2\}=\{{\bf y}_1^*, {\bf y}_3^*\}.$$ By Lemma 4.3, we can conclude that $(P+{\bf x}_1^*)\cap(P+{\bf y}_2^*)$ has a belt $B(M,{\bf x}_1^*,{\bf y}_2^*)$ and $(P+{\bf y}_1^*)\cap(P+{\bf y}_3^*)$ has a belt $B(M,{\bf y}_1^*,{\bf y}_3^*)$. In other words, $(P+{\bf x}_3)\cap(P+{\bf y}_2)$ has a belt $B(M,{\bf x}_3,{\bf y}_2)$, and it should have two facets parallel to $F'$, which is impossible since $B(E,{\bf x}_3,{\bf y}_2)$ doesn't have any facet parallel to $F'$.

\smallskip\noindent
{\bf Subcase 1.2.} $F'=F_3^*+{\bf x}_5^*$, we have $relint(F')\subset int(P)+{\bf x}$ for ${\bf x}\in\{{\bf x}_2^*,{\bf y}_1^*, {\bf y}_2^*, {\bf y}_3^*\}$, $(P+{\bf x}_2^*)\cap(P+{\bf y}_1^*)$ and $(P+{\bf y}_2^*)\cap(P+{\bf y}_3^*)$ are both centrally symmetric about the center of $F'$. Then we have $$\{{\bf x}_3, {\bf y}_2\}=\{{\bf x}_2^*,{\bf y}_1^*\}$$ or $$\{{\bf x}_3, {\bf y}_2\}=\{{\bf y}_2^*,{\bf y}_3^*\}.$$ Similar to Subcase 1.1 above, $(P+{\bf x}_3)\cap(P+{\bf y}_2)$ should have two facets parallel to $F'$, which is impossible.

\smallskip\noindent
{\bf Case 2.} $F'\in\{F_5+{\bf x}_2, F_4+{\bf x}_2\}$ and $F'=F_4+{\bf x}_2$. Observing the dihedral adjacent wheel at $M^*$, we have $F'=F_2^*+{\bf x}_5^*$ since there need three translates containing $F'$ in their interior, and $$\{{\bf x}_5,{\bf y}_2\}=\{{\bf x}_3^*, {\bf y}_3^*\}.$$ Since $(P+{\bf x}_3^*)\cap(P+{\bf y}_3^*)=(P+{\bf x}_5)\cap(P+{\bf y}_2)$ has two belts $B(M,{\bf x}_3^*,{\bf y}_3^*)$ and $B(E,{\bf x}_5,{\bf y}_2)$, it has two facets parallel to $F'$ and this is impossible since $B(E,{\bf x}_5,{\bf y}_2)$ doesn't have any facet parallel to $F'$.

For $P=P_2$, by the symmetry of translates at $M^*$ above, we can delete FIG \ref{P8_M2}(b). For (a), we can conclude that this is impossible by the similar discussion in $P=P_1$. For (c), we only need to consider $F'=F_5^*+{\bf x}_1^*$ by symmetry, similarly, we can also conclude that this is impossible.

Therefore we have finished the proof of Lemma 7.2. \hfill{$\Box$}

\medskip
According to Lemma 4.4 and Lemmas 7.1-7.2, we can obtain the consequence following.

\medskip\noindent
{\bf Lemma 7.3.} {\it If $P+X$ is a proper fivefold translative tiling relate to $E$ for $\#B(E)=10$, then $P$ must be a cylinder over a two-dimensional fivefold decagonal tile.}

\bigskip\noindent
{\Large\bf 7.2 $P+X$ for $\#B(E)=8$}

\bigskip\noindent
As introduced in Subsection 6.2, we still use the notations before. See in FIG \ref{P8_E1} and FIG \ref{P8_E2} for $\#B(E)=8$ and we write $P$ as $P_1$ and $P_2$ by their projections along $E$, respectively. Let {\bf p} be a proper point on $E'\in A(E)+X$ and let the dihedral adjacent wheel at {\bf p} consist of $P+{\bf x}_1, P+{\bf x}_2, P+{\bf x}_3, P+{\bf x}_4,P+{\bf x}_5$ in clock order where ${\bf p}\in relint(F_1)+{\bf x}_4$. And let $P+{\bf y}_1, P+{\bf y}_2, P+{\bf y}_3$ be the I-type translates at {\bf p}. Since $P+X$ is a proper fivefold translative tiling relate to $E$, we also call $P+{\bf x}_i$ and $P+{\bf y}_j$ the translates at $E'$ for convenience. Without loss of generality, let $E'=E_8$ and ${\bf x}_1={\bf o}$.

By observing the translates at $E_8$, we can get that $F'\in\{F_3+{\bf x}_2, F_2+{\bf x}_2, F_1+{\bf x}_4, F_4+{\bf x}_5\}$ is a translate of $F$. Let $M^*$ be the edge of $F'$ that is a translate of $M$. Without loss of generality, we can only analyze FIG \ref{P8_E1} and FIG \ref{P8_E2}(a)(b) enough and other cases can be obtained by the similar method. Then we have the following results.

\medskip\noindent
{\bf Lemma 7.4.} {\it If $P$ is a fivefold translative tile in $\mathbb{E}^3$ with $\#B(E)=8$, then $\#B(M)\not=8$.}

\medskip\noindent
{\bf Proof.} For the contrary, suppose that $P$ is a fivefold translative tile with $\#B(E)=8$ and $\#B(M)=8$, then we have that $P+X$ is both a proper fivefold translative tiling relate to $E$ and a proper fivefold translative tiling relate to $M$ based on the results above. We consider $P=P_1$ and $P=P_2$, respectively.

\smallskip\noindent
{\bf Case 1.} $P=P_1$. Similar to the analysis in Lemma 7.2, we also have the following two subcases.

\smallskip\noindent
{\bf Subcase 1.1.} $F'\in\{F_3+{\bf x}_2, F_4+{\bf x}_5,F_2+{\bf x}_2\}$. Similar to the discussion in Lemma 7.2, we can conclude that this case is impossible.

\smallskip\noindent
{\bf Subcase 1.2.} $F'=F_1+{\bf x}_4$. It is easy to see that $6\le\#F'\le10$ since $P_1+X$ is a non-proper fivefold translative tiling relate to $M$ if $\#F'=4$ which contradicts our discussion. Since $P+X$ is also a proper fivefold translative tiling relate to $M$ with $\#B(M)=8$, without loss of generality, we can assume that $$M^*\subset(F_1+{\bf x}_4)\cap(F_5+{\bf x}_3)$$ or $$M^*\subset(F_1+{\bf x}_4)\cap(F_5+{\bf x}_5).$$

\smallskip\noindent
{\bf Subcase 1.2.1.} If $M^*\subset(F_1+{\bf x}_4)\cap(F_5+{\bf x}_5)$, for convenience, write $F^*=(F_1+{\bf x}_4)\cap(F_5+{\bf x}_5)$, then we have $relint(F^*)\subset int(P)+{\bf x}$ for $\{{\bf x}_2, {\bf y}_3, {\bf y}_1, {\bf y}_2\}$, and $(P_1+{\bf x}_2)\cap(P_1+{\bf y}_3)$ and $(P_1+{\bf y}_1)\cap(P+{\bf y}_2)$ are both centrally symmetric about the center of $F^*$. Now let us observe the translates at $M^*$. If the projection of $P_1$ along $M^*$ is the projection of $P_1$ along $E$ by a suitable affine transformation, then we have $$\{{\bf x}_2, {\bf y}_3\}=\{{\bf x}_5^*,{\bf y}_3^*\}$$ or $$\{{\bf x}_2, {\bf y}_3\}=\{{\bf y}_1^*,{\bf y}_2^*\},$$ and $(P_1+{\bf x}_2)\cap(P_1+{\bf y}_3)$ should have facets parallel to $F^*$ which is impossible. If the projection of $P_1$ along $M^*$ is the projection of $P_2$ along $E$ by a suitable affine transformation, we can only consider FIG \ref{P8_M2}(a)(b) by the symmetry, similarly above and Lemma 7.2, we can also conclude that this case is impossible.

\smallskip\noindent
{\bf Subcase 1.2.2.} If $M^*\subset(F_1+{\bf x}_4)\cap(F_5+{\bf x}_3)$, for convenience, write $F^*=(F_1+{\bf x}_4)\cap(F_5+{\bf x}_3)$, then we have $relint(F^*)\subset int(P_1)+{\bf x}$ for ${\bf x}\in\{{\bf x}_1, {\bf y}_1, {\bf y}_3\}$, and $(P_1+{\bf x}_1)\cap(P+{\bf y}_1)$ and $P_1+{\bf y}_3$ are both centrally symmetric about the center of $F^*$. Observing the translates at $M^*$ in FIG \ref{P8_M1} and FIG \ref{P8_M2}, there is only $F^*=(F_1^*+{\bf x}_2^*)\cap(F_5^*+{\bf x}_1^*)$ in FIG \ref{P8_M1} for the projection of $P_1$ along $M^*$ is the projection of $P_1$ along $E$ by a suitable affine transformation satisfying this case. Then $(P_1+{\bf x}_1)\cap(P_1+{\bf y}_1)=(P_1+{\bf x}_4^*)\cap(P_1+{\bf y}_1^*)$ has two belts $B(E,{\bf x}_1,{\bf y}_1)$ and $B(M,{\bf x}_1,{\bf y}_1)$  and it contains one facet $F^\bullet$ having the translates of $E$ and $M$ as its edges. By observation, we have $F^\bullet$ is a translate of $(F_1+{\bf x}_4)\cap(F_5+{\bf x}_5)$.

Apparently, if $(F_1+{\bf x}_4)\cap(F_5+{\bf x}_5)$ doesn't have the translates of $M$ as its edges, we can have that this case doesn't exist.

If $(F_1+{\bf x}_4)\cap(F_5+{\bf x}_5)$ does have the translate $M^{**}$ of $M$ as its edge, as shown in Subcase 1.2.1 above, this case is also impossible.

\smallskip\noindent
{\bf Case 2.} $P=P_2$. According to the consequence in $P=P_1$ above, we have that the projection of $P_2$ along $M$ is only a suitable affine transformation of the projection of $P_2$.

\smallskip\noindent
{\bf Subcase 2.1.} $F'\in\{F_4+{\bf x}_5, F_3+{\bf x}_2, F_2+{\bf x}_2\}$. Similar to the discussion in Lemma 7.2, we can conclude that this case is impossible.

\smallskip\noindent
{\bf Subcase 2.2.} $F'=F_1+{\bf x}_4$. If the projection of the dihedral adjacent wheel at $E$ is shown in FIG \ref{P8_E2}(b), the analysis is similar to that in $P=P_1$, then we can conclude that it's impossible. Thus we can only consider the translates shown in FIG \ref{P8_E2}(a), and without loss of generality, we assume that the translates at $M^*$ are shown in FIG \ref{P8_M2}. Similar to Subcase 1.2 above, we also have the following two subcases.

\smallskip\noindent
{\bf Subcase 2.2.1.} $M^*\subset(F_1+{\bf x}_4)\cap(F_5+{\bf x}_5)$, for convenience, write $F^*=(F_1+{\bf x}_4)\cap(F_5+{\bf x}_5)$ then we have that $relint(F^*)\subset int(P_2)+{\bf x}$ for ${\bf x}\in\{{\bf x}_2, {\bf y}_3, {\bf y}_1, {\bf y}_2\}$, and $(P_2+{\bf x}_2)\cap(P+{\bf y}_3)$ and $(P_2+{\bf y}_1)\cap(P_2+{\bf y}_2)$ are centrally symmetric about the center of $F^*$. Moreover, $(P_2+{\bf x}_2)\cap(P_2+{\bf y}_3)$ has one belt determined by the translate of $E$ but no belt determined by the translate of $M$. Similar to Case 2 in Lemma 7.2, this case is impossible for both $F^*=(F_1^*+{\bf x}_2^*)\cap(F_5^*+{\bf x}_3^*)$ and $F^*=(F_1^*+{\bf x}_2^*)\cap(F_5^*+{\bf x}_1^*)$ by observing the structures of $(P_2+{\bf x}_5^*)\cap(P_2+{\bf y}_1^*)$, $(P_2+{\bf y}_2^*)\cap(P_2+{\bf y}_3^*)$, $(P_2+{\bf x}_4^*)\cap(P_2+{\bf y}_3^*)$ and $(P_2+{\bf y}_2^*)\cap(P_2+{\bf y}_1^*)$, respectively.

\begin{figure}[h!]
\includegraphics[scale=0.6]{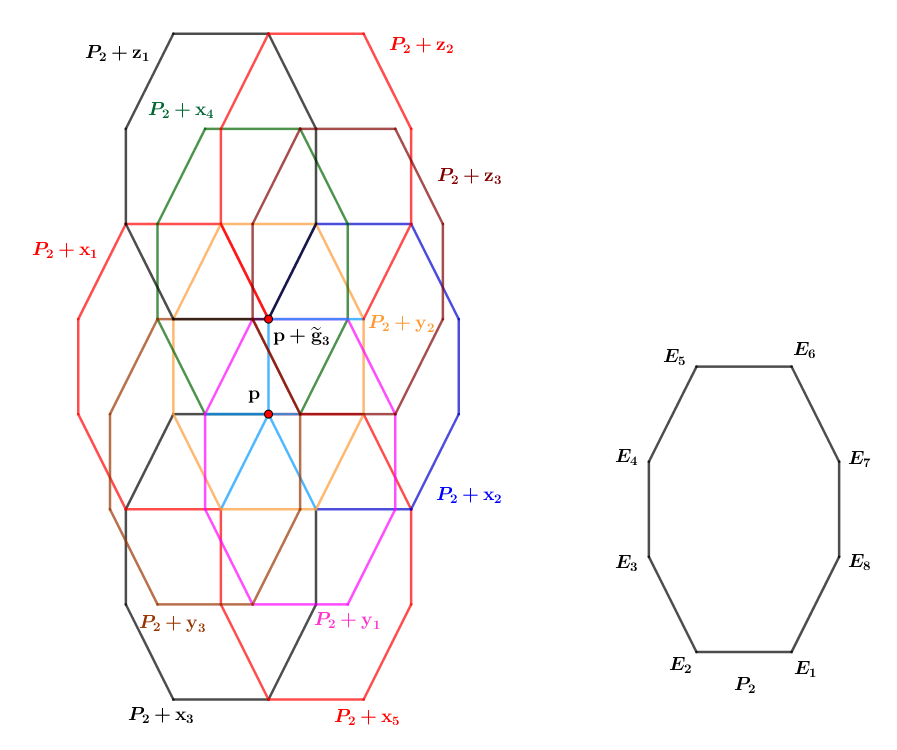}
\caption{The translates at $E+\widetilde{\bf g}_3$}\label{P8_2E+g33}
\end{figure}

\smallskip\noindent
{\bf Subcase 2.2.2.} $M^*\subset(F_1+{\bf x}_4)\cap(F_5+{\bf x}_3)$, for convenience, write $F^*=(F_1+{\bf x}_4)\cap(F_5+{\bf x}_3)$ then we have that $relint(F^*)\subset int(P_2)+{\bf x}$ for ${\bf x}\in\{{\bf x}_1, {\bf y}_1, {\bf y}_2, {\bf y}_3\}$, and $(P_2+{\bf x}_1)\cap(P+{\bf y}_1)$ and $(P_2+{\bf y}_2)\cap(P_2+{\bf y}_3)$ are centrally symmetric about the center of $F^*$. Moreover, $(P_2+{\bf x}_1)\cap(P_2+{\bf y}_1)$ and $(P_2+{\bf y}_2)\cap(P_2+{\bf y}_3)$ have belts $B(E,{\bf x}_1,{\bf y}_1)$ and $B(E,{\bf y}_2,{\bf y}_3)$, respectively. If $F^*=(F_1^*+{\bf x}_2^*)\cap(F_5^*+{\bf x}_3^*)$, it is easy to see this case is impossible by the analysis above. If $F^*=(F_1^*+{\bf x}_2^*)\cap(F_5^*+{\bf x}_1^*)$, by observing the translates at $M^*$ and relying on Subcase 2.2.1 above, we have that $$\{{\bf x}_1, {\bf y}_1\}=\{{\bf x}_4^*,{\bf y}_3^*\}$$ or $$\{{\bf x}_1, {\bf y}_1\}=\{{\bf y}_1^*,{\bf y}_2^*\},$$ and $(F_1+{\bf x}_4)\cap(F_5+{\bf x}_5)$ doesn't have the translate of $M^*$ as its edge otherwise we can return Subcase 2.2.1 above. If $B(E,{\bf x}_1,{\bf y}_1)$ or $B(E,{\bf y}_2,{\bf y}_3)$ doesn't have two facets with the translate of $M^*$ as its edge, then we can conclude that this case is impossible. Then we consider that $B(E,{\bf x}_1,{\bf y}_1)$ and $B(E,{\bf y}_2,{\bf y}_3)$ both have facets with the translate of $M^*$ as their edge. Now let us observe the translates at $E_8+\widetilde{\bf g}_3$, as shown in FIG \ref{P8_2E+g33}.

We can see that $(F_1+{\bf z}_2)\cap(F_5+{\bf y}_1)$ is contained in a translate of one facet in $B(E,{\bf y}_2,{\bf y}_3)$ and it also has the translate $M^{**}$ of $M$ as its edge. For convenience, let $F^{**}=(F_1+{\bf z}_2)\cap(F_5+{\bf y}_1)$, and $relint(F^{**})\subset int(P_2)+{\bf x}$ for ${\bf x}\in\{{\bf x}_2,{\bf x}_4, {\bf y}_2, {\bf z}_3\}$. It is easy to see that both $(P_2+{\bf x}_2)\cap(P_2+{\bf x}_4)$ and $(P_2+{\bf y}_2)\cap(P_2+{\bf z}_3)$ are centrally symmetric about the center of $F^{**}$. And $(P_2+{\bf x}_2)\cap(P_2+{\bf x}_4)$ has a belt $B(E,{\bf x}_2,{\bf x}_4)$ which doesn't have facets with the translate of $M^{**}$ as its edge since $(F_1+{\bf x}_4)\cap(F_5+{\bf x}_5)$ doesn't have the translate of $M^*$ as its edge. Then by considering the translates at $M^{**}$, we can conclude that this case is also impossible.

Therefore we have finished the proof of Lemma 7.4. \hfill{$\Box$}

\medskip
According to Lemma 4.5 and Lemma 7.4, we can obtain the consequence following.

\medskip\noindent
{\bf Lemma 7.5.} {\it If $P+X$ is a proper fivefold translative tiling relate to $E$ for $\#B(E)=8$, then $P$ must be a cylinder over a two-dimensional fivefold octagonal tile.}

\vspace{1cm}\noindent
{\LARGE\bf 8. Proof of Theorem 1.1}

\medskip\noindent
{\bf Proof of Theorem 1.1.} Let $P+X$ be a fivefold translative tiling of $\mathbb{E}^3$. By Lemmas 2.3-2.4, we have known that every belt of $P$ has four, six, eight or ten facets. On the one hand, if each belt of $P$ only contains four or six facets, then $P$ must be a parallelohedron by Lemmas 2.1-2.2 and Lemma 2.4. One the other hand, if there is one belt of $P$ containing eight or ten facets, then $P$ must be a cylinder over a two-dimensional fivefold translative tile by Theorem 5.1, Lemmas 6.1-6.3, Lemma 7.3 and Lemma 7.5.

Thus, Theorem 1.1 is proved.   \hfill{$\Box$}

\vspace{1cm}\noindent
{\LARGE\bf 9. A non-trivial multiple lattice tile}

\medskip\noindent
For convenience, we call a multiple tile in $\mathbb{E}^3$ a {\it non-trivial multiple tile}, if it is neither a parallelohedron nor a cylinder. In this section, our purpose is to find a non-trivial multiple tile with multiplicity as small as possible. The idea of our construction is based on the following result.

\medskip\noindent
{\bf Lemma 9.1 (Gravin, Robins and Shiryaev \cite{grs}; Lev and Liu \cite{lev-liu}).} {\it In $\mathbb{E}^n$, a centrally symmetric convex polytope $P$, with centrally symmetric facets, if all the vertices of $P$ belongs to a lattice $\Lambda$, then $P$ admits a multiple lattice tiling by translating along with $\Lambda$.}

\medskip
As shown in Section 3, let $P$ be a three-dimensional zonotope with the generator set $\mathbb{S}_P=\{S_0,S_1,...,S_w\}$, where $\mathbb{S}_P$ is a set of pairwise linear independent segments centered at {\bf o}. And let ${\bf v}_{S_i}=(v_1^i, v_2^i, v_3^i)$ be an endpoint of $S_i$ for each $0\le i\le w$, where all $v_j^i$ are the coordinate components of ${\bf v}_{S_i}$.

By Lemma 9.1, we have the following example for a non-trivial multiple lattice tile.

\medskip\noindent
{\bf Example 9.1.} Let $P$ be a zonotope in $\mathbb{E}^3$ with $\{S_0, S_1, S_2, S_3, S_4\}$ as its generator set and let ${\bf v}_{S_0}=(0,0,\frac{1}{2})$, ${\bf v}_{S_1}=(0,\frac{1}{2},0)$, ${\bf v}_{S_2}=(\frac{1}{2},0,0)$, ${\bf v}_{S_3}=(\frac{1}{2},\frac{1}{2},\frac{1}{2})$ and ${\bf v}_{S_4}=(\frac{1}{2},1,\frac{1}{2})$, as shown in Figure \ref{k=10}.

\begin{figure}[h!]
\includegraphics[scale=0.6]{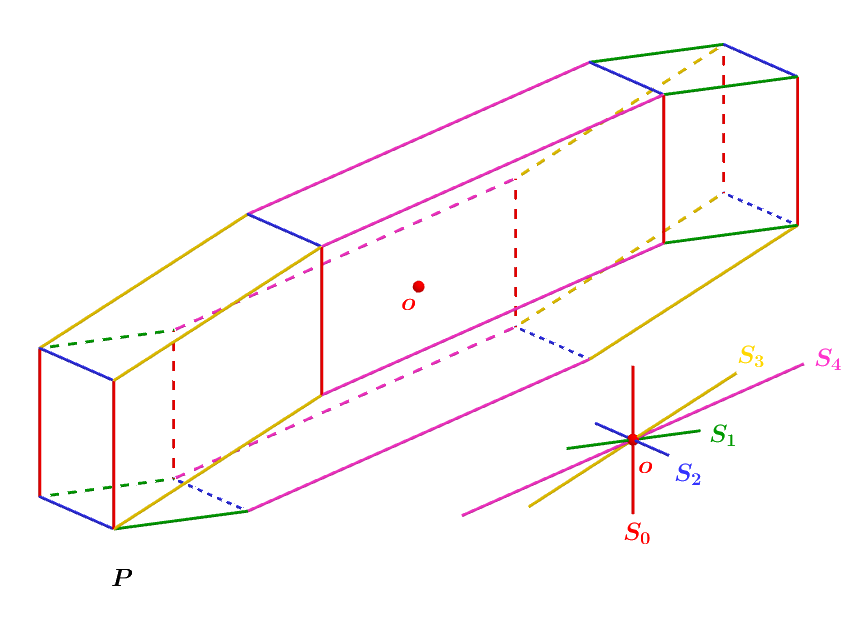}
\caption{A non-trivial multiple lattice tile $P$}\label{k=10}
\end{figure}

It is easy to get that $P$ is a multiple lattice tile in $\mathbb{E}^3$ and the volume of $P$ is 10, then we have that $P+\mathbb{Z}^3$ is a ten-fold lattice tiling of $\mathbb{E}^3$. Thus, $P$ is a non-trivial multiple tile in $\mathbb{E}^3$ with multiplicity at most 10.

\vspace{0.6cm}\noindent
{\bf Acknowledgements.} This work is supported by the National Natural Science Foundation of China (NSFC12226006, NSFC11921001), the National Key Research and Development Program of China (2018YFA0704701), 973 Program 2013CB834201, and a special fund for postgraduate education of Tianjin University (B1-2021-001).}

\vspace{0.6cm}
\noindent
Mei Han, Center for Applied Mathematics, Tianjin University, Tianjin 300072, China.

\noindent
Email: ludy$_-$han@163.com

\medskip\noindent
Kirati Sriamorn, Department of Mathematics and Computer Science, Chulalongkorn University, Thailand.

\noindent
Email: Kirati.S@chula.th

\medskip\noindent
Qi Yang, College of Mathematics and Statistics, Chongqing University, China.

\noindent
Email: yangqi7@cqu.edu.cn

\medskip\noindent
Chuanming Zong, Center for Applied Mathematics, Tianjin University, Tianjin 300072, China.

\noindent
Email: cmzong@tju.edu.cn

\end{document}